\documentclass{amsart}
\usepackage{amsmath,amsfonts,amscd,latexsym,amssymb,epsfig,float,
afterpage}

\newtheorem{theorem}{Theorem}[section]
\newtheorem{lemma}[theorem]{Lemma}
\newtheorem{proposition}[theorem]{Proposition}

\theoremstyle{definition}
\newtheorem{definition}[theorem]{Definition}

\newtheorem{corollary}[theorem]{Corollary}
\newtheorem{remark}[theorem]{Remark}
\newtheorem{fact}[theorem]{Fact}
\newtheorem{proposition+definition}[theorem]{Proposition \& Definition}
\newtheorem{warning}[theorem]{Warning}

\theoremstyle{remark}

\numberwithin{equation}{section}

\newcommand{\cpxm}{{\mathbb{C}}}
\newcommand{\ratm}{{\mathbb{Q}}}
\newcommand{\itgm}{{\mathbb{Z}}}

\newcommand{\natm}{{\mathbb{N}}}
\newcommand{\reellm}{{\mathbb{R}}}
\newcommand{\korpm}{{\mathbb{K}}}
\newcommand{\prom}{{\mathbb{P}}}

\newcommand{\maxid}{{\mathfrak{m}}}

\newcommand{\im}{\operatorname{im}}
\newcommand{\Hom}{\operatorname{Hom}}

\newcommand{\Aut}{\operatorname{Aut}}
\newcommand{\Res}{\operatorname{Res}}

\newcommand{\Extmhs}{\operatorname{Ext_{{\scriptscriptstyle MHS}}}}

\newcommand{\lcm}{\operatorname{lcm}}

\newcommand{\topa}[2]{\genfrac{}{}{0pt}{}{#1}{#2}}

\newcommand{\bigsearrow}{\text{\Large $\searrow$}}
\newcommand{\bigswarrow}{\text{\Large $\swarrow$}}

\begin{document}
\renewcommand{\labelenumi}{(\roman{enumi})}

\title{New Period Mappings for Plane Curve Singularities}

\author{Rainer H.~Kaenders}

\address{Department of Mathematics, Utrecht University,
Postbus 80.010, 3508 TA Utrecht, The Netherlands} 
\email{kaenders@math.ruu.nl}

\thanks{This work was partly supported by grant ERBCHBICT930403 (HCM) from the
Eurpean Community and The Netherlands Organisation for Scientific
Research (NWO)}



\begin{abstract}
For plane curve singularities we construct a mixed Hodge structure (MHS)
over $\itgm$ on the fundamental group of the Milnor fiber. The
concept {\it nearby fundamental group} is introduced and we
develop a theory of iterated integrals along elements of this group.
We give an example for which the MHS on the nearby fundamental group 
detects a modulus which is invisible for the MHS on the vanishing
cohomology. 
\end{abstract}

\maketitle

 
\setcounter{tocdepth}{1}

\tableofcontents

\section*{Introduction}

The building block for many invariants of plane curve singularities 
like the monodromy, the intersection form, the Seifert form, 
Coxeter-Dynkin diagrams, the mixed Hodge structure on the vanishing cohomology 
and the spectral pairs is the first integral homology group of the Milnor 
fiber or -- which is the same -- the abelianized fundamental group of the 
Milnor fiber. 

The restriction to the abelianized version of this fundamental group 
causes some serious limitations. E.~g.~L{\^e} D{\~u}ng Tr{\'a}ng 
\cite{Le-algebraische-Knoten}
proved that for irreducible plane curve singularities the monodromy 
(on homology) has finite order. However A'Campo 
\cite{A'Campo-infinite-order} was the first to discover
that the monodromy on the fundamental group of the Milnor fiber 
(where the one boundary component is collapsed to one point, the basepoint)
has infinite order, whenever the singularity has at least two Puiseux pairs.
This shows that the monodromy on the fundamental group is richer than 
on homology. 

For the study of analytical invariants of plane curve 
singularities, the first (co)ho\-mo\-lo\-gy group is of particular interest
because of Hodge theory. 

Deligne \cite{HodgeII} introduced in 1970/71 the notion of mixed Hodge 
structure (MHS)
and showed that the cohomology of any complex algebraic variety carries
such a MHS in a natural way. Schmid \cite{Schmid} and Steenbrink 
\cite{Steenbrink-limits} determined a limit-MHS of a 
variation of Hodge structures coming from families of compact algebraic
manifolds. This provided the background for the 
definition of the MHS on the vanishing cohomology of hypersurface singularities 
by Steenbrink \cite{Steenbrink-Oslo}. In the case of plane curve singularities 
this is a MHS on the first cohomology group of the Milnor fiber. 

However, due to the fact that for irreducible plane curve singularities the 
monodromy on homology is finite, families of mutually 
right-left -- but not necessarily right-equivalent plane curve singularities
all have isomorphic MHSs on the vanishing cohomology. This holds, for example,
for the family
\begin{equation} \label{Liszt}
f_\lambda(x,y)=\left(y^2-x^3\right)^2-\left(4\lambda x^5 y+
\lambda^2 x^7\right),\quad \lambda \neq 0,
\end{equation}
which is the singularity with Puiseux expansion 
\[
y= x^{\frac{3}{2}}+\sqrt{\lambda}\, x^{\frac{7}{4}}.
\]

In 1987 Hain \cite{Hain-de-Rham-homotopy}, \cite{Hain-the-geometry}
used Chen's iterated integrals to generalize Deligne's construction of a
MHS on cohomology to a construction of a MHS on the homotopy Lie algebra 
for any complex algebraic variety with base point. In particular,
this object contains a MHS which is called {\it the MHS on the fundamental 
group}. Already before, in 1978, Morgan \cite{Morgan} had used Sullivan's 
minimal models to put a MHS on the homotopy Lie algebra of {\it smooth} complex 
algebraic varieties.

Hain also showed, in \cite{Hain-de-Rham-homotopy-II}, that 
the local system of homotopy groups associated with a topologically
trivial family of smooth pointed varieties underlies a good
variation of $\ratm$-MHSs, which implies by the earlier work of Steenbrink 
and Zucker \cite{Steenbrink-Zucker} that there is a limit $\ratm$-MHS.

This provides the background of the problem\footnote{This problem was 
proposed to us by J.~H.~M.~Steenbrink.} we study in this paper:
{\sf Is there a generalization of the MHS on the vanishing (co)homology to 
something like `a MHS on the fundamental group of the general Milnor fiber'
in the case of plane curve singularities?} 

In this article we show that the construction
of such a generalization is possible and we will call it {\it the mixed
Hodge structure on the nearby fundamental group of a plane
curve singularity}. 
We develop a theory that allows us to take iterated integrals 
along certain paths in the central fiber. These paths mimic paths 
in the regular fiber. The given approach seems to be even new for the 
integration of elements in the vanishing cohomology. Due to this
integration theory, the constructed MHS is defined over the integers. 
For instance in the
family $f_\lambda$ defined on page \pageref{Liszt} this $\itgm$-MHS detects 
the modulus. The above formulated
problem is the guideline of the paper. 

We consider first the situation of a degenerating 
family $h:(Z,D^+)\rightarrow(\Delta,0)$ of compact Riemann surfaces with
$r$ boundary components over the disk $\Delta\subset \cpxm$, where the
singular fiber, the fiber over $0$, is a divisor with normal crossings
$D=D^+\cup D_0\cup\cdots\cup D_{r-1}$. 
Its components are compact closed Riemann surfaces, except 
for the $r$ components $D_0,\ldots, D_{r-1}$, which are disks. 
Such a degenerating family can 
be constructed from a plane curve singularity $f$ by a process
that is called `semistable reduction' (see \ref{Semistable}).
Each of the disks $D_0,\ldots, D_{r-1}$ corresponds to a branch of $f$. 

\begin{figure}[ht]
\centering\epsfig{file=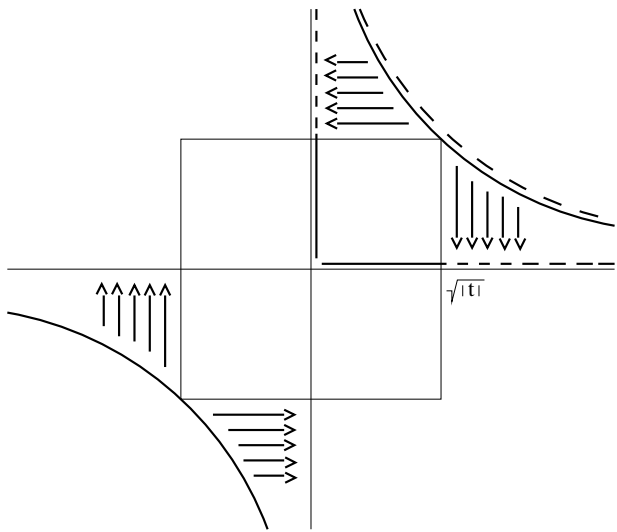,height=7cm}
\caption{Identification of a path in the regular fiber of the map 
$(x,y)\mapsto x\cdot y=t$ with a path in the singular fiber.} \label{Madonna}
\end{figure}

In this situation, we want to mimic a path in a regular fiber by a path 
in the central fiber. Let us give a finite version of the ideas that we 
develop in an infinitesimal way in Chapter \ref{Rietgans}. 
Intuitively it is clear that the difficulties for this mimicry arise 
at the double points of the divisor with normal crossings $D$, which
are all given locally by an equation: $x\cdot y=t$. If $t\neq 0$, 
this equation gives within the central fiber (locally) a 1-1 relation 
between points off the
double point in one component and points off the double point in the other.
This is obviously not true for $t=0$. For $t\neq 0$, we consider 
then paths $\gamma$ that approach the double point in one component 
and leave it in the other component in such a way that locally for all
$\varepsilon >0$ holds:
\[
x\circ\gamma(\tau_0-\varepsilon) \cdot y\circ \gamma(\tau_0+\varepsilon)=t,
\]
where $\gamma(\tau_0)$ is the double point. By the identification that is
sketched in Figure \ref{Madonna}, these paths can be thought of as being in the fiber over $t$. 

In this picture one can also observe that if $|t|$ is big, the path in 
the central fiber has to make a big detour compared to the corresponding
path in the regular fiber. Therefore, the speed with which the path in the
central fiber approaches and leaves the double point is a measure for 
the distance between the central fiber and the regular fiber in which the 
`regular path' lives. We give an infinitesimal version of this idea in
\ref{Verspannung}. Given a non zero tangent vector $\vec{v}$ in the
tangent plane $T_0\Delta$ of the base space $\Delta$ and a non zero 
tangent vector
$\vec{w}$ in the tangent plane $T_{p_0}D_0$ at the disk $D_0$,
we introduce the notion of {\it a path over $\vec{v}$ based at $\vec{w}$}.
Moreover, we define a {\it homotopy over $\vec{v}$ based at $\vec{w}$}
or shorter, a {\it nearby homotopy}, in such a way that the set of all
{\it nearby homotopy classes} becomes a group that we denote by:
\[
\pi_1(Z_{\vec{v}},\vec{w}).
\]
The $\{\pi_1(Z_{\vec{v}},\vec{w})\}_{\vec{v}\in(T_{0}\Delta)^*,
\vec{w}\in(T_{p_0}D_0)^*}$ form a local system on $(T_{0}\Delta)^*
\times(T_{p_0}D_0)^*$, which is isomorphic to a corresponding local 
system that comes from the fundamental groups on the regular fibers
(see \ref{Dr. Yoghurt} on page \pageref{Dr. Yoghurt}).

Denote by $\vec{J}=J_{\vec{v},\vec{w}}$ the augmentation ideal of the 
group ring $\itgm\pi_1(Z_{\vec{v}},\vec{w})$. In the Chapters \ref{Quark} 
and \ref{Bl"omscher}, we will define for all $s\ge 0$ a $\itgm$-MHS on
\begin{equation} \label{Braunb"ar}
\itgm\pi_1(Z_{\vec{v}},\vec{w})\big/ J_{\vec{v},\vec{w}}^{s+1}.
\end{equation}
This will be done by putting a $\itgm$-MHS on
$({\vec{J}}\big/{\vec{J}^{s+1}})^*$, the dual of 
${\vec{J}}\big/{\vec{J}^{s+1}}$ $\subset$
$\itgm\pi_1(Z_{\vec{v}},\vec{w})\big/ J_{\vec{v},\vec{w}}^{s+1}$. 

The main ingredient for the definition of these $\itgm$-MHSs is a differential
graded algebra (dga) $A^\bullet$, which is made up from differential
forms on the central fiber but whose cohomology is isomorphic to the cohomology
of the regular fiber. It was suggested to us by Hain\footnote{The dga 
$A^\bullet$ is in principle a modified version of the complex 
$P_\cpxm^\bullet(X^*)[\theta]$ of \cite{Hain-de-Rham-homotopy-II}
for this special situation with non compact fibers.}
to use a dga like $A^\bullet$ and he proposed the problem of finding
a way to integrate elements of this dga along some kind of 
paths in the central fiber. We discovered
that the right notion of {\it path}
for this purpose turned out to be the notion of a path over $\vec{v}$
(based at $\vec{w}$). We will show that closed elements in $A^1$ can be
integrated along paths in $Z_0$ over $\vec{v}$ in a natural way. 
This yields an integral
structure on $H^1(A^\bullet)$, which is by this integration dual to
${\vec{J}}\big/ {\vec{J}^{\;2}}$. On $A^\bullet$ we have filtrations $W_\bullet$
and $F^\bullet$ such that the induced filtrations together with the integral
structure define a $\itgm$-MHS on $H^1(A^\bullet)$ (see Theorem
\ref{Rothalsgans} on page \pageref{Rothalsgans}).

Since $A^\bullet$ has an {\it algebra} structure, we can also define {\it
iterated}
line integrals along paths over $\vec{v}$. We shall do this in Chapter
\ref{Bl"omscher}. This will give a way (see \ref{Novalis}) to describe 
$\Hom_\itgm({\vec{J}}\big/ {\vec{J}^{s+1}},\cpxm)$ in terms of the
dga $A^\bullet$, similar to Chen's classical $\pi_1$-De-Rham theorem 
(cf.~(4.4) in \cite{Hain-the-geometry}).

These iterated integrals will lead to the definition of a $\itgm$-MHS on
$({\vec{J}}\big/ {\vec{J}^{\,s+1}})^*$, and hence to a $\itgm$-MHS on 
\eqref{Braunb"ar} (see \ref{Blaumeise}).
And furthermore, we will see that the variation of $\itgm$-MHS for all 
$s\ge 1$
\[
\left\{J_{\vec{v},\vec{w}}\big/J_{\vec{v},\vec{w}}^{s+1}\right\}_{\vec{v}\in
(T_{0}\Delta)^*,\;\vec{w}\in(T_{p_0}D_0)^*}
\]
is a nilpotent orbit of MHS (see \ref{Maibaum} and \ref{Blindfisch}).

Let $\pi_1(Z_0,p_0)$ be the fundamental group of the central fiber and let 
$\itgm\pi_1(Z_0,p_0)$ be its group ring with augmentation ideal $J_0$. 
We define a MHS on the fundamental group of the central fiber.
This MHS is more or less a special case of the general construction
of a MHS on the fundamental group of a complex algebraic varietiy as 
introduced by Hain \cite{Hain-de-Rham-homotopy}. 
We show that the obvious group homomorphism 
\[
c:\pi_1(Z_{\vec{v}},\vec{w})\rightarrow \pi_1(Z_{0},p_0)
\]
induces an inclusion of MHSs:
\[
c^*:\big(J_0\big/J_0^{s+1}\big)^*\rightarrow 
\big(\vec{J}\big/\vec{J}^{s+1}\big)^*.
\]

In Chapter \ref{Kranenburg} we show how all these considerations can
be applied to study plane curve singularities
\[
f=f_0\cdots f_{r-1}:(\cpxm^2,0)\rightarrow (\cpxm,0).
\]
Let $t$ be a coordinate on $(\cpxm,0)$.
The tangent vector $\frac{\partial}{\partial t}$ in $T_0\cpxm$ yields
a finite number of tangent vectors in the base space of the 
`semistable reduction'. 
If $f_0$ is a distinguished branch of $f$ and $D_0$ the disk 
corresponding to $f_0$, then $\frac{\partial}{\partial t}$ 
defines a finite number
of tangent vectors in $T_{p_0}D_0$ in a natural way, which we call 
{\it the monstrance of $(f,f_0)$}. 
We use this to define 
{\it the MHS on the nearby fundamental group
of a plane curve singularity $f$}.
This MHS depends on the tangent vector $\frac{\partial}{\partial t}$,
which is preserved under right-equivalence. It might be 
interesting to note in this 
context that if a right-left-equivalence between two plane
curve singularities preserves $\frac{\partial}{\partial t}$ that
then the singularities are right-quivalent. This is a consequence of
a result of J.~Brian\c{c}on and H.~Skoda (see \ref{Geburtszange}).  

Finally we consider the 
{\it MHS on the nearby fundamental group} of the example
\eqref{Liszt} on page \pageref{Liszt} and show that the extension
of $H_1$ by the part of weight $-3$ of $\vec{J}\big/\vec{J}^4$
detects the modulus that is hidden for the vanishing cohomology.
Here we focus on the effects due to the
infinite order of the monodromy. 

In general we expect also interesting information from other
parts of this MHS on the nearby fundamental group. Even the extension
of $H_1$ by $H_1\otimes H_1$ given by $\vec{J}\big/\vec{J}^{\;3}$ might 
contain interesting information. 
Also the decomposition of the monodromy into a unipotent
and a semisimple part has still to be studied.
The ideas described here for plane curve singularities have a counterpart
for degenerating families of compact Riemann 
surfaces with a section (representing the basepoint in each fiber).
Mutatis mutandis this is a special case of our construction. 

In the first two sections we chose to gather a few technical 
proofs in appendices in order not to cloud the simplicity of 
the construction. We recommend to skip these appendices for the 
first reading.

\subsection*{Acknowledgements:}
This paper evolved from my thesis which was accepted at the 
Catholic University of Nijmegen. I am grateful to my advisor,
Joseph Steenbrink for his continuous encouragement. 
I owe very much to many valuable conversations with Richard Hain and
David Reed during a stay of three months at Duke University in spring
1996 and during a short visit of one week in spring 1997.
As mentioned in the introduction,
Richard Hain made some very useful suggestions for the problem I was 
looking at. Finally I would like to thank many people with whom I 
had interesting discussions and who gave encouragement during the preparation
of my thesis. 

\pagebreak


\section{The Nearby Fundamental Group} \label{Rietgans}

In this section we describe our concept of {\it nearby fundamental group}.

\subsection{The Setting}  \label{Setzauf}

Before we study plane curve singularities let us consider the folllowing 
setting.  We are going to use this setting throughout the whole first 
three sections of the article. 

Let 
\[
h:(Z,\; D^+) \longrightarrow (\Delta,\; 0),
\]
be a map of space germs, where
\begin{itemize}
\item
$Z$ is a complex manifold of dimension 2,
\item
$\Delta =\{z\in \cpxm \big| |z|<1\}$,
\item
$Z_0:=h^{-1}(0)=\bigcup_{i\ge 0} D_i$ is a connected
reduced divisor with normal
crossings (DNC) with components $D_i$ for ${i\ge 0}$. The components 
intersect
each other mutually in either one or no point. If it exists, the 
intersection point
between $D_k$ and $D_l$ for $k<l$ is denoted by $p_{kl}$. 
\item 
$D^+:=\bigcup_{i> r-1} D_i$ is the union of closed compact Riemann
surfaces $D_i$. 
\item
$D_0,\ldots,D_{r-1}$ are disjoint open disks ($r\ge 1$). Each disk $D_i$ 
intersects $D^+$ in one point $p_i$. Note that each point $p_i\in
\{p_0,\ldots,p_{r-1}\}$ is in particular an intersection point of two 
components,
i.~e.~$p_i=p_{il}\}$ for some $l>r-1$. 
\item
Let $Z^*:=Z \setminus h^{-1}(0)$ and $\Delta^* := \Delta \setminus \{ 0\}$
as well as $h^*:=h_{\big| Z^*}$. Then
$h^*:Z^* \longrightarrow \Delta^*$ is a
locally trivial $C^\infty$-fibration, where the closure of each fiber
(in a bigger representative of $(Z,D^+)$)
is a compact Riemann surface with $r$ boundary components.
\end{itemize}
Moreover let us assume that we are given two nonzero tangent vectors
\[
\vec{v}\in T_0 \Delta \qquad\text{ and } \qquad \vec{w}\in T_{p_0} D_0.
\]

We make the following convention:\\
\vspace{-0.7cm} 
\begin{tabbing}
???????????? \= ??????????????????????????????? \kill \\
{\bf ``[k$<$l]"} \qquad \> means ``all pairs $k,\, l$ for which 
$D_k\cap D_l \neq \emptyset$ and $k<l$". \\
{\bf ``[a$<$k$<$l]"} \>  means ``all pairs $k,\, l$ for which
$D_k\cap D_l \neq \emptyset$ and $a<k<l$".
\end{tabbing}
For $i\ge 0$ define
\[
P_i:=\bigcup_{j\neq i} D_i\cap D_j 
\qquad \text{ and } \qquad D_i^*:= D_i\setminus P_i.
\]
In section \ref{Kranenburg} we will extract such data
like above with 
additional actions of finite cyclic groups from
a plane curve singularity modulo right-equivalence.

\begin{remark} \rm
The construction which we are going to present can similarly be
considered in the case $r=0$. That is, a degenerating family of
compact Riemann surfaces, where the role of tangent vector 
$\vec{w}\in T_{p_0} D_0$  
is played by a point ${w}\in D^+ \setminus \bigcup_{[k<l]}\{p_{kl}\}$.
The requirement that two components intersect in at most one point as well as 
that the components $D_i$ have no self-intersection can easily be dropped. 
In our exposition it keeps the notation simple and it will be satisfied
in the case of curve singularities. 
\end{remark}

\subsection[$\pi_1$ of the Nearby Fiber over a Tangent Vector]
{Fundamental Group of the Nearby Fiber over a Tangent Vector}\label{Verspannung}

Let us explain in this subsection, what we mean by 
{\it the fundamental group of the nearby fiber over a tangent vector}. 

\subsubsection{The Hessian at a Double Point}

Consider the Hessian of $h$ in a double point $p_{kl}=D_k\cap D_l$. It is
a symmetric bilinear form:
\[
H(h): T_{p_{kl}} Z \times T_{p_{kl}} Z \longrightarrow T_0 \Delta.
\]
Since $D_k$ and $D_l$ intersect transversally we have
$ T_{p_{kl}} Z = T_{p_{kl}} D_k \oplus T_{p_{kl}} D_l $.
Note that $H(h)(T_{p_{ij}} D_i \times T_{p_{ij}} D_i)=0 $ for
$\{i,j\}=\{k,l\}$, that is to say that $T_{p_{kl}} D_k$ and $T_{p_{kl}} D_l$
are isotropic subspaces of $H(h)$ in $T_{p_{kl}} Z$.
Therefore, the following restriction of the
Hessian determines it completely:
\[
\eta_{kl}: T_{p_{kl}} D_k \times T_{p_{kl}} D_l \longrightarrow T_0 \Delta.
\]
Putting all such maps for all double points together, we obtain a map:
\[
\langle \cdot,\cdot\rangle :
\coprod\limits_{[k<l]} T_{p_{kl}} D_k \times T_{p_{kl}} D_l
\longrightarrow T_0 \Delta.
\]
When we are given local coordinates $(x,y)$ around a double point $p_{kl}$
in $Z$ and a coordinate $t$ on $\Delta$ such that in these coordinates
\[
h(x,y)=xy=t,
\]
then we find
\[
\left\langle a\frac{\partial}{\partial x},\,
b\frac{\partial}{\partial y}\right\rangle = a\cdot b \;
\frac{\partial}{\partial t}.
\]

\subsubsection{Paths in $Z_0$ over $\vec{v}$}

Recall that $\vec{v}$ is a tangent vector in $T_0\Delta$.
We define a {\it path in $Z_0$ over $\vec{v}$} or a
{\it path in $Z_{\vec{v}}$}
to be a continuous map
\[
\gamma: [a,b]\longrightarrow Z_0,
\]
which, when considered as map
to $Z$, is piecewise smooth and which satisfies the following conditions:
\begin{enumerate}
\item
$\gamma$ is right- and left differentiable at all double points 
$p_{kl}$, $[k<l]$.
\item
For each double point $p_{kl}$ with $[0<k<l]$ there are
two possibilities:
\begin{enumerate}
\item
either $\gamma$ stays in the same component while passing $p_{kl}$, or
\item
it changes from one component to another. 
\end{enumerate}
\end{enumerate}
Assume that $\tau_0\in
{[a,b]}$ is such that $\gamma(\tau_0) = p_{kl}$.
If $\gamma$ stays in the same component we require (and also for $p_0$)
\[
 -\dot{\gamma}^{\le \tau_0}(\tau_0) = \dot{\gamma}_{\ge \tau_0}(\tau_0)\neq 0
\]
and if $\gamma$ changes from one component to another we impose
\[
\left\langle -\dot{\gamma}^{\le \tau_0}(\tau_0),\;
\dot{\gamma}_{\ge \tau_0}(\tau_0) \right\rangle = \vec{v},
\]
where ${\gamma}^{\le \tau_0}:[a;\tau_0]\rightarrow Z_0$ and
${\gamma}_{\ge \tau_0}:[\tau_0;b]\rightarrow Z_0$ denote the
respective restrictions of $\gamma$. If $\tau_0$ is $a$
resp.~$b$, then we define $\dot{\gamma}^{\le a}(a) := \dot{\gamma}^{\le b}(b)$
resp.~$\dot{\gamma}_{\ge b}(b):=\dot{\gamma}_{\ge a}(a)$.

\subsubsection{Closed Paths in $Z_{0}$ over $\vec{v}$ based at $\vec{w}$}

Let $\vec{w}$ be a tangent vector in $T_{p_0} D_0$. We define
{\it a path over $\vec{v}$ based at $\vec{w}$} to be a path
$\gamma:[a,b]\rightarrow D^+\subset Z_0$
over $\vec{v}$ with $\gamma(a)=\gamma(b)=p_0$ such that holds:
\[
\langle \vec{w},\, \dot{\gamma}_{\ge a}(a) \rangle = 
\vec{v} \quad \text{ and }\quad
\langle -\dot{\gamma}^{\le b}(b),\, \vec{w}\rangle = \vec{v}.
\]
It follows that: $\dot{\gamma}_{\ge a}(a) = -\dot{\gamma}^{\le b}(b)$.

\subsubsection[Homotopy of Paths in $Z_0$ over $\vec{v}$]{Homotopy of Paths 
in $Z_0$ over $\vec{v}$ Relative to their Endpoints}

Let $p$ and $q$ be two points in $Z_0$ and let
$\gamma_0$ and ${\gamma_1}$ be two paths over $\vec{v}$ both starting in
$p$ and ending in $q$.
We say that there is a {\it homotopy over $\vec{v}$ relative
to their endpoints $p$ and $q$} between $\gamma_0$ and $\gamma_1$,
if there are continuous maps
$a,\, b:I\rightarrow \reellm$ from the unit interval $I$ to $\reellm$ 
with $a(s)<b(s)$ for all $s\in I$ and
\[
H:\left\{(t,s)\big|s\in I,\,t\in [a(s),b(s)]\right\}
\longrightarrow Z_0
\]
such that $H(\cdot, 0)=\gamma_0$ and $H(\cdot, 1)={\gamma_1}$ as well as
that for any $s\in I$ the map $H(\cdot, s):[a(s),b(s)]\rightarrow Z_0$
is a path over $\vec{v}$ with endpoints $p$ and $q$.
If the context is clear we will call two paths in $Z_0$ over $\vec{v}$
with the same endpoints {\it nearby homotopic relative to their endpoints},
if there is a homotopy over $\vec{v}$ relative to their endpoints
between them. Note that the notion of nearby homotopy relative to the
endpoints $p$ and $q$ as defined above yields an equivalence relation on
the set of all paths in $Z_0$ over $\vec{v}$ with endpoints $p$ and $q$.

\subsubsection{Homotopy of Paths in $Z_0$ over $\vec{v}$ based at $\vec{w}$}

We say that there is a {\it homotopy over $\vec{v}$ based at $\vec{w}$}
between two paths $\gamma_0$ and ${\gamma_1}$ in $Z_0$ over $\vec{v}$
based at $\vec{w}$, if there are continuous maps
$a,\, b:I\rightarrow \reellm$ with $a(s)<b(s)$ for all $s\in I$ and
\[
H:\left\{(t,s)\big|s\in I,\,t\in [a(s),b(s)]\right\}
\longrightarrow Z_0
\]
such that $H(\cdot, 0)=\gamma_0$ and $H(\cdot, 1)={\gamma_1}$ as well as
that for any $s\in I$ the map $H(\cdot, s):[a(s),b(s)]\rightarrow Z_0$
is a path over $\vec{v}$ based at $\vec{w}$.
If there is no chance of confusion we will call two paths in
$Z_0$ over $\vec{v}$ based at
$\vec{w}$ {\it nearby homotopic} when there is a homotopy over $\vec{v}$
based at $\vec{w}$ between them.
Note that the notion of {\it nearby homotopy} as defined above gives
an equivalence relation on the set of all paths in $Z_0$ over $\vec{v}$
based at $\vec{w}$.

Moreover, remark that a homotopy over $\vec{v}$ based at $\vec{w}$
between two paths over $\vec{v}$ based at $\vec{w}$ is in particular
a homotopy over $\vec{v}$ relative to the basepoint $p_0$ between these paths.

\afterpage{\clearpage
\begin{figure}[H]
\qquad\quad\qquad\qquad \epsfig{file=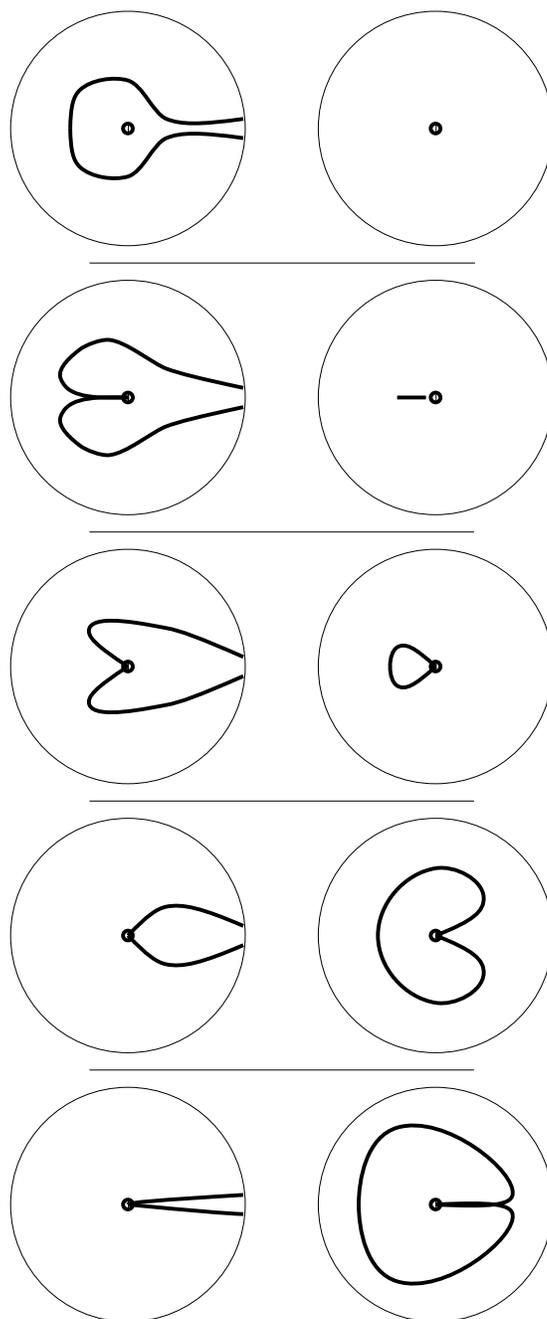,height=17.5cm,clip=}
\caption{The disks on the left (right) show a coordinate plot in a
component $D_k$ ($D_l$) with coordinate $x$ ($y$). In these coordinates
$h$ is given by $h(x,y)=x\cdot y$. The centers of the disks represent the
intersection point $p_{kl}=D_k\cap D_l$.} 
\end{figure}}

\begin{proposition} \label{Huibuh} $\text{    }$
\begin{enumerate}
\item
Every path $\gamma:[a,b]\rightarrow  Z_0$ over $\vec{v}$ (based at $\vec{w}$)
is nearby homotopic relative to its endpoints (resp.~based at $\vec{w}$)
to a path $\gamma_I:I\rightarrow Z_0$
over $\vec{v}$ (based at $\vec{w}$).
\item
Assume, two paths over $\vec{v}$, say
$\gamma_0,\,{\gamma_1}:I\rightarrow Z_0$ with $\gamma_0(0)=\gamma_1(0)$
and $\gamma_0(1)=\gamma_1(1)$ (resp.~both based at $\vec{w}$),
are nearby homotopic relative to their endpoints (resp.~based at $\vec{w}$).
Then there is a nearby homotopy
\[
H:I\times I\longrightarrow Z_0
\]
such that $H(\cdot,0)=\gamma_0$ and $H(\cdot,1)={\gamma_1}$ as well as
$H(\cdot,s)$ is a path over $\vec{v}$ with the endpoints
$\gamma_0(0)=\gamma_1(0)$ and $\gamma_0(1)=\gamma_1(1)$ (resp.~based
at $\vec{w}$) for all $s\in I$.
\end{enumerate}
\end{proposition}
\noindent
For the proof of this proposition one has to reparametrize the paths 
and the homotopies in such a way that the speed stays unchanged 
just in small neighbourhoods around the double points. For a detailled
proof we refer to \cite{Doktorarbeit}, Prop.~5.3.  

\pagebreak
\noindent
{\bf Remark:}
Observe that for example each $\gamma$ for which there is a
$\lambda\in [0,\frac{1}{2}[$ such that
\[
\gamma(\frac{1}{2}-\tau) = \gamma(\frac{1}{2} + \tau) \qquad
\forall\;\tau\in [0;\lambda]
\]
is {\it homotopic over $\vec{v}$ based at $\vec{w}$} with the path that
does not go all the way to $\gamma(\frac{1}{2})$ but returns already from
$\gamma(\frac{1}{2}-\lambda)=\gamma(\frac{1}{2}+\lambda)$, that is the path
\[
\topa{\tilde{\gamma}:}{ } \topa{[0;\,1-2\lambda]}{ t }\topa{\longrightarrow}{
\longmapsto} \topa{Z_0}{\left\{
\topa{\gamma(t)}{\gamma(t+2\lambda)}
\topa{\text{ for } t\le \frac{1}{2}-\lambda}{\text{ for } t\ge
\frac{1}{2}-\lambda }\right.}.
\]
Let $a(s)\equiv 0$ and $b(s)= 1-2s\lambda$. Then the homotopy is
given by
\[
\topa{H}{ } \topa{\left\{(t,s)\big|s\in I,\, t\in [a(s),b(s)]\right\}}{t }
\topa{\longrightarrow}{\longmapsto} 
\topa{Z_0}{\left\{
\topa{\gamma(t)}{\gamma(t+2s\lambda)}
\topa{\text{ for } t\le \frac{1}{2}-s\lambda}{\text{ for } t\ge
\frac{1}{2}-s\lambda }\right.}.
\]

We denote the set of all `homotopy over $\vec{v}$ based at $\vec{w}$'-classes
of paths in $Z_0$ over $\vec{v}$ based at $\vec{w}$ by:
\[
\pi_1(Z_{\vec{v}},\,\vec{w}).
\]
Again, if the context is clear we will refer to elements of this set
as to {\it nearby homotopy classes}. $\pi_1(Z_{\vec{v}},\,\vec{w})$ itself
will be called {\it the nearby fundamental group
(over $\vec{v}$ based at $\vec{w}$)}. The term `group' is justified by the
following proposition.

\begin{proposition}
$\pi_1(Z_{\vec{v}},\,\vec{w})$ with the composition
of `paths in $Z_{0}$ over $\vec{v}$ based at $\vec{w}$' is a group.
\end{proposition}

\noindent
{\bf Proof:} Note first that the composition of two paths in $Z_0$
over $\vec{v}$
based at $\vec{w}$ is again a path in $Z_0$ over $\vec{v}$ based at $\vec{w}$.
Furthermore the composition of such paths gives a well-defined operation
on the set of homotopy classes, which is associative.

All paths $\gamma:I\rightarrow Z_0$ over ${\vec{v}}$ based at $\vec{w}$ with
$\gamma(\frac{1}{2}-\tau) = \gamma(\frac{1}{2} + \tau)$ for all
$\tau\in [0;\frac{1}{2}]$ are mutually nearby homotopic.
Call their nearby homotopy class $e$ and note that it is the
neutral element with respect to the above defined operation on homotopy
classes.
Given an element $[\gamma]\in \pi_1(Z_{\vec{v}},\,\vec{w})$, the inverse
is represented by the path $t\mapsto \gamma(1-t)$ for $t\in I$.
\qed

\subsection[The Local System of Nearby $\pi_1$'s]
{The Local System of Nearby Fundamental Groups} \label{Lennaert}

The association of $\pi_1(Z_{\vec{v}},\,\vec{w})$ to a pair of non zero
tangent vectors $\vec{v}\in (T_0\Delta)^* := T_0\Delta\setminus\{0\}$
and $\vec{w}\in (T_{p_0} D_0)^* := T_{p_0} D_0 \setminus\{0\}$
gives rise to a local system\footnote{By a {\it local system of groups 
(vector spaces etc.) on a topological space $Y$} we mean a contravariant
functor from the fundamental groupoid $\Pi(Y)$ of $Y$ to the category 
of groups (vector spaces etc.). Note that we multiply paths in their
natural order, that is $\alpha * \beta$ means `traverse $\alpha$, then 
$\beta$'.} 
of groups on $(T_0\Delta)^* \times
(T_{p_0} D_0)^*$ as we want to point out in the sequel.
We call this local system {\it the local system of nearby fundamental
groups} and denote it by:
\[ 
\left\{\pi_1(Z_{\vec{v}},\,\vec{w})   
\right\}_{(\vec{v},\vec{w})\in (T_0 \Delta)^*\times (T_{p_0} D_0)^*}.
\]

For a homotopy class modulo endpoints represented by a path
\[
\vec{\eta}=\vec{\alpha}\times\vec{\beta}:I\rightarrow (T_0\Delta)^* \times
(T_{p_0} D_0)^*,
\]
i.~e.~a morphism of the fundamental groupoid of
$(T_0\Delta)^* \times (T_{p_0} D_0)^*$, we define
\[
\vec{\eta}_*: \pi_1(Z_{\vec{\alpha}(0)},\,\vec{\beta}(0))\longrightarrow
	\pi_1(Z_{\vec{\alpha}(1)},\,\vec{\beta}(1))
\]
in the following way. Let $\gamma$, a path over $\vec{\alpha}(0)$
based at $\vec{\beta}(0)$, represent an element in
$\pi_1(Z_{\vec{\alpha}(0)},\,\vec{\beta}(0))$. Moreover, let
$H:I\times I\rightarrow Z_0$ be a homotopy with the property that
$H(\cdot,0)=\gamma$ and $H(\cdot,s)$ is a path over
$\vec{\alpha}(s)$ based at $\vec{\beta}(s)$ for all $s\in I$.
Note that such a homotopy always exists. Then define:
\[
\vec{\eta}_*\left([\gamma]\right) := [H(\cdot,1)].
\]
In the following proposition we state that this $[H(\cdot,1)]$
is well-defined and does not depend on the special choice of $H$.

\begin{proposition} \label{Firstencel}
Let $\vec{\eta}_i=\vec{\alpha}_i\times\vec{\beta}_i:I\rightarrow
(T_0\Delta)^* \times (T_{p_0} D_0)^*$ for $i=1,\,2$ be two paths
which are homotopic modulo endpoints. Let moreover $H_i:I\times I
\rightarrow Z_0$ for $i=1,\,2$ be two homotopies such that
$H_i(\cdot,s)$ is a path over $\vec{\alpha}_i(s)$ based at
$\vec{\beta}_i(s)$ for any $s\in I$. 
If there is a homotopy over $\vec{\alpha}_1(0) = \vec{\alpha}_2(0)$
based at $\vec{\beta}_1(0) = \vec{\beta}_2(0)$ between $H_1(\cdot,0)$
and $H_2(\cdot,0)$, then there is a homotopy over
$\vec{\alpha}_1(1) = \vec{\alpha}_2(1)$ based at
$\vec{\beta}_1(1) = \vec{\beta}_2(1)$ between $H_1(\cdot,1)$
and $H_2(\cdot,1)$.
\end{proposition}
\noindent
The proof of this proposition can be found in the Appendix \ref{proofs}.

\subsubsection{The Local System of Fundamental Groups} 

The locally trivial fibration $h:Z^*\rightarrow \Delta^*$ together with
a section $\sigma:\Delta\rightarrow Z$ defines a local system on $\Delta^*$
by associating $\pi_1(Z_t,\sigma(t))$ to a point $t\in\Delta^*$. For a
homotopy class modulo endpoints represented by a path in $\Delta^*$
one considers trivializations of the fibration over this path, which have
$\sigma$ as a constant section. Such a homotopy induces a map from
the fundamental group over the starting point to the fundamental group
over the endpoint.

However, here we would like to consider a local system on
$\Delta^*\times D_0^*$. The idea is that such a section $\sigma$
corresponds to a point $p$ in $Z_0$ by taking $p=\sigma(0)$. We want to
consider those sections $\sigma$, where $\sigma(0)\in D_0$.
In order to be able to define morphisms we perform the following
construction.

Let $t$ resp.~$p$ be coordinates on $\Delta$ resp.~$D_0$ such that
$\Delta$ resp.~$D_0$ are open disks with radius $1$ in these
coordinates. Consider the map
\[
\topa{\pi:}{ }\topa{\Delta\times D_0^*}{(t,p)}
\topa{\longrightarrow}{\longmapsto}
\topa{\Delta}{|p|\cdot t.}
\]

Observe that all fibers of this map are punctured disks although $\pi$
is not a locally trivial fibration. Nonetheless, if we define for any
$\delta>0$ the annulus
$D_0^\delta:=\left\{p\in D_0 \big| |p|>\delta \right\}$, then
the restriction of $\pi$ to $\Delta\times D_0^\delta$ is a locally trivial
fibration.

\afterpage{\clearpage
\begin{figure}[H]
\centering\epsfig{file=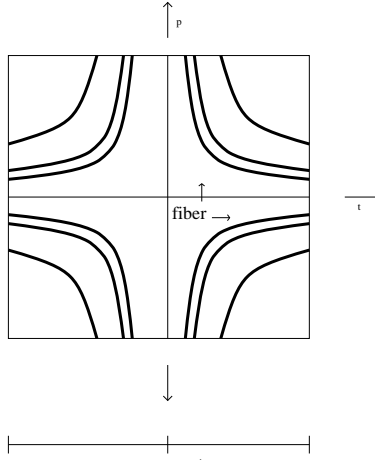,width=5cm,clip=}
\caption{Picture of the map $\pi$ (restricted to the real coordinates
in $\Delta\times D_0^*$ and $\Delta$).}
\end{figure}}

Let now $P:\Delta\times D_0^* \rightarrow Z$ be a continuous map
such that the following diagram commutes:
\[
\begin{CD} 
\Delta\times D_0^* @>P>> Z \\
& \!\!\!\!\!
\pi\bigsearrow \quad \bigswarrow h \\
& \!\!\!\!\Delta & 
\end{CD}
\]
and such that $P_{\big| \{0\}\times D_0^*}:D_0^*\rightarrow D_0^*$
is the identity.

For $\delta >0$ we define
$\Delta_\delta:=\left\{t\in \Delta\big| |t|<\delta \right\}$ and
$\Delta^*_\delta:=\Delta_\delta\setminus\{0\}$. 
Then for any $p\in D_0^\delta$ we have a section of
$\pi_{\big| \Delta\times D_0^\delta}$
\[
\topa{\check{\sigma}_p:}{}\topa{\Delta_\delta}{\tau}
\topa{\longrightarrow}{\longmapsto}
\topa{\Delta\times D_0^\delta}{\left(\frac{\tau}{|p|}, p \right).}
\]
Composed with $P$ we obtain a section of $h_{\big| h^{-1}(\Delta_\delta)}$,
which we denote by $\sigma_p:\Delta_\delta\rightarrow Z$. Note 
$\sigma_p(0) = p$.

In the following, we will for any $0<\delta <1$ define a local system of
groups on $\Delta_\delta^* \times D_0^\delta$. This system will be called
{\it the local system of fundamental groups} and we denote 
it by:
\[
\left\{\pi_1(Z_{t},\,\sigma_p(t))   
\right\}_{(t,p)\in \Delta_\delta^* \times D_0^\delta}.
\]

To a point $(t,p)$ in $\Delta_\delta^* \times D_0^\delta$ we associate
the group $\pi_1 \left(Z_t,\sigma_p(t)\right)$.
With a homotopy class modulo endpoints represented by a path
\[
\eta=\alpha\times\beta: I\longrightarrow \Delta_\delta^* \times D_0^\delta
\]
we associate a group isomorphism
\[
\eta_*: \pi_1 \left(Z_{\alpha(0)},\sigma_{\beta(0)}(\alpha(0))\right)
\longrightarrow
\pi_1 \left(Z_{\alpha(1)},\sigma_{\beta(1)}(\alpha(1))\right).
\]
in the following way.

Let $\gamma:I\rightarrow Z_{\alpha(0)}$ represent an element of
$\pi_1 \left(Z_{\alpha(0)},\sigma_{\beta(0)}(\alpha(0))\right)$ and
let $H:I\times I\rightarrow Z$ be a continuous map such that
$H(\cdot,0)=\gamma$ and $H(\cdot,s)$ is a path in $Z_{\alpha(s)}$
based at $\sigma_{\beta(s)}\left(\alpha(s)\right)$.
Then define:
\[
\eta_*\left([\gamma]\right):= \left[H(\cdot,1)\right]\in
\pi_1 \left(Z_{\alpha(1)},\sigma_{\beta(1)}(\alpha(1))\right).
\]
Define 
\[
\stackrel{\circ}{\Delta}_\delta:=
\left\{ \left(\rho,e^{2\pi i \theta}\right)\in\reellm^{\ge 0}\times S^1
\big| \rho\cdot e^{2\pi i \theta}
\in \Delta_\delta \right\}. 
\]
In the Appendix \ref{Kolk} we will define a local system on
$\stackrel{\circ}{\Delta}_\delta\times D_0^\delta$ (and prove that it is 
well-defined) such that its restriction to
$\Delta_\delta^*\times D_0^\delta$ by means of the obvious embedding
$\Delta_\delta^*\times D_0^\delta \hookrightarrow
\stackrel{\circ}{\Delta}_\delta\times D_0^\delta$ is the local
system above.

Note that the above defined local system depends on the 
choice of the map $P$. However, it will be a consequence of 
Theorem \ref{Dr. Yoghurt} below that two different such maps $P$ 
yield isomorphic local systems.  

\subsubsection{Comparison of the two Local Systems}

Here we compare the local system of nearby fundamental
groups with the local system of fundamental groups on the fibration.
We state the following theorem.
Note that $\left(T_0\Delta\right)^*$, $\left(T_{p_0}D_0\right)^*$,
$\Delta^*_\delta$ and $D_0^\delta$ are oriented by their complex structures.

\begin{theorem} \label{Dr. Yoghurt}
Any pair of orientation preserving homotopy equivalences 
\[
\left(T_0\Delta\right)^*\longrightarrow \Delta^*_\delta \quad\text{   and   }
\quad \left(T_{p_0}D_0\right)^*\longrightarrow D_0^\delta
\]
induces an isomorphism of the local system of nearby fundamental
groups on the product 
$\left(T_0\Delta\right)^*\times\left(T_{p_0}D_0\right)^*$
with the local system of fundamental groups on 
$\Delta^*_\delta\times D_0^\delta$.
\end{theorem}
The proof of this theorem will be given in Appendix \ref{proofs}.
The link between the two local systems will be provided by a local
system which we call {\it the local system of the real blow-up}. 
Here we will make use of a
construction of A'Campo \cite{A'Campo--La-fonction}.

\subsection{Appendix} \label{proofs}

This Appendix consists of two parts. In the first part we give the 
proof of Proposition \ref{Firstencel}.
In the second we compare the two local systems of \ref{Lennaert}
by constructing a third one 
and show that this is isomorphic to both of them. 

\subsubsection{Proof of \ref{Firstencel}:} 

\afterpage{\clearpage
\begin{figure}[H]
\centering\epsfig{file=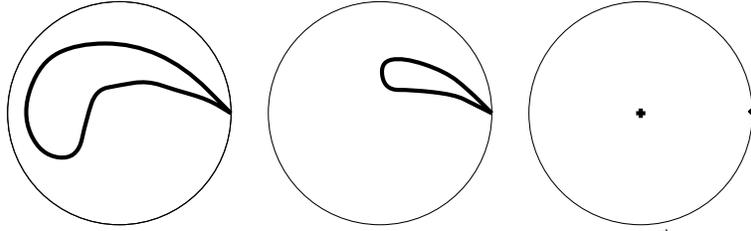,width=10cm,clip=}
\caption{Illustration of a homotopy between a path $\vec{\lambda}$ in
$T_0\Delta$ and the constant path $\vec{\lambda}(0)=\vec{\lambda}(1)$.}
\end{figure}

\begin{figure}[H]
\centering\epsfig{file=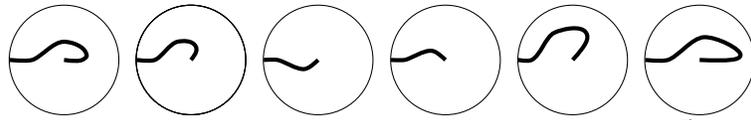,width=10cm,clip=}
\caption{Illustration of a homotopy $\mathcal H$ (over the path $\vec{\lambda}$)
at a double point (other than $p_0$) in the x-coordinate.}
\end{figure}

\begin{figure}[H]
\centering\epsfig{file=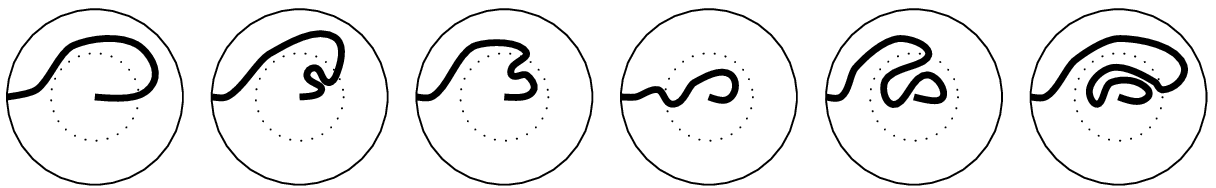,width=10cm,clip=}
\caption{Illustration of the homotopy $\tilde{\tilde{\mathcal H}}$ at a
double point (other than $p_0$) in the x-coordinate.}
\end{figure}

\begin{figure}[H]
\centering\epsfig{file=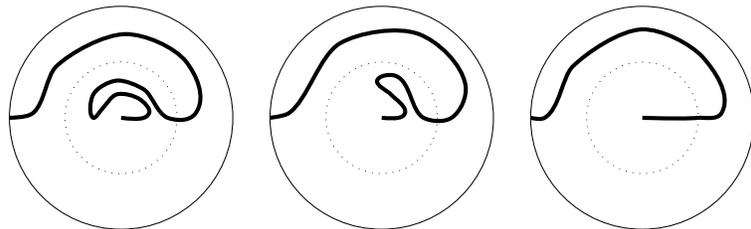,width=10cm,clip=}
\caption{Illustration of the homotopy between ${\mathcal H}(\cdot,1)$ and
$\tilde{\tilde{\mathcal H}}(\cdot,1)$ at a double point (other than $p_0$)
in the x-coordinate.}
\end{figure}}

We know that there is a homotopy, say ${\mathcal H}:I\times I\rightarrow Z_0$,
between $H_1(\cdot,1)$ and $H_2(\cdot,1)$
by using first $H_1$ as homotopy between $H_1(\cdot,1)$ and $H_1(\cdot,0)$,
then using the homotopy between $H_1(\cdot,0)$ and $H_2(\cdot,0)$ and
finally using $H_2$ as homotopy between $H_2(\cdot,0)$ and $H_2(\cdot,1)$.
The question is: Can such a homotopy be realized as homotopy
over $\vec{\alpha}_1(1) = \vec{\alpha}_2(1)$ based at
$\vec{\beta}_1(1) = \vec{\beta}_2(1)$?
Define $\vec{\lambda}:=(\vec{\alpha}_1)^{-1}\star\vec{\alpha}_i(0)\star
\vec{\alpha}_2$ and $\vec{\mu}:=(\vec{\beta}_1)^{-1}\star\vec{\beta_i(0)}\star
\vec{\beta}_2$. Then we can formulate the question more precisely:
Given closed
paths $\vec{\lambda}:I\rightarrow (T_0\Delta)^*$ and
$\vec{\mu}:I\rightarrow (T_{p_0}D_0)^*$, which are homotopic modulo endpoints
to the respective constant paths and given a homotopy
${\mathcal H}:I\times I\rightarrow Z_0$ such that ${\mathcal H}(\cdot,s)$
is a path over $\vec{\lambda}(s)$ based at $\vec{\mu}(s)$, can
we then find a homotopy $\tilde{\mathcal H}:I\times I\rightarrow Z_0$
between ${\mathcal H}(\cdot,0)$ and ${\mathcal H}(\cdot,1)$ such that
${\mathcal H}(\cdot,s)$ is a path over $\vec{\lambda}(0) = \vec{\lambda}(1)$
based at $\vec{\mu}(0) = \vec{\mu}(1)$ for all $s\in I$?

We give a positive answer to this question by doing the following
calculation locally around each double point of $Z_0$.

Let $t$ be a coordinate on $\Delta$ such that $\frac{\partial}{\partial t}
= \vec{\lambda}(0) = \vec{\lambda}(1)$. Write $\vec{\lambda}(\tau)= l(\tau)\,
\frac{\partial}{\partial t}$.
Choose around each double point $p_{kl}$ local coordinates $(x,y)$ in $Z$
such that the function $h$ in the coordinates $(x,y)$ and $t$ looks like:
$h(x,y)=x\cdot y$.
For the double point $p_0$ we may assume that $y$ is the coordinate on
$D_0$ and that $\frac{\partial}{\partial y} = \vec{\mu}(0) = \vec{\mu}(1)$.
Write $\vec{\mu}(\tau)=m(\tau)\,\frac{\partial}{\partial y}$.

Let $\vec{L}:I\times I\rightarrow (T_0\Delta)^*$ resp.~$\vec{M}:
I\times I\rightarrow (T_{p_0} D_0)^*$ be the homotopies modulo endpoints
between $\vec{\lambda}$ resp.~$\vec{\mu}$ and the constant path
$\vec{\lambda}(0) = \vec{\lambda}(1)$ resp.~$\vec{\mu}(0) = \vec{\mu}(1)$.
Write $\vec{L}(\tau,s)= {L}(\tau,s)\;\frac{\partial}{\partial t}$
and $\vec{M}(\tau,s)= {M}(\tau,s)\;\frac{\partial}{\partial y}$.

By $\sigma:\reellm^{\ge 0}\rightarrow [0;1]$ we denote a bump function which is
$\equiv 1$ near $0$ and $\equiv 0$ on $\reellm^{\ge r_0}$ for a sufficiently
small $r_0$ and strictly decreasing inbetween.

For $(\tau,s)\in I\times I$ where ${\mathcal H}(\tau,s)$ lies in an
x-coordinate, we describe $x\circ {\mathcal H}$ by polar coordinates:
$x\circ {\mathcal H}(\tau,s)=r(\tau,s)\; e^{2\pi i \varphi(\tau,s)}$.

Then modify the homotopy $\mathcal H$ to a homotopy $\tilde{\tilde{\mathcal H}}$ by
setting:
\[
x\circ \tilde{\tilde{\mathcal H}}(\tau,s):=
\frac{1}{l(\sigma(r(\tau,s))\cdot s)} \;
x\circ {\mathcal H}(\tau,s).
\]
Around the point $p_0$ we define $\tilde{\tilde{\mathcal H}}$ by:
\[
x\circ \tilde{\tilde{\mathcal H}}(\tau,s):=
\frac{m(\sigma(r(\tau,s))\cdot s)}
{l(\sigma(r(\tau,s))\cdot s)}
\; x\circ {\mathcal H}(\tau, s).
\]
Note that $\tilde{\tilde{\mathcal H}}$ is a homotopy over
$\vec{\lambda}(0) = \vec{\lambda}(1)$ based at $\vec{\mu}(0) = \vec{\mu}(1)$.

It remains to show that there is a homotopy over
$\vec{\lambda}(0) = \vec{\lambda}(1)$ based at $\vec{\mu}(0) = \vec{\mu}(1)$
between $\tilde{\tilde{\mathcal H}}(\cdot,1)$ and ${\mathcal H}(\cdot,1)$.
This homotopy is realized by the the following deformation of
${\mathcal H}(\cdot,1)$ in the neighbourhood of each double point:
\begin{align*}
\qquad\qquad\qquad &(\tau,s)\mapsto
\frac{1}{L(\sigma(r(\tau,1)), s)} \; x\circ {\mathcal H}(\tau,1) \\
\intertext{respectively around $p_0$:}
& (\tau, s)\mapsto
\frac{M(\sigma(r(\tau,1)),s)}{L(\sigma(r(\tau,1)), s)}
\; x\circ {\mathcal H}(\tau,1). \qquad \qquad\qquad \qed 
\end{align*}

\subsubsection{The Real Blow-up of the Fibration} \label{Kolk}

In this part of the Appendix we will compare the local system of 
nearby fundamental groups with the local system of fundamental groups.
The link between these two will be given by a local system, which
we call {\it the local system of the real blow-up} and which 
we are going to construct in the sequel.

We want to define an extension of $Z^*$ to a topological space
$\stackrel{\circ}{Z}$ and a map
\[
\stackrel{\circ}{h}:\stackrel{\circ}{Z}\longrightarrow \stackrel{\circ}{\Delta},
\]
in such a way that the diagram
\[
\begin{array}{rcccl}
 & \stackrel{\circ}{Z} & \hookleftarrow & Z^* & \\
\stackrel{\circ}{h} & \Big\downarrow & & \Big\downarrow & h^* \\
& \!\! \stackrel{\circ}{\Delta} & \hookleftarrow & \Delta^* &
\end{array}
\]
commutes and that $\stackrel{\circ}{h}:\stackrel{\circ}{Z}\rightarrow
\stackrel{\circ}{\Delta}$ is a locally trivial topological fibration.

Choose a coordinate $t$ on $\Delta$ and choose local coordinates
$(x,y):W_{kl}= U_{kl}^k\times U_{kl}^l\rightarrow \cpxm^2$ around any
double point
$p_{kl}\in Z$ for $[k<l]$ such that the function $h$ in these coordinates
looks like
\[
h(x,y)= x\cdot y \qquad\text{  for all } |x|<2,\, |y|<2.
\]

Consider the map
\[
\topa{\chi:}{ }\topa{\left(\reellm^{\ge 0}\times S^1 \right)}{
\left( \left( R_1, e^{2\pi i \phi_1} \right) \right. }
\topa{\times}{ ,}
\topa{\left(\reellm^{\ge 0}\times S^1 \right)}{
\left. \left( R_2, e^{2\pi i \phi_2} \right)\right)}
\topa{\longrightarrow}{\longmapsto}
\topa{\left(\reellm^{\ge 0}\times S^1 \right)}{
\left( R_1 R_2, e^{2\pi i (\phi_1 + \phi_2)} \right).}
\]
Then define for any pair $k$, $l$ with $[k<l]$:
\[
\stackrel{\circ}{U}_{kl}^k:=\left\{ \left( R, e^{2\pi i \phi}\right) \in
\reellm^{\ge 0}\times S^1 \Big|
R e^{2\pi i \phi}\in x(U_{kl}^k) \right\}
\]
and
\[
\stackrel{\circ}{U}_{kl}^l:=\left\{ \left( R, e^{2\pi i \phi}\right) \in
\reellm^{\ge 0}\times S^1 \Big|
R e^{2\pi i \phi}\in y(U_{kl}^l) \right\}.
\]
There are the obvious maps $\stackrel{\circ}{U}_{kl}^k\rightarrow U_{kl}^k$
and
$\stackrel{\circ}{U}_{kl}^l\rightarrow U_{kl}^l$.

We construct $\stackrel{\circ}{Z}$ as follows. Recall that there are
obvious maps $\stackrel{\circ}{\Delta}\rightarrow\Delta$ and
$\Delta^*\hookrightarrow\stackrel{\circ}{\Delta}$.
Define first:
\[
\breve{Z}:= Z\times_{\Delta}\stackrel{\circ}{\Delta}.
\]
In $\breve{Z}$ we have the subspaces
\begin{eqnarray*}
\breve{W}_{kl} & :=  &W_{kl}\times_{\Delta} \stackrel{\circ}{\Delta}
= \left( U_{kl}^k\times U_{kl}^l \right) \times_{\Delta}
\stackrel{\circ}{\Delta}\\
\breve{W}^*_{kl} & := &
\left( {U^*}_{kl}^k\times {U^*}_{kl}^l \right) \times_{\Delta^*} \Delta^*.
\end{eqnarray*}
Now we use  $\chi_{\big|{\stackrel{\circ}{U_{kl}^k}\times
\stackrel{\circ}{U_{kl}^l}}}$
and the identity map to define
\[
\stackrel{\circ}{W}_{kl}:=
\left( \stackrel{\circ}{U_{kl}^k}\times\stackrel{\circ}{U_{kl}^l} \right)
\times_{\stackrel{\circ}{\Delta}} \stackrel{\circ}{\Delta}.
\]
Finally we construct $\stackrel{\circ}{Z}$ by replacing $\breve{W}_{kl}$
in $\breve{Z}$ by $\stackrel{\circ}{W}_{kl}$, glued along the
common subspace $W^*_{kl}$.
The projection to the second factor in the above mentioned fibered
products yields the well-defined map we are looking after:
\[
\stackrel{\circ}{h}:\stackrel{\circ}{Z}
\longrightarrow\stackrel{\circ}{\Delta}.
\]
Note that the maps $\stackrel{\circ}{U}_{kl}^k\rightarrow U_{kl}^k$ and
$\stackrel{\circ}{U}_{kl}^l\rightarrow U_{kl}^l$ induce
a map $\stackrel{\circ}{Z} \rightarrow Z$ over $\stackrel{\circ}{\Delta}
\rightarrow \Delta$.

\noindent
{\bf Remark:}
It seems that this construction was considered first by 
A'Campo in \cite{A'Campo--La-fonction}.
There is a coordinate-free way to construct $\stackrel{\circ}{Z}$
by using the theory of log-structures (cf.~ e.~ g.~ \cite{Kato-Nakayama},
\cite{Illusie}, \cite{Kato}, \cite{log-Steenbrink}).
We decided not to introduce  the language of
log-structures here since we use the space $\stackrel{\circ}{Z}$
just as a link between the two earlier defined local systems
and because we need the coordinates anyway.
The following theorem is proved in greater generality by Usui
\cite{Usui}.

\begin{theorem}
$\stackrel{\circ}{h}:\stackrel{\circ}{Z} \rightarrow \stackrel{\circ}{\Delta}$
is a locally trivial topological fibration.
\end{theorem}

This theorem and its proof can be found in in A'Campo 
\cite{A'Campo--La-fonction} p.~241. For a direct proof in this special 
situation we refer to \cite{Doktorarbeit}, Thm.~5.6, p.~81.

Having this theorem at hand we can define a third local system of
groups -- this time on $\stackrel{\circ}{\Delta}_\delta\times D_0^\delta$.
The purpose of this local system here is that it serves as link
between the local system of nearby fundamental groups and the
local system of fundamental groups on the fibration. This third
local system will be refered to as 
{\it the local system of the real blow-up}. 
We denote it by:
\[
\left\{ \pi_1\left(\stackrel{\circ}{Z}_{(\rho, e^{2\pi i\theta})},\,
\stackrel{\circ}{\sigma}_{p}(\rho, e^{2\pi i\theta})\right)   
\right\}_{\big((\rho, e^{2\pi i\theta}),p\big)\in
\stackrel{\circ}{\Delta}_\delta\times D_0^\delta}
.\]
It is an extension
of the local system of fundamental groups on the fibration on
$\Delta^*_\delta \times D_0^\delta$ to
$\stackrel{\circ}{\Delta}_\delta\times D_0^\delta$ and it is defined as 
follows:

Observe that for $p\in D_0^\delta$ the section
$\sigma_p:\Delta_\delta\rightarrow Z$ can be pulled back in a unique way
to a section $\stackrel{\circ}{\sigma}_p:
\stackrel{\circ}{\Delta}_\delta\rightarrow \stackrel{\circ}{Z}$.
To any $\left((\rho, e^{2\pi i\theta}),p\right)\in
\stackrel{\circ}{\Delta}_\delta\times D_0^\delta$ we associate the group
\[
\pi_1\left(\stackrel{\circ}{Z}_{(\rho, e^{2\pi i\theta})},\,
\stackrel{\circ}{\sigma}_{p}\left(\rho, e^{2\pi i\theta}\right)
\right)
\]
and to a homotopy class modulo endpoints represented by a path
\[
\eta=\alpha\times\beta: I\longrightarrow
\stackrel{\circ}{\Delta}_\delta\times D_0^\delta
\]
we associate similarly as before a group isomorphism
\[
\eta_*:
\pi_1\left(\stackrel{\circ}{Z}_{\alpha(0)},\,
\stackrel{\circ}{\sigma}_{\beta(0)}\left(\alpha(0)\right)\right)
\longrightarrow
\pi_1\left(\stackrel{\circ}{Z}_{\alpha(1)},\,
\stackrel{\circ}{\sigma}_{\beta(1)}\left(\alpha(1)\right)\right).
\]

Let $\gamma: I\longrightarrow \stackrel{\circ}{Z}$
be an element of $\pi_1\left(\stackrel{\circ}{Z}_{\alpha(0)},\,
\stackrel{\circ}{\sigma}_{\beta(0)}\left(\alpha(0)\right)\right)$ and let
$H:I\times I\longrightarrow \stackrel{\circ}{H}$ be a continuous
map such that $H(\cdot,0)=\gamma$ and for any $s\in I$ is
$H(\cdot,s)$ a path in $\stackrel{\circ}{Z}_{\alpha(s)}$
based at $\stackrel{\circ}{\sigma}_{\beta(s)}\left(\alpha(s)\right)$.
Then define:
\[
\eta_*\left([\gamma]\right):=
\left[H(\cdot,1)\right] \in \pi_1\left(\stackrel{\circ}{Z}_{\alpha(1)},\,
\stackrel{\circ}{\sigma}_{\beta(1)}\left(\alpha(1)\right)\right).
\]

The following proposition justifies this definition and moreover
the definition of the local system of fundamental groups of the fibration.

\begin{proposition}
Let $\eta_i=\alpha_i\times \beta_i:I \rightarrow
\stackrel{\circ}{\Delta}_\delta\times D_0^\delta$ for $i=1,\,2$ be two paths,
which are homotopic modulo endpoints. Let moreover
$H_i:I\times I\rightarrow \stackrel{\circ}{Z}$
for $i=1,\, 2$ be two homotopies such
that $H_i(\cdot,s)$ is a path in $\stackrel{\circ}{Z}_{\alpha_i(s)}$
based at $\stackrel{\circ}{\sigma}_{\beta_i(s)}\left(\alpha_i(s)\right)$
for any $s\in I$.
If $H_1(\cdot,0)$ and $H_2(\cdot,0)$ represent the same homotopy class
in
\[
\pi_1\left(\stackrel{\circ}{Z}_{\alpha_i(0)},\,
\stackrel{\circ}{\sigma}_{\beta_i(0)}\left(\alpha_i(0)\right)\right)
\]
with $i=1$
or $2$, then $H_1(\cdot,1)$ and $H_2(\cdot,1)$ represent the same homotopy
class in
\[
\pi_1\left(\stackrel{\circ}{Z}_{\alpha_i(1)},\,
\stackrel{\circ}{\sigma}_{\beta_i(1)}\left(\alpha_i(1)\right)\right)
\]
with $i=1$ or $2$.
\end{proposition}

\noindent
{\bf Proof:}
Similar as in the proof of Proposition \ref{Firstencel} the assertion may 
be reduced to:
If $\lambda:I\rightarrow \stackrel{\circ}{\Delta}_\delta$ and
$\mu:I\rightarrow D_0^\delta$ are closed paths, which are homotopic
modulo endpoint to the respective constant paths and given a homotopy
${\mathcal H}:I\times I\rightarrow \stackrel{\circ}{Z}$
such that ${\mathcal H}(\cdot,s)$ is a path in
$\stackrel{\circ}{Z}_{\lambda(s)}$ based at
$\stackrel{\circ}{\sigma}_{\mu(s)}\left(\lambda(s)\right)$
for any $s\in I$, then there is a
$\tilde{\mathcal H}:I\times I\rightarrow \stackrel{\circ}{Z}$, a homotopy
between ${\mathcal H}(\cdot,0)$ and ${\mathcal H}(\cdot,1)$ such that
$\tilde{\mathcal H}(\cdot,s)$ is a path in
$\stackrel{\circ}{Z}_{\lambda(0)}=\stackrel{\circ}{Z}_{\lambda(1)}$
based at $\stackrel{\circ}{\sigma}_{\mu(0)}\left(\lambda(0)\right)=
\stackrel{\circ}{\sigma}_{\mu(1)}\left(\lambda(1)\right)$
for any $s\in I$.

Observe first that there is a homotopy $\tilde{\tilde{\mathcal H}}:I\times I
\rightarrow \stackrel{\circ}{Z}_{\lambda(0)}=\stackrel{\circ}{Z}_{\lambda(1)}$
between ${\mathcal H}(\cdot,0)$ and ${\mathcal H}(\cdot,1)$ lying completely in the
one fiber over $\lambda(0)=\lambda(1)$ but not necessarily leaving
the base point $\stackrel{\circ}{\sigma}_{\mu(0)}\left(\lambda(0)\right)=
\stackrel{\circ}{\sigma}_{\mu(1)}\left(\lambda(1)\right)$ fixed. This is
true for the following reason.
Let $L:I\times I \rightarrow \stackrel{\circ}{\Delta}_\delta$ be the homotopy
between $\lambda$ and the constant path $\lambda(0)=\lambda(1)$. Then
the pull-back bundle $L^*\stackrel{\circ}{Z}$ is globally trivial, i.~e.~there
is a trivialization
\[
\begin{array}{ccc} L^*\stackrel{\circ}{Z} &
\stackrel{\Phi}{\longrightarrow} &
\!\!\!
\left( I\times I \right)\times \stackrel{\circ}{Z}_{\lambda(0)}   \\
 & \!\!\! \bigsearrow L^*h \,\,\,\,\,\,
 \bigswarrow pr_2 & \\
 & & \\
& \!\! I\times I & \end{array} .
\]

The homotopy $\mathcal H$ can be considered as a map to $L^*\stackrel{\circ}{Z}$
(since ${\mathcal H}(\cdot,s)$ is a path in $\stackrel{\circ}{Z}_{\lambda(s)}$
for $s\in I$) and therefore to $\left( I\times I \right)\times
\stackrel{\circ}{Z}_{\lambda(0)}$. The projection to the second factor
is a homotopy in $\stackrel{\circ}{Z}_{\lambda(0)}$ between
${\mathcal H}(\cdot,0)$ and ${\mathcal H}(\cdot,1)$.
If the path, which is formed by the trace of the base point is zero-homotopic,
then we may find a homotopy between ${\mathcal H}(\cdot,0)$ and ${\mathcal H}(\cdot,1)$
leaving the base point $\stackrel{\circ}{\sigma}_{\mu(0)}
\left(\lambda(0)\right)=\stackrel{\circ}{\sigma}_{\mu(1)}
\left(\lambda(1)\right)$ fixed.

Note that a path in $\left( I\times I \right)\times
\stackrel{\circ}{Z}_{\lambda(0)}$ is zero homotopic if and only if its
projection to $\stackrel{\circ}{Z}_{\lambda(0)}$ is.
The above mentioned trace of the base point is the projection
to $\stackrel{\circ}{Z}_{\lambda(0)}$ of the path
\[
\phi\left({\mathcal H}(0,\cdot)\right)=\phi\left({\mathcal H}(1,\cdot)\right).
\]
But this path is zero-homotopic if and only if 
${\mathcal H}(0,\cdot)={\mathcal H}(1,\cdot)$
is. Note that ${\mathcal H}(0,s)={\mathcal H}(1,s)
=\stackrel{\circ}{\sigma}_{\mu(s)}\left(\lambda(s)\right)$.
Therefore, a homotopy to the constant path
$\stackrel{\circ}{\sigma}_{\mu(0)}\left(\lambda(0)\right)=
\stackrel{\circ}{\sigma}_{\mu(1)}\left(\lambda(1)\right)$ is given by:
$(s,\tau)\longmapsto \stackrel{\circ}{\sigma}_{M(s,\tau)}
\left(L(s,\tau)\right)$.
\qed

\bigskip
Let us conclude this section by proving Theorem \ref{Dr. Yoghurt},
which asserts that 
any pair of orientation preserving homotopy equivalences
\[
\left(T_0\Delta\right)^*\longrightarrow \Delta^*_\delta \quad\text{   and   }
\quad \left(T_{p_0}D_0\right)^*\longrightarrow D_0^\delta
\]
induces an isomorphism of the local system of nearby fundamental
groups on the product
$\left(T_0\Delta\right)^*\times\left(T_{p_0}D_0\right)^*$
with the local system of fundamental groups on
$\Delta^*_\delta\times D_0^\delta$.

\noindent
{\bf Proof of \ref{Dr. Yoghurt}:}
Observe that $\Delta^*_\delta\subset\stackrel{\circ}{\Delta}_\delta$
and call $S:=\stackrel{\circ}{\Delta}_\delta\setminus\Delta^*_\delta$.
By definition of both local systems, the inclusion
$\Delta^*_\delta\times D_0^\delta\subset\stackrel{\circ}{\Delta}_\delta
\times D_0^\delta$ induces an isomorphism of the local system of
fundamental groups on the fibration with the local system of the real blow up.
The latter local system is also equivalent to its restriction to
$S\times D_0^\delta$ by the inclusion $S\times D_0^\delta\subset
\stackrel{\circ}{\Delta}_\delta\times D_0^\delta$. We call this restriction
to $S\times D_0^\delta$ {\it the local system on the soul of the real blow-up}.

Hence, the assertion of the theorem is equivalent with:
Any pair of orientation preserving homotopy equivalences
\[
a: \left(T_0\Delta\right)^*\longrightarrow S \quad \text{   and   }\quad
b: \left(T_{p_0}D_0\right)^*\longrightarrow D_0^\delta
\]
induces an isomorphism of the local system of nearby fundamental groups to
the local system on the soul of the real blow-up.

We have to show that there is
an isomorphism between the group associated to one pair of tangent vectors
$(\vec{v}_0,\vec{w}_0)\in \left(T_0\Delta\right)^*\times
\left(T_{p_0}D_0\right)^*$
and the group associated to
$\left(a(\vec{v}_0),b(\vec{w}_0)\right)\in S\times D_0^\delta$,
which is equivariant with respect to the action of
$\pi_1 \left( \left(T_0\Delta\right)^*\times\left(T_{p_0} D_0\right)^*,\,
(\vec{v}_0,\vec{w}_0) \right)$
and the action of
$\pi_1 \left(S\times D_0^\delta,\,
\left(a(\vec{v}_0),b(\vec{w}_0)\right) \right)$.

Let us first construct the isomorphism and show that it is equivariant.
Associated to each double point $p_{kl}$ for $[k<l]$ we have from the
definition of $\stackrel{\circ}{Z}$ coordinates
$\left(R_1,e^{2\pi i\phi_1}\right)$, $\left(R_2,e^{2\pi i\phi_2}\right)$
in $\stackrel{\circ}{Z}$ and coordinates $\left(\rho,e^{2\pi i\theta}\right)$
on $\stackrel{\circ}{\Delta}$ such that in these coordinates
\[
\stackrel{\circ}{h}\left(
\left(R_1,e^{2\pi i\phi_1}\right),\,\left(R_2,e^{2\pi i\phi_2}\right)\right)
= \left(R_1 R_2,e^{2\pi i(\phi_1+\phi_2})\right).
\]
Here we may assume that $a(\vec{v}_0)=(0,1)\in S$ and $b(\vec{w}_0)=(1,0)
\in \stackrel{\circ}{\bar{D}}_0$.

We will define an isomorphism
\[
{\mathcal J}: \pi_1 \left(Z_{\vec{v}_0},\vec{w}_0 \right)\longrightarrow
\pi_1 \left(\stackrel{\circ}{Z}_{(0,1)},\stackrel{\circ}{\sigma}_{(1,0)}
\left(0,1 \right)\right).
\]
Let $[\gamma]$ be a nearby homotopy class in
$\pi_1 \left(Z_{\vec{v}_0},\vec{w}_0 \right)$ represented
by the path $\gamma$ over $\vec{v}_0$ based at $\vec{w}_0$.
Note that any piecewise smooth path in $\cpxm$ has a unique lifting 
with respect to the map
\[
\begin{matrix}\stackrel{\circ}{\cpxm}:= & 
\reellm^{\ge 0}\times S^1 & \rightarrow & \cpxm \\
 & (r,e^{2\pi i \varphi}) & \mapsto &
r\cdot e^{2\pi i \varphi} . \end{matrix}
\]

Therefore, there is a unique lifting of $\gamma \in Z_0$ to a path
${\stackrel{\bullet}{\gamma}}$ in $\stackrel{\circ}{Z}$. The condition
that $\gamma$ lies over $\vec{v}_0$ implies that
${\stackrel{\bullet}{\gamma}}$ is even in $\stackrel{\circ}{Z}_{(0,1)}$.
and that it is connected. Define then $\stackrel{\circ}{\gamma}:=
\sigma \star {\stackrel{\bullet}{\gamma}}\star \sigma^{-1}$,
where $\sigma$ is the straight path in $\stackrel{\circ}{D}_0$
defined by $t\mapsto (t, e^{2\pi i 0})= (1-t,1)$, then we obtain a path
in $\stackrel{\circ}{Z}_{(0,1)}$, which is based at
$\stackrel{\circ}{\sigma}_{(1,0)}\left( 0,1 \right)
=(1,0) \in \stackrel{\circ}{\bar{D}}_0$. \label{Hodgepodge}

The map $\gamma\mapsto \stackrel{\circ}{\gamma}$ induces a well-defined map
from the group $\pi_1 \left(Z_{\vec{v}_0},\vec{w}_0 \right)$ to
$\pi_1 (\stackrel{\circ}{Z}_{(0,1)},\stackrel{\circ}{\sigma}_{(1,0)}
\left(0,1 \right))$, which is certainly a group homomorphism.
Define $\mathcal J$ by $\mathcal J(\gamma):=\stackrel{\circ}{\gamma}$. 
It remains to show that $\mathcal J$ is an {\it iso}morphism.

\noindent
{\bf $\mathcal J$ is surjective:}
Call the sets
$\left\{(0,e^{2\pi i \varphi}),\,(0,e^{-2\pi i \varphi})\big|
\varphi \in [0;1[ \right\}\subset \stackrel{\circ}{P}_{kl}\cap
\stackrel{\circ}{Z}_{(0,1)}$ for
$[k<l]$ {\it vanishing cycles}. We can give $\stackrel{\circ}{Z}_{(0,1)}$
a differentiable structure in such a way that it coincides with
the already given differentiable structure off the vanishing cycles
coming from $Z_0\setminus \bigcup_{[k<l]} \{p_{kl} \}$.
It is well-known that any element in the fundamental group of a smooth
manifold (with boundary) can be represented by a smooth path
(cf.~e.~g.~\cite{tom-Dieck}, II, 1.10). Moreover observe that by a
local construction any such smooth path is homotopic to a smooth
path, which intersects the vanishing cycles transversally in a finite
number of parameter values, the first and the last of which are the only
intersections with the vanishing cycle corresponding to $p_0$.
Moreover, we may assume that $\stackrel{\odot}{\gamma}$ can be written
as $\stackrel{\odot}{\gamma}=\sigma \star
{\stackrel{\bullet}{\gamma}}\star \sigma^{-1}$ for some closed path
${\stackrel{\bullet}{\gamma}}$, which lies, apart from the
base point, entirely in $D^+$.

A path $\stackrel{\odot}{\gamma}$ in $\stackrel{\circ}{Z}_{(0,1)}$
with these properties, composed with the map
$\stackrel{\circ}{Z}_{(0,1)}\rightarrow {Z}_{0}$
(collapsing the vanishing cycles)
yields a path $\gamma$ in $D^+\subset Z_0$. By reparametrizing
the path $\stackrel{\odot}{\gamma}$ appropriately we may assume that
it lies over $\vec{v}_0$ and is based at $\vec{w}_0$.
And hence we find: $\stackrel{\circ}{\gamma}=\stackrel{\odot}{\gamma}$.

\noindent
{\bf $\mathcal J$ is injective:}
Here we have to convince ourselves from the following fact. If we
are given two homotopic paths $\stackrel{\circ}{\gamma}_1$ and
$\stackrel{\circ}{\gamma}_2$ in $\stackrel{\circ}{Z}_{(0,1)}$
based at $(1,0)\in\stackrel{\circ}{\bar{D}}_0$,
which are images under $\mathcal J$
of two paths $\gamma_1$ and $\gamma_2$ over $\vec{v}_0$ based at
$\vec{w}_0$, then these
latter paths are homotopic over $\vec{v}_0$ based at $\vec{w}_0$.

Let $H:I\times I\rightarrow \stackrel{\circ}{Z}_{(0,1)}$ be a homotopy
between the paths $\stackrel{\circ}{\gamma}_1$ and
$\stackrel{\circ}{\gamma}_2$ in $\stackrel{\circ}{Z}_{(0,1)}$
based at $(1,0)\in\stackrel{\circ}{\bar{D}}_0$.
Using again the differentiable
structure on $\stackrel{\circ}{Z}_{(0,1)}$ like above we may assume
that $H$ is smooth (cf.~e.~g.~\cite{tom-Dieck}, II, 1.10) and
that for each $s\in I$ the path $H(\cdot,s)$ intersects the
the vanishing cycles transversally in a finite
number of parameter values, the first and the last of which are the only
intersections with the vanishing cycle corresponding to $p_0$.
Moreover, we may assume that $H(\cdot,s)$ can be written
as $H(\cdot,s)=\sigma \star
\eta_s \star \sigma^{-1}$ for some closed path
$\eta_s$, which lies, apart from the
base point, entirely in $D^+$.

Composed with
$\stackrel{\circ}{Z}_{(0,1)}\rightarrow {Z}_{0}$
we obtain a homotopy between $\gamma_1$ and $\gamma_2$ over $\vec{v}_0$
based at $\vec{w}_0$.
Finally, we study the monodromy actions. Let
\[
\vec{\eta}=\vec{\alpha}\times\vec{\beta}: I\longrightarrow
\left(T_0\Delta\right)^*\times\left(T_{p_0}D_0\right)^*
\]
be a path representing a morphism in the fundamental groupoid of
$\left(T_0\Delta\right)^*\times\left(T_{p_0}D_0\right)^*$. Then define
$\eta:=\left( a\circ\vec{\alpha}\right)\times\left(b\circ\vec{\beta}\right):
I\longrightarrow S \times D_0^\delta$.
And let $\gamma$ be a path over $\vec{v}_0$ based at $\vec{w}_0$ representing
an element in $\pi_1 \left(Z_{\vec{v}_0},\vec{w}_0 \right)$.
It is to show:
${\mathcal J}\circ \vec{\eta}_*\left([\gamma]\right) =
{\eta}_*\circ {\mathcal J}\left([\gamma]\right)$.

Let $H:I\times I\rightarrow Z_0$ be a homotopy with the property that
$H(\cdot,0)=\gamma$ and $H(\cdot, s)$ is a path over $\vec{\alpha}(s)$
based at $\vec{\beta}(s)$ for every $s\in I$.
Like the construction of $\stackrel{\circ}{\gamma}$ from $\gamma$ on
page \pageref{Hodgepodge} we construct a homotopy $\stackrel{\circ}{H}
:I\times I \rightarrow \stackrel{\circ}{Z}$ from $H$ such that
$\stackrel{\circ}{H}(\cdot,0)=\stackrel{\circ}{\gamma}$ and such that
$\stackrel{\circ}{H}(\cdot, s)$ is a path
$\stackrel{\circ}{Z}_{a\circ\vec{\alpha}(s)}$ based at
$b\circ\vec{\beta}(s)$. This can be achieved by the definition:
$\stackrel{\circ}{H}(\tau, s):=\left( H(\cdot,s)\right)^\circ (\tau)$.
Now compute:
\[
{\mathcal J}\circ \vec{\eta}_*\left([\gamma]\right) =
{\mathcal J}\left(\left[H(\cdot,1)\right]\right) =
\left[\stackrel{\circ}{H}(\cdot,1)\right]
\text{ and }
{\eta}_*\circ {\mathcal J}\left([\gamma]\right)=
{\eta}_*\left([\stackrel{\circ}{\gamma}]\right)=
\left[\stackrel{\circ}{H}(\cdot,1)\right].
\]
\qed


\section[Line Integrals on $\pi_1(Z_{\vec{v}},\vec{w})$]
{Line Integrals on the Nearby Fundamental Group} \label{Quark}

In Section \ref{Rietgans} we defined the group
\[
\pi_1\left(Z_{\vec{v}},\vec{w}\right).
\]
Consider the group ring $\itgm\pi_1\left(Z_{\vec{v}},\vec{w}\right)$
and let $J_{\vec{v}\vec{w}}$ or simply $\vec{J}$ be the augmentation ideal.
By definition, this is the set of elements in 
$\itgm\pi_1\left(Z_{\vec{v}},\vec{w}\right)$,
whose sum of coefficients vanishes. Our aim is to put a $\itgm$-MHS on
\[
{\itgm\pi_1\left(Z_{\vec{v}},\vec{w}\right)}\Big/{\vec{J}^{s+1}}
\]
for any $s\ge 1$. This will be done by putting a MHS on
$({\vec{J}}\big/{\vec{J}^{s+1}})$.
Note that multiplication induces an isomorphism
${\vec{J}^s}\big/ {\vec{J}^{s+1}} \xrightarrow{\cong}
({\vec{J}}\big/ {\vec{J}^{\;2}})^{\otimes s}$ 
(recall that: $\pi_1\left(Z_{\vec{v}},\vec{w}\right) \cong
\pi_1\left(Z_{t},p\right)$ for some $t\in\Delta^*$ and $p\in Z_t$).
The short exact sequence
\begin{equation} \label{Hausgans}
0\rightarrow {\vec{J}^s}\big/ {\vec{J}^{s+1}}\rightarrow
{\vec{J}}\big/ {\vec{J}^{s+1}}\rightarrow {\vec{J}}\big/ {\vec{J}^{s}}
\rightarrow 0
\end{equation}
will be a short exact sequence of MHSs. 

The main ingredient for the definition of these MHSs is a differential
graded algebra (dga) $A^\bullet$, which is made up from differential
forms on the central fiber but whose cohomology is isomorphic to the cohomology
of a regular fiber. It was suggested to us by Hain\footnote{The dga 
$A^\bullet$ is in principle a modified version of the complex 
$P_\cpxm^\bullet(X^*)[\theta]$ of \cite{Hain-de-Rham-homotopy-II}
for this special situation with non compact fibers.}
to use a dga like $A^\bullet$ and he proposed the problem of 
finding a way to integrate elements of this dga somehow on the
central fiber. We will see that the notion of {\it a path over $\vec{v}$
(based at $\vec{w}$)} will be the right notion of path
for this purpose. We show that closed elements in $A^1$ can be
integrated along paths in $Z_0$ over $\vec{v}$ in a natural way. 
This yields an integral
structure on $H^1(A^\bullet)$, which is by this integration dual to
${\vec{J}}\big/ {\vec{J}^{\;2}}$. On $A^\bullet$ we have filtrations $W_\bullet$
and $F^\bullet$ such that the induced filtrations together with the integral
structure define a $\itgm$-MHS on $H^1(A^\bullet)$. It will
turn out that the integral lattice of this $\itgm$-MHS depends on
$\vec{v}$ (and not on $\vec{w}$). The dependence upon $\vec{v}$ is such that
$\{ (J_{\vec{v},\vec{w}}\big/ J_{\vec{v},\vec{w}}^{\;2})^*\}_{\vec{v}\in
(T_0\Delta)^*}$ is a nilpotent orbit of MHSs. 
 
Since $A^\bullet$ has an algebra-structure, we can also define {\it iterated
line integrals} along paths over $\vec{v}$. We shall do this in section
\ref{Bl"omscher}. This will give a way to describe
$({\vec{J}}\big/ {\vec{J}^{s+1}})^*$ in terms of the
dga $A^\bullet$ similar to Chen's $\pi_1$-De Rham theorem 
(cf.~\cite{Chen-covering}, \cite{Chen} or \cite{Hain-the-geometry}).
Finally it will lead to the definition of a MHS on
$({\vec{J}}\big/ {\vec{J}^{s+1}})^*$ such that
the short exact sequence \eqref{Hausgans} becomes a short exact sequence
of MHSs.

\subsection{A Differential Graded Algebra of Differential 
Forms on the Central Fiber} \label{Ultraschall}

Let $(\Lambda^\bullet\,,d):=
({\bigwedge}^\bullet(\frac{dx}{x},\frac{dy}{y})[\log t],\,d)$
be the dga defined in the following way. Let
${\bigwedge}^\bullet(\frac{dx}{x},\frac{dy}{y})$ be the free graded-commutative
$\cpxm$-algebra with unit, generated by the symbols $\frac{dx}{x}$,
$\frac{dy}{y}$ in degree 1.

Let then $\Lambda^\bullet:=
{\bigwedge}^\bullet(\frac{dx}{x},\frac{dy}{y})[\log t]$
be the algebra ${\bigwedge}^\bullet(\frac{dx}{x},\frac{dy}{y})$, 
where the symbol\footnote{We will use $u$ for
convenience or if confusion with the function `$\log$' might occur.}
$\log t$ or $u$ has been added in degree 0.
The differential $d$ on $\Lambda^\bullet$ is defined by the Leibniz-rule,
$d\big(\frac{dx}{x}\big)=d\big(\frac{dy}{y}\big)=0$ and
\[
d\,\left(\log t\right):=\frac{dx}{x}+\frac{dy}{y}.
\]
Note that $(\Lambda^\bullet,\,d)$ computes the cohomology of the fibers
of the map $\Delta^*\times \Delta^*\rightarrow \Delta^*$ given by 
$(x,y)\mapsto xy$.

For each component $D_i$ of $D^+$ let
\[
E^\bullet(D_i \log P_i):=\Gamma\left(\Omega^\bullet(D_i \log P_i)
\otimes_{\Omega^0(D_i)} {\mathcal E}^{0,\bullet}(D_i)\right)
\]
be the logarithmic de Rham complex. Here ${\mathcal E}^{0,\bullet}(D_i)$ denotes
the sheaf of $C^\infty$-forms of type $(0,\bullet)$. Note that
$\left(E^\bullet(D_i \log P_i),\,d\right)$ computes the cohomology of
$D_i\setminus P_i$.

Finally let ${\bigwedge}^\bullet (\frac{dp}{p})$ be the differential
graded-commutative $\cpxm$-algebra with unit generated by the symbol
$\frac{dp}{p}$ in degree 1, where the differential $d$ is defined to be
always the zero-map. $\left({\bigwedge}^\bullet (\frac{dp}{p}),\,d\right)$
computes the cohomology of $\cpxm^*$. 
Now we consider $r$ copies $\left({\bigwedge}^\bullet_i  (\frac{dp}{p}),\,d\right)\, i=0,\ldots,r-1$ of $\left({\bigwedge}^\bullet
(\frac{dp}{p}),\,d\right)$ and think of an element $g_i$ in
${\bigwedge}^0_i (\frac{dp}{p})=\cpxm$ as a function on the point $p_i$,
i.~e.~$g_i(p_i):=g_i$. We refer to the complex number $\rho$ of an element
$\rho\frac{dp}{p}\in{\bigwedge}^1_i(\frac{dp}{p})$ as 
{\it the residue in} $p_i$
of this element. Write $\mbox{Res}_{p_i}(\rho\frac{dp}{p})=\rho$.

The dga $A^\bullet$, which we are going to define will be a
sub-dga of
\[ 
B^\bullet := \bigoplus_{0\le i\le r-1}{\bigwedge}^\bullet_i(\frac{dp}{p})\oplus
\bigoplus_{i>r-1}E^\bullet(D_i\log P_i)\oplus\bigoplus_{[k<l]}
E^\bullet(\Delta^1)\otimes_\cpxm\Lambda^\bullet,
\]
where $\Delta^1$ denotes the one-dimensional standard-simplex 
(the unit interval) and
$E^\bullet(\Delta^1)=\cpxm[\xi]\oplus\cpxm[\xi]d\xi$ the dga of real 
polynomial differential forms on $\Delta^1$.

The elements of $A^\bullet$ in $B^\bullet$ will be called {\it the
compatible elements} in $B^\bullet$. Hence, we will define $A^\bullet$
now by defining, when an element of $B^\bullet$ is compatible. Observe
that $B^n=0$ for $n\ge 4$.

We shall use $P_{kl}$, $K_{kl}$, $L_{kl}$, $H_{kl}$,
$R_{kl}$, $S_{kl}$, $T_{kl}$ and $U_{kl}$ to denote elements in
$\cpxm[\xi,u]$. If the context is clear we omit the indices $k,l$.
Moreover it will be useful to expand these elements of $\cpxm[\xi,u]$
as polynomials in $u$ with coefficients in $\cpxm[\xi]$. That means
for instance
\[
P_{kl}=P_{kl} (\xi,u)=P=P_0+P_1 u+P_2 u^2+ \cdots + P_m u^m
\mbox{ and } P'=\frac{\partial P}{\partial \xi}.
\]

\noindent
$\bf A^0 :$
An element $f=\sum_{i\ge 0} g_i + \sum_{[k<l]} P_{kl} \in B^0$, where
$g_i\in{\bigwedge}^0_i(\frac{dp}{p})$ for $i\le r-1$ and $g_i\in E^0(D_i)$ 
for $i>r-1$ is called {\it compatible} iff for $[k<l]$ holds:
\[
P_{kl}(0,u)=g_k(p_{kl}) \text{ and }
P_{kl}(1,u)=g_l(p_{kl}).
\]

\noindent
$\bf A^1 :$
We call an element in $B^1$,
\[
\varphi=\sum_{i\ge 0} \omega_i + \sum_{[k<l]} \left(
K_{kl}\frac{dx}{x}
+L_{kl}\frac{dy}{y} + H_{kl} d\xi \right),
\]
where $\omega_i\in {\bigwedge}^1(\frac{dp}{p})$ for $i\le r-1$
and $\omega_i\in
E^1(D_i \log P_i)$ for $i>r-1$, a {\it compatible} element iff
for $[k<l]$ holds:
\begin{eqnarray*}
K_{kl}(0,u) = &\mbox{Res}_{p_{kl}} \omega_k &,\quad L_{kl}(0,u)=0, \\
K_{kl}(1,u) = & 0 &,\quad L_{kl}(1,u)=\mbox{Res}_{p_{kl}} \omega_l.
\end{eqnarray*}

\noindent
$\bf A^2 :$
Let us call an element
\[
\phi=\sum_{i>r-1} \Omega_i + \sum_{[k<l]} \left(R_{kl} \,d\xi\wedge\frac{dx}{x}
+S_{kl}\, d\xi\wedge\frac{dy}{y} + T_{kl}\, \frac{dx}{x}\wedge\frac{dy}{y}
\right) \in B^2,
\]
with $\Omega_i\in E^2(D_i\log P_i)$ for $i>r-1$, {\it compatible} iff for
$[k<l]$ holds:
\[
T_{kl}(0,u)=T_{kl}(1,u)=0.
\]

\noindent
$\bf A^3 :$
We define $A^3:=B^3$. That is, all elements of $B^3$ are {\it compatible}.

Observe that holds $dA^i\subset A^{i+1}$ for $i\ge 0$ and that $A^\bullet$
is a dga.

\begin{definition} \label{Weisser Schwan} \rm
On $A^\bullet$ we define an {\it augmentation map} $a$ by:
\[
\begin{matrix} a: & A^0 & \rightarrow & \cpxm \\
 & \sum_{i\ge 0} g_i + \sum_{[k<l]} P_{kl} & \mapsto & g_0. \end{matrix}
\]
\end{definition}

An alternative definition of $A^\bullet$ can be given as follows.
Consider the surjective map of complexes
\[
\Phi: B^\bullet \longrightarrow
\bigoplus_{[k<l]} \Lambda^\bullet\oplus\Lambda^\bullet,
\]
which sends $\sum_{i\ge 0} g_i + \sum_{[k<l]} P_{kl}\in B^0$ to
\[
\sum_{[k<l]}\left\{\left(P_{kl}(0,u)-g_k(p_{kl})\right)\oplus
\left(P_{kl}(1,u)-g_l(p_{kl})\right)\right\}
\]
and $\sum_{i\ge 0} \omega_i + \sum_{[k<l]} \left(K_{kl}\frac{dx}{x}
+L_{kl}\frac{dy}{y} + H_{kl} d\xi\right)\in B^1$ to
\begin{align*}
\sum\limits_{[k<l]} & \, \left\{ \left[\left(K_{kl}(0,u)-\mbox{Res}_{p_{kl}}
\omega_k\right)\frac{dx}{x} 
+ L_{kl}(0,u)\frac{dy}{y} \right] \right. \\
\oplus & \left. \quad
\left[K_{kl}(1,u)\frac{dx}{x}+\left(L_{kl}(1,u)-\mbox{Res}_{p_{kl}}\omega_l
\right)\frac{dy}{y} \right] \right\}.
\end{align*}
Moreover, $\Phi$ maps
\[
\sum_{i>r-1} \Omega_i + \sum_{[k<l]} R_{kl}\, d\xi\wedge\frac{dx}{x}
+S_{kl}\, d\xi\wedge\frac{dy}{y} + T_{kl}\, \frac{dx}{x}\wedge\frac{dy}{y}
\in B^2
\]
to $\sum_{[k<l]} T_{kl}(0,u) \frac{dx}{x}\wedge\frac{dy}{y} \oplus
T_{kl}(1,u) \frac{dx}{x}\wedge\frac{dy}{y}$. It is easy to check that $\Phi$ is
indeed a surjective map of complexes. We find:
\[
A^\bullet = \ker \Phi.
\]

Define $C^\bullet:= \sum_{[k<l]} \Lambda^\bullet\oplus\Lambda^\bullet$ and
note that we have the short exact sequence
\[
0\longrightarrow A^\bullet\longrightarrow B^\bullet\longrightarrow C^\bullet
\longrightarrow 0,
\]
which yields a long exact cohomology sequence
(as $H^i(B^\bullet)=0$ for $i\ge 2$)
\begin{multline}\label{Saatgans}
0\rightarrow H^0(A^\bullet)\rightarrow H^0(B^\bullet)
\rightarrow H^0(C^\bullet)\\ \rightarrow H^1(A^\bullet)
\rightarrow H^1(B^\bullet)\rightarrow H^1(C^\bullet)
\rightarrow H^2(A^\bullet)\rightarrow 0.
\end{multline}

The complex $(A^\bullet,\,d)$ computes the cohomology of the regular fiber
as the following theorem asserts.

\begin{theorem} \label{Graugans}
$H^\bullet(A^\bullet)\cong H^\bullet(Z_t;\cpxm)$ for any $t\in \Delta^*$.
\end{theorem}
The proof of this theorem can be found in the Appendix \ref{Menno} on page 
\pageref{Menno}.

\subsection[Integration along Paths in the Nearby Fiber] {Integration
in the Nearby Fiber along Paths over a Tangent Vector}

Here we are going to define integrals of a closed element
\[
\varphi=\sum_{i\ge 0} \omega_i +\sum_{[k<l]} \left( K_{kl} \frac{dx}{x}
+L_{kl} \frac{dy}{y} +H_{kl} d\xi\right) \in A^1
\]
along paths over $\vec{v}$.
We define those integrals in four steps by defining them first in
special cases.
\begin{description}
\item[(a)] \label{Kormoran}
Let $\gamma:[a;b]\rightarrow Z_0$ be a path over $\vec{v}$, 
which meets the set
of double points only once with parameter value $\tau_0\in ]a;b[$, where
it changes from $D_k$ to $D_l$. Let $(x,y):W_{kl}=U_{kl}^k\times U_{kl}^l
\rightarrow \cpxm^2$ and $t:\Delta\rightarrow\cpxm$ be coordinates such that
$h(x,y)=x\cdot y$ and $\vec{v}=\frac{\partial}{\partial t}$.
Define $\gamma_x(\tau)$ (resp.~$\gamma_y(\tau)$) to be $x(\gamma(\tau))$
(resp.~$y(\gamma(\tau))$) for all $\tau\in [a;b]$ with $\gamma(\tau)
\in U_{kl}^k$ (resp.~$\gamma(\tau)\in U_{kl}^l$). Then define:
\[
\int\limits_\gamma \varphi:= \lim_{\varepsilon\to 0} \left\{
\int\limits_{\gamma^{\le \tau_0-\varepsilon}} \!\!\!\!\!\!\omega_k
+\!\!\!\int\limits_0^1 H\left(\xi, \log\left(\gamma_x(\tau_0-\varepsilon)
\gamma_y(\tau_0+\varepsilon)\right)\right) d\xi +
\!\!\!\!\!\!\!\!\int\limits_{\gamma_{\ge \tau_0+\varepsilon}} \!\!\!\!\!\!
\omega_l
\right\}.
\]
\item[(b)]
Let $\gamma:[a;b]\rightarrow Z_0$ be a path over $\vec{v}$, 
which meets the set
of double points only once with parameter value $\tau_0\in ]a;b[$, where
it stays in one component $D_k$. Then define:
\[
\int\limits_\gamma \varphi:= \lim_{\varepsilon\to 0} \left\{
\int\limits_{\gamma^{\le \tau_0-\varepsilon}} \omega_k
+\int\limits_{\gamma_{\ge \tau_0+\varepsilon}} \omega_k
\right\}.
\]
\item[(c)]
Now let $\gamma:[a;b]\rightarrow Z_0$ be a path over $\vec{v}$
such that $\gamma(a)$ and $\gamma(b)$ are no double points.
Then there is a finite number of parameter values $\tau_1,\ldots,\tau_N$
with $a<\tau_1<\cdots<\tau_N<b$, which are mapped onto double points.
Choose a $\tau_i^*\in ]\tau_i;\tau_{i+1}[$ for $i=1,\ldots,N-1$ and let
$\tau_0^*:=a$ and $\tau_N^*:=b$. Then define for $i=1,\ldots,N$ the paths
$\gamma_i:= \gamma_{\big| [\tau_{i-1}^*,\tau_i^*]}$ and we define:
\[
\int\limits_\gamma \varphi:= \sum_{i=1}^{N}
\int\limits_{\gamma_i} \varphi.
\]
\item[(d)]
Finally let $\gamma:[a;b]\rightarrow Z_0$ be a path in $Z_0$ 
over $\vec{v}$
based at $\vec{w}$. The coordinate $p:D_0\rightarrow\cpxm$ allows us to 
consider $D_0$ as part of the complex plane. There, in the complex plane,
we have the differential form (and not the symbol) $\frac{dp}{p}$.
Let $\sigma$ be the straight
path $\tau\mapsto (1-\tau)$ in $\cpxm$, then we define 
\[
\int\limits_\gamma \varphi:= 
\int\limits_{\sigma\star\gamma\star \sigma^{-1}} \varphi.
\]
\end{description}

\begin{theorem} \label{Schleiereule}
For a path $\gamma:[a;b]\rightarrow Z_0$ over $\vec{v}$ based at $\vec{w}$
and a closed $\varphi\in A^1$ holds:
\begin{enumerate}
\item
$\int_\gamma \varphi$ is well-defined and finite,
\item
$\int_\gamma \varphi$ does not depend on the choice of coordinates,
\item
$\int_\gamma \varphi =0$, if $\varphi$ is exact.
\end{enumerate}
\end{theorem}

Before proving \ref{Schleiereule} let us recall the following well-known 
lemma.
\begin{lemma}
On a local coordinate $(U,z)$ of $D_i$ around a double point $p_{kl}$ 
any closed $\omega\in E^1(D_i \log P_i)$ can be written as
$\omega=\mbox{Res}_{p_{kl}}\omega\cdot\frac{dz}{z} +\psi$,
where $\psi$ is a smooth 1-form in $E^1(U)$.
\end{lemma}

\noindent
{\bf Proof of \ref{Schleiereule}:}
Let $\varphi=\sum_{i\ge 0} \omega_i + \sum_{[k<l]} K_{kl}\frac{dx}{x}
+L_{kl}\frac{dy}{y} + H_{kl} d\xi$. To prove (i), we have to show that
integrals of type (a) and (b) are well-defined, finite and do not
depend on the choice of coordinates (note that (d) is similar to (b)).
Assume that we are given coordinates $x$, $y$ and $t$ like in (a) and
write
\[
\omega_k=\mbox{Res}_{p_{kl}}\omega_k\frac{dx}{x}+\psi_k\quad\mbox{   and   }
\quad\omega_l=\mbox{Res}_{p_{kl}}\omega_l\frac{dx}{x}+\psi_l.
\]

\noindent
{\bf ad (a)}
That $\varphi$ is closed implies $K'=\frac{\partial H}{\partial u}$,
$L'=\frac{\partial H}{\partial u}$ and $\frac{\partial L}{\partial u}
=\frac{\partial K}{\partial u}$. Let $H=H_0+H_1u +\cdots+H_m u^m $.
Then $iH_i=K_{i-1}'=L_{i-1}'$ for $i=1,\ldots,m$. Since $K(0,u)$,
$K(1,u)$, $L(0,u)$, $L(1,u)$ are constants in $\cpxm$ we have
$K_i(0)=K_i(1)=L_i(0)=L_i(1)=0$ for $i\ge 1$. Moreover
$\mbox{Res}_{p_{kl}}\omega_k=-\mbox{Res}_{p_{kl}} \omega_l$
and because of $H_1=K_0'=L_0'$ we obtain $\int_0^1 H_1(\xi)d\xi=
\mbox{Res}_{p_{kl}}\omega_l$. Therefore we derived:
\[
\int\limits_0^1 H(\xi,u)d\xi = \int\limits_0^1 H_0(\xi)d\xi
+\mbox{Res}_{p_{kl}}\omega_l\cdot u.
\]
Now we may compute as follows: recall $\gamma(\tau_0)=p_{kl}$ and
choose $\varepsilon_0$ such that
$\gamma([\tau_0-\varepsilon_0;\tau_0+\varepsilon_0])\subset W_{kl}$
and $\gamma_x([\tau_0-\varepsilon_0;\tau_0])$ as well as
$\gamma_y([\tau_0;\tau_0+\varepsilon_0])$ lie within an open
set of $\cpxm$, where the logarithm is univalent. Then we find, when
we split up $\omega_k= \Res_{p_{kl}}\omega_k \frac{dx}{x}
+\psi_k$ and
$\omega_l= \Res_{p_{kl}}\omega_l\frac{dy}{y} +\psi_l$:
\begin{align*}
& \int\limits_\gamma \varphi= \lim_{\varepsilon\to 0} \left\{
 \int\limits_{\gamma^{\le \tau_0-\varepsilon}} \!\!\!\!\!\!\omega_k
 +\!\!\!\int\limits_0^1 H\left(\xi, \log\left(\gamma_x(\tau_0-\varepsilon)
 \gamma_y(\tau_0+\varepsilon)\right)\right) d\xi +
 \!\!\!\!\!\!\!\!\int\limits_{\gamma_{\ge \tau_0+\varepsilon}} \!\!\!\!\!\!
 \omega_l
 \right\}
\\
& =\int\limits_{\gamma^{\le \tau_0-\varepsilon_0}} \!\!\!\!\!\!\omega_k
 +\int\limits_{\gamma^{\le \tau_0}_{\ge \tau_0-\varepsilon_0}} \!\!\!\!\!\!\psi_k
 +\int\limits_{\gamma_{\ge \tau_0+\varepsilon_0}} \!\!\!\!\!\!\omega_l
 +\int\limits_{\gamma_{\ge \tau_0}^{\le\tau_0+\varepsilon_0}} \!\!\!\!\!\!\psi_l
 +\int\limits_0^1 H_0 (\xi) d\xi
 \\
& + \lim_{\varepsilon\to 0} \left\{
 \int\limits_{\gamma^{\le \tau_0-\varepsilon}_{\ge \tau_0-\varepsilon_0}}
 \!\!\!\!\!\!
 \Res_{p_{kl}}\omega_k \frac{dx}{x}
 +\Res_{p_kl} \omega_l \log\left(\gamma_x(\tau_0-\varepsilon)
 \gamma_y(\tau_0+\varepsilon)\right) +
 \!\!\!\!\!\!\!\!
 \int\limits_{\gamma_{\ge \tau_0+\varepsilon}^{\le \tau_+\varepsilon_0}}
 \!\!\!\!\!\!
 \Res_{p_{kl}}\omega_l \frac{dy}{y}
 \right\}.
\end{align*}

And then remark that
\begin{align*}
 & \lim_{\varepsilon\to 0} \left\{
\int\limits_{\gamma^{\le \tau_0-\varepsilon}_{\ge \tau_0+\varepsilon_0}}
\!\!\!\!\!\!
\frac{dx}{x} -\log\left(\gamma_x(\tau_0-\varepsilon)\right)+
\int\limits_{\gamma_{\ge \tau_0+\varepsilon}^{\le \tau_0 +\varepsilon_0}}
\!\!\!\!\!\!
-\frac{dy}{y} + \log\left(\gamma_y(\tau_0+\varepsilon)\right) \right\} \\
= &-\log\left(\gamma_x(\tau_0-\varepsilon_0)\right)+
\log\left(\gamma_y(\tau_0+\varepsilon_0)\right).
\end{align*}
This shows that the integral in (a) is finite.

Once we know that $\int_\gamma \varphi$ converges, we may write
$\int_\gamma \varphi$ as
\[
\int\limits_0^1 H_0 d\xi
+ \lim_{\varepsilon\to 0} \left\{
\int\limits_{\gamma^{\le \tau_0-\varepsilon}}
\!\!\!\!\!\!
\omega_k  +\mbox{Res}_{p_{kl}}\omega_l
\log\left(\gamma_x(\tau_0-\varepsilon)
\gamma_y(\tau_0+\varepsilon)\right) +
\!\!\!\!\!\!\!\!
\int\limits_{\gamma_{\ge \tau_0+\varepsilon}}
\!\!\!\!\!\!
\omega_l
\right\}
\]
\begin{multline*}
=\int\limits_0^1 H_0 d\xi
+ \lim_{\varepsilon\to 0} \Bigg\{
\int\limits_{\gamma^{\le \tau_0-\varepsilon}}
\!\!\!\!\!\!
\omega_k  +\mbox{Res}_{p_{kl}}\omega_l
\log\left(\varepsilon^2\right) +
\!\!\!\!\!\!\!\!
\int\limits_{\gamma_{\ge \tau_0+\varepsilon}}
\!\!\!\!\!\!
\omega_l \\
+\mbox{Res}_{p_{kl}}\omega_l
\log\left(\frac{\gamma_x(\tau_0-\varepsilon)}{\varepsilon}
\cdot\frac{\gamma_y(\tau_0-\varepsilon)}{\varepsilon}\right) \Bigg\}.
\end{multline*}

Observe that
\[
\lim_{\varepsilon\to 0} \left\{
\frac{\gamma_x(\tau_0-\varepsilon)}{\varepsilon}
\cdot\frac{\gamma_y(\tau_0+\varepsilon)}{\varepsilon}\right\}
=-{\dot{\gamma}_x}^{\le\tau_0} \cdot \dot{\gamma}_{y\;\ge\tau_0}
\]
and
\[
-{\dot{\gamma}_x}^{\le\tau_0} \cdot \dot{\gamma}_{y\;\ge\tau_0}
\frac{\partial}{\partial t}= \langle
-{\dot{\gamma}}_{\le\tau_0}, {\dot{\gamma}}_{\ge\tau_0}\rangle
=\frac{\partial}{\partial t}.
\]
Hence, the second summand vanishes, i.~e.~:
\[
\int\limits_\gamma \varphi
=\int\limits_0^1 H_0 d\xi
+ \lim_{\varepsilon\to 0} \left\{
\int\limits_{\gamma^{\le \tau_0-\varepsilon}}
\!\!\!\!\!\!
\omega_k  +\mbox{Res}_{p_{kl}}\omega_l
\log\left(\varepsilon^2\right) +
\!\!\!\!\!\!\!\!
\int\limits_{\gamma_{\ge \tau_0+\varepsilon}}
\!\!\!\!\!\!
\omega_l\right\}\; .
\]
This expression obviously does not depend upon the choice of the coordinates
$(x,y)$ and $t$.
If $\varphi=d(\sum_{i\ge 0} g_i +\sum_{[k<l]} P_{kl}) = 
\sum_{i\ge 0} dg_i +\sum_{[k<l]} \frac{\partial P_{kl}}{\partial u}
( \frac{dx}{x}+\frac{dy}{y} ) + P_{kl}'\;d\xi$, then
\begin{eqnarray*}
\int_\gamma \varphi & = & P_{kl}(1,u)-P_{kl}(0,u)\\
 &+ &\lim_{\varepsilon\to 0} \left\{g_k(\gamma(\tau_0-\varepsilon))
-g_k(\gamma(0)) + g_l (\gamma(1))-
g_l(\gamma(\tau_0+\varepsilon))\right\}  \\
&=& g_l(\gamma(1))-g_k(\gamma(0)).
\end{eqnarray*}

\noindent
{\bf ad (b)}
The definition does not depend on the choice of coordinates.
We just have to show that the integral converges.
\begin{align*}
\int\limits_\gamma \varphi
 &= \int\limits_{\gamma^{\le \tau_0-\varepsilon_0}} \!\!\!\!\!\!\omega_k
+\int\limits_{\gamma^{\le \tau_0}_{\ge \tau_0-\varepsilon_0}} \!\!\!\!\!\!\psi_k
+\int\limits_{\gamma_{\ge \tau_0+\varepsilon_0}} \!\!\!\!\!\!\omega_l
+\int\limits_{\gamma_{\ge \tau_0}^{\le\tau_0+\varepsilon_0}}
\!\!\!\!\!\!\psi_l\\
&+ \Res_{p_{kl}}\omega_k \lim_{\varepsilon\to 0} \left\{
\int\limits_{\gamma_{\ge \tau_0-\varepsilon_0}^{\le\tau_0-\varepsilon}}
\frac{dx}{x}
+\int\limits_{\gamma_{\ge \tau_0+\varepsilon}^{\le\tau_0+\varepsilon_0}}
\frac{dx}{x}
\right\}.
\end{align*}
And further
\begin{gather*}
\lim_{\varepsilon\to 0} \left\{
\int\limits_{\gamma_{\ge \tau_0-\varepsilon_0}^{\le\tau_0-\varepsilon}}
\frac{dx}{x}
+\int\limits_{\gamma_{\ge \tau_0+\varepsilon}^{\le\tau_0+\varepsilon_0}}
\frac{dx}{x}
\right\} \\
=\log(\gamma_x(\tau_0+\varepsilon_0))-\log(\gamma_x(\tau_0-\varepsilon_0))
+ \lim\limits_{\varepsilon\to 0}
\left(
\log(\gamma_x(\tau_0-\varepsilon))-\log(\gamma_x(\tau_0+\varepsilon))
\right).
\end{gather*}
Now is
\begin{eqnarray*}
 & & \lim_{\varepsilon\to 0} \left\{
\log(\gamma_x(\tau_0-\varepsilon))-\log(\gamma_x(\tau_0+\varepsilon))
\right\} \\
& = &  \lim_{\varepsilon\to 0} \left\{
\log\left(\frac{\gamma_x(\tau_0-\varepsilon)}{\varepsilon}\right)
-\log\left(\frac{\gamma_x(\tau_0+\varepsilon)}{\varepsilon}\right)
\right\} \\
& = &
\log\left(-{\dot{\gamma}_x}^{\le\tau_0}\right)-
\log\left(\dot{\gamma}_{x\;\ge\tau_0}\right) \\
 & = & 0 \quad\mbox{ by the definition of a path over } \vec{v}.
\end{eqnarray*}

If $\varphi$ is exact, we have
$\int_\gamma \varphi=g_k(\gamma(1))- g_k(\gamma(0))$.
\hspace*{\fill} $\Box$

\begin{remark} \rm
We see from the proof of Theorem \ref{Schleiereule} that we
have in the case of an integral  of type (a) like on page
\pageref{Kormoran}:
\[
\int\limits_\gamma \varphi
=\int\limits_0^1 H_0 d\xi
+ \lim_{\varepsilon\to 0} \left\{
\int\limits_{\gamma^{\le \tau_0-\varepsilon}}
\!\!\!\!\!\!
\omega_k  +\mbox{Res}_{p_{kl}}\omega_l
\log\left(\varepsilon^2\right) +
\!\!\!\!\!\!\!\!
\int\limits_{\gamma_{\ge \tau_0+\varepsilon}}
\!\!\!\!\!\!
\omega_l\right\} \; ,
\]
\begin{eqnarray*}
&=& \lim_{\varepsilon\to 0} \left\{
\int\limits_{\gamma^{\le \tau_0-\varepsilon}}
\!\!\!\!\!\! \omega_k
 -\mbox{Res}_{p_{kl}}\omega_k
 \log\left(\gamma_x(\tau_0-\varepsilon)\right)\right\} \\
 & + & \lim_{\varepsilon\to 0} \left\{
\int\limits_{\gamma_{\ge \tau_0+\varepsilon}}
\!\!\!\!\!\!
\omega_l + \Res_{p_{kl}}\omega_l
\log\left(\gamma_y(\tau_0+\varepsilon)\right) \right\}\\
& +& \int\limits_0^1 H_{kl}(\xi,0) d\xi .
\end{eqnarray*}
\end{remark}

\begin{remark} \rm
Observe that if one considers in (a) the closed
\[
\Xi:= K\frac{dx}{x} +L\frac{dy}{y} + H\,d\xi \in
\left({\bigwedge}^\bullet (\frac{dx}{x},\frac{dy}{y})[\log t]
\otimes E^\bullet(\Delta^1)\right)^1
\]
as element in
\[
\left({\bigwedge}^\bullet (\frac{dx}{x},\frac{dy}{y})[\log x,\log y]
\otimes E^\bullet(\Delta^1)\right)^1,
\]
which is defined in the obvious way (for the precise definition see
\ref{Tommi Engel}), by setting
\[
d\,\log x =\frac{dx}{x}\qquad \text{   and   }\qquad d \log y=\frac{dy}{y},
\]
then $\Xi$ is even exact. If $P\in\cpxm[\xi,\log x,\log y]$ is such
that $dP=\Xi$, then
\[
\int_{[0;1]} H(\xi,\log t) d\xi = P(1,\log x,\log y)- P(0,\log x,\log y).
\]
\end{remark}

Next we state that $\int_{\gamma} \varphi$ does not depend on the
choice of a path $\gamma$ within a homotopy class in
$\pi_1(Z_{\vec{v}},\vec{w})$.

\begin{theorem} \label{Hoehner}
Let $\varphi\in A^1$ be closed and let $\gamma_0$, $\gamma_1:I
\rightarrow Z_0$ be two paths over $\vec{v}$ based at $\vec{w}$,
which are nearby homotopic. Then
\[
\int_{\gamma_0} \varphi = \int_{\gamma_1} \varphi.
\]
\end{theorem}

\noindent
We will prove the corresponding theorem for iterated integrals in general
as Theorem \ref{Blaek Foes}. The Theorem \ref{Hoehner} above is then just
a special case of \ref{Blaek Foes}.

\begin{theorem}  \label{Steinkauz}
The integration of closed elements in $A^1$ along paths over $\vec{v}$
based at $\vec{w}$ defines an isomorphism
\[
H^1(A^\bullet)\stackrel{\cong}{\longrightarrow}
\Hom_\itgm\left(\vec{J}\big/\vec{J}^{\;2};\cpxm\right).
\]
\end{theorem}

\noindent
{\bf Proof:}
Since we know that both $\cpxm$-vector spaces have the same dimension,
it suffices to show that the integration map is injective. Note that 
\[
\Hom_\itgm\left(\vec{J}\big/\vec{J}^{\;2};\cpxm\right)
\hookrightarrow \Hom_\itgm\left(\vec{J};\cpxm\right)
\cong \Hom_\itgm\left(\pi_1(Z_{\vec{v}},\vec{w});\cpxm\right).
\]
Assume $\varphi=\sum_{i\ge 0} \omega_i + \sum_{[k<l]} K_{kl}\frac{dx}{x}
+L_{kl}\frac{dy}{y} + H_{kl} d\xi\in A^1$ is closed and that
\[
\int_\gamma \varphi = 0 \qquad\mbox{ for any } [\gamma]
\in\pi_1(Z_{\vec{v}},\vec{w}).
\]
We claim that $\varphi$ is exact. First, observe that all residues
of $\varphi$ are necessarily zero, as one sees by integration along a path over
$\vec{v}$ based at $\vec{w}$, which approaches a double point,
turns once around this double point and returns as it came.
Moreover all the $\omega_i$ are exact since the expression
$g_i(p):= \int_{\vec{w}}^p \varphi$ for $p\in D_i$
makes sense and $dg_i=\omega_i$. Now define 
\[
f:=\sum_{i\ge 0} g_i + \sum_{[k<l]} \int_0^\xi H(\xi, u)d\xi\in B^0
\]
and note that it satisfies the compatibility condition and that
$df=\varphi$.  \hspace*{\fill} $\Box$

The Theorem \ref{Steinkauz} allows us to define an integral lattice
$H^1(A^\bullet)_\itgm$ within $H^1(A^\bullet)$ in the following way.

\begin{definition}
Define $H^1(A^\bullet)_\itgm$ to be the lattice of those cohomology classes
which correspond to 
$\Hom_\itgm(\vec{J}\big/\vec{J}^{\;2};\itgm)$
via the integration map. $H^1(A^\bullet)_\ratm$ is then defined as
$H^1(A^\bullet)_\ratm:= 
\Hom_\itgm(\vec{J}\big/\vec{J}^{\;2};\ratm)
\subset H^1(A^\bullet).$
\end{definition}

\noindent
{\bf Remark:}
Note that in contrast to $H^1(A^\bullet)$ the lattice $H^1(A^\bullet)_\itgm$
depends on the tangent vectors $\vec{v}$ and $\vec{w}$. 

\subsection[The Hodge- and Weight Filtration]{The Hodge- and Weight Filtration}

Next we define a decreasing filtration $F^\bullet$ and an increasing filtration
$W_\bullet$ on $A^\bullet$ by defining such filtrations on $B^\bullet$.
On $\Lambda^\bullet$ we set
\begin{align*}
F^p\Lambda^\bullet &:= \bigoplus\limits_{n+m\ge p} 
{\bigwedge}^n(\frac{dx}{x},\; \frac{dy}{y})(\log t)^m \\
W_l\Lambda^\bullet & :=\bigoplus\limits_{n+2\cdot m\le l} 
{\bigwedge}^n(\frac{dx}{x},\;
\frac{dy}{y})(\log t)^m.
\end{align*}
In particular we have: $\log t \in W_2\cap F^1$. 

Then define on $E^\bullet(\Delta^1)\otimes\Lambda^\bullet$ the filtrations
$F^\bullet$ and $W_\bullet$ by:
\begin{eqnarray*}
F^p\left(E^\bullet(\Delta^1)\otimes\Lambda^\bullet\right) & := &
E^\bullet(\Delta^1)\otimes F^p\Lambda^\bullet \\
W_l\left(E^\alpha(\Delta^1)\otimes\Lambda^\beta\right) & := &
E^\alpha(\Delta^1)\otimes W_{l+\alpha}\Lambda^\beta .
\end{eqnarray*}
On $E^\bullet(D_i\log P_i)$ for $i>0$ we have the classical {\it Hodge
filtration}:
\[
F^p E^\bullet(D_i\log P_i):= \bigoplus_{n\ge p} \Gamma\left(
\Omega^n(D_i\log P_i)
\otimes_{\Omega^0(D_i)} {\mathcal E}^{0,\bullet-n}(D_i)\right)
\]
and the {\it weight filtration} $W_\bullet$. This filtration $W_\bullet$ is  
given by: 
\[
W_{l} \; E^p(D_i\log P_i)= E^l(D_i\log P_i)\otimes E^{p-l}(D_i) \quad 
\text{with } \quad 0\le l \le p.
\]
Finally, define for $i=0,\ldots,r-1$ on ${\bigwedge}^\bullet_i(\frac{dp}{p})
= {\bigwedge}^\bullet(\frac{dp}{p})$:
\begin{gather*}
W_l {\bigwedge}^\bullet(\frac{dp}{p}):=\bigoplus_{n\le l}
{\bigwedge}^n(\frac{dp}{p}) \\
F^p {\bigwedge}^\bullet(\frac{dp}{p}):=\bigoplus_{n\ge p}
{\bigwedge}^n(\frac{dp}{p}).
\end{gather*}

Putting all these filtrations together we obtain filtrations $W_\bullet$
and $F^\bullet$ on $A^\bullet$. In particular, $W_\bullet$
and $F^\bullet$ allow us to define filtrations on the  
cohomology groups in the following way.

\begin{definition}\rm
The {\it weight filtration} 
$W_\bullet$ on $H^\bullet(A^\bullet)$ is defined by:
\[
W_{l+m} H^m(A^\bullet):= \im  \left\{ H^m(W_l A^\bullet)\rightarrow 
H^m(A^\bullet)\right\}.
\]
And define the {\it Hodge filtration} 
$F^\bullet$ on $H^\bullet(A^\bullet)$ as:
\[
F^p H^m(A^\bullet):= \im  \left\{ H^m(F^p A^\bullet)\rightarrow 
H^m(A^\bullet)\right\}.
\]
On the $\ratm$-vector space $H^1(A^\bullet)_\ratm$, the weight
filtration $W_\bullet$ on $H^1(A^\bullet)$ induces the filtration:
\[
W_l H^1(A^\bullet)_\ratm := W_l H^1(A^\bullet)\cap H^1(A^\bullet)_\ratm.
\]
\end{definition}
Observe that $Gr_l^W\left(H^1(A^\bullet)_\ratm\right)=
\mbox{im}\left\{W_l H^1 (A^\bullet)_\ratm\rightarrow 
Gr_l^W H^1(A^\bullet)\right\}$.

We want to finish this subsection by proving the following theorem.

\begin{theorem} \label{Rothalsgans}
\[
( J_{\vec{v}\vec{w}}/ J_{\vec{v}\vec{w}}^2)^*
:= \left(H^1(A^\bullet)_\itgm,\;(H^1(A^\bullet)_\ratm, W_\bullet),\;
(H^1(A^\bullet), W_\bullet, F^\bullet)\right)
\]
is a $\itgm$-MHS of possible weights 0, 1, and 2.
\end{theorem}

\begin{remark} \rm
If one computes the weight spectral sequence $(A^\bullet,\;W_\bullet)$ 
one finds infinitely many $E_1$-terms, which are nonzero, each of which is
an infinite dimensional pure Hodge structure. Nevertheless, the spectral
sequence degenerates at $E_2$. In the language of 
\cite{Hain-de-Rham-homotopy} the $E_1$-term is 
the image under a certain 
De Rham functor applied to the finite dimensional $E_1$-terms of the spectral 
sequences of $( E^\bullet(D_i\log P_i),W_\bullet)$ for $i>r-1$, of
$( {\bigwedge}^\bullet_i(\frac{dp}{p}) ,W_\bullet)$ for $i\le r-1$ and of
$(\Lambda^\bullet ,W_\bullet)$ (cf.~\cite{Hain-de-Rham-homotopy}, 
the proof of (5.5.1)).

However, if one just wants to prove that $H^1(A^\bullet)$ is a MHS,
the weight spectral sequence of $(A^\bullet,\;W_\bullet)$ suggests a 
shortcut which avoids going all the way through the spectral sequences.
This is how we will prove \ref{Rothalsgans}.

Furthermore one could consider $(A^\bullet,\;W_\bullet,\; F^\bullet)$ as 
part of a mixed Hodge complex, 
which takes the integral structure (coming from integration along paths over 
$\vec{v}$ based at $\vec{w}$) into account. We omit that 
construction, since we do not need it for our further investigations.
\end{remark}

\noindent
{\bf Proof:}
Since the category of MHSs is abelian, we can prove Theorem \ref{Rothalsgans}
if we show that $W_2 H^1(A^\bullet)=H^1(A^\bullet)$ and that there are 
exact sequences
\begin{enumerate} 
\item 
$ 0 \rightarrow Gr_2^W H^1(A^\bullet)\stackrel{\alpha}{\rightarrow}
\bigoplus\limits_{[k<l]} H^0(p_{kl})(-1) \stackrel{\beta}{\rightarrow} 
\bigoplus\limits_{i>r-1} H^2(D_i) \rightarrow 0 $,
\item 
$ 0\rightarrow Gr_1^W H^1(A^\bullet) \stackrel{\cong}{\rightarrow} 
\bigoplus\limits_{i>r-1} H^1(D_i)  \rightarrow  0 $,
\item 
$\bigoplus\limits_{i\ge 0} H^0(D_i) \stackrel{\gamma}{\rightarrow}
\bigoplus\limits_{[k<l]} H^0(p_{kl}) \stackrel{\delta}{\rightarrow} 
Gr_0^W H^1(A^\bullet) \rightarrow 0 $,
\end{enumerate}
where all the maps are defined over $\ratm$ and are strict with 
respect to the respective Hodge filtrations. $Gr_l^W H^1(A^\bullet)$ is then
a Hodge structure of weight $l$ since all the other terms are Hodge structures
of appropriate weights. 

Note that if we are given a $[\varphi]=[\sum_{i\ge 0} \omega_i 
+\sum_{[k<l]} K_{kl}\frac{dx}{x}+ L_{kl}\frac{dy}{y} +H_{kl} d\xi]
\in H^1(A^\bullet)$, then
we can always write with complex numbers $\rho_{kl}$, $c_{kl}\in\cpxm$
\[
[\varphi]=[\sum_{i\ge 0} \omega_i 
+\sum_{[k<l]} \frac{1}{2\pi i} \rho_{kl} \,\Theta  + c_{kl} d\xi],
\]
where $\Theta:= (1-\xi)\frac{dx}{x}-\xi \frac{dy}{y} - \log t \, d\xi \in
W_1\left[\left(E^\bullet(\Delta^1)\otimes\Lambda^\bullet\right)^1\right]$.

This shows that $W_2H^1(A^\bullet)=H^1(A^\bullet)$ and the following two
equivalences:
\[
[\varphi]\in W_1H^1(A^\bullet)\Leftrightarrow \;\forall [k<l]:\;\rho_{kl}=0 
\mbox{ and }[\varphi]\in W_0H^1(A^\bullet)\Leftrightarrow 
\;\forall i\ge 0:\;\omega_i=0.
\]

\noindent
{\bf ad (i):}
We identify $H^0(p_{kl})(-1)$ with $\cpxm(-1)$ such that 
$H^0(p_{kl})(-1)_\itgm$ 
becomes $\frac{1}{2\pi i}\itgm\subset \cpxm$ and $Gr_F^1 \cpxm =\cpxm$.
Then define the map $\alpha$ by 
\[
\alpha\left([\sum_{i\ge 0} \omega_i 
+\sum_{[k<l]} \frac{1}{2\pi i} \rho_{kl} \,\Theta]\right):=
\sum_{[k<l]} \frac{1}{2\pi i} \rho_{kl}.
\]
If $[\sum_{i\ge 0} \omega_i +\sum_{[k<l]} \frac{1}{2\pi i} \rho_{kl} \,\Theta]$
represents an element in $H^1(A^\bullet)_\itgm$ then it is mapped to 
$H^0(p_{kl})(-1)_\itgm$ under $\alpha$, since all the numbers $\rho_{kl}$
have to be integers (they are the values of an integral along a path around 
a double point).

Now define the map $\beta$ in the following way. Given an element 
$\sum_{[k<l]} \frac{1}{2\pi i} \rho_{kl} \in
\bigoplus_{[k<l]} H^0(p_{kl})(-1)$, choose a $\sum_{i\ge 0} \omega_i 
+\sum_{[k<l]} \frac{1}{2\pi i} \rho_{kl} \,\Theta\in A^1$ such that 
for $i>r-1$ all the $d\omega_i$ have no singularities\footnote{This 
is always possible. 
It can be done similarly to the following: define a bump-function 
$b:\cpxm\rightarrow [0;1]$ around $0\in\cpxm$ such that $b(z)\equiv 1$ in a 
small neighbourhood of $0\in\cpxm$. 
Then $b(z)\rho\frac{dz}{z}$ has residue $\rho$ and 
$d(b(z)\rho\frac{dz}{z})\in E^2(\cpxm)$.}, 
i.~e.~$d\omega_i\in E^2(D_i)$. Then define:
\[
\beta\left(\sum_{[k<l]} \frac{1}{2\pi i} \rho_{kl}\right):=
\sum_{i > r-1} [d\omega_i].
\]
Note that $\alpha$ and $\beta$ are well-defined and that they yield a short
exact sequence (i), where all the maps are 
strict with respect to the filtrations $F^\bullet$.

\noindent
{\bf ad (ii):}
The isomorphism sends $[\sum_{i\ge 0} \omega_i]\in W_1 H^1(A^\bullet)$ to 
$\sum_{i\ge r-1} [\omega_i]\in \bigoplus\limits_{i>0} H^1(D_i)$. All the 
remaining verifications are obvious.

\noindent
{\bf ad (iii):}
The maps $\gamma$ and $\delta$ are defined in the following way:
\[
\begin{array}{cccccc}
\bigoplus\limits_{i\ge 0} H^0(D_i) & \stackrel{\gamma}{\rightarrow} &
\bigoplus\limits_{[k<l]} H^0(p_{kl}) & \stackrel{\delta}{\rightarrow} &
Gr_0^W H^1(A^\bullet) & \rightarrow 0 \\
\sum_{i\ge 0} g_i & \mapsto & \sum_{[k<l]}(g_l-g_k) & & & \\
& & \sum_{[k<l]} c_{kl} & \mapsto & [\sum_{[k<l]} c_{kl}\;d\xi] & . 
\end{array}
\]
Note that the kernel of $\gamma$ has dimension 1.
Also this sequence is exact, defined over $\ratm$ and respects the 
filtrations $F^\bullet$. This accomplishes the proof. \hspace*{\fill} $\Box$

\begin{remark} \label{sleutel}
Using the short exact sequences (i), (ii) and (iii) from the proof of
\ref{Rothalsgans} we may derive that, if $g(D_i)$ denotes the genus of
$D_i$ for $i>r-1$:
\begin{eqnarray*}
\dim_\cpxm\; Gr_2^W H^1(A^\bullet) & = & (\sum_{[k<l]} 1 )
- ( \sum_{i>r-1} 1 ), \\
\dim_\cpxm\; Gr_1^W H^1(A^\bullet) & = & \sum_{i> r-1} 2 g(D_i) , \\
\dim_\cpxm\; Gr_0^W H^1(A^\bullet) & = & (\sum_{[k<l]} 1 )
- ( \sum_{i\ge 0} 1 ) + 1.
\end{eqnarray*}
\end{remark}

\subsection[The Variation is a Nilpotent Orbit] {The Variation of
Mixed Hodge Structure is a Nilpotent Orbit} \label{Maibaum}

Now it is not difficult anymore to investigate the monodromy
of the local system of nearby fundamental groups and to show
that the variation of Hodge structure 
\[
\left({J}_{\vec{v},\vec{w}}
\big/ J_{\vec{v},\vec{w}}^2\right)^*_{\vec{v}\in (T_0\Delta)^*}
\] 
is a nilpotent orbit of Hodge structure. 

\begin{remark} \rm 
The term {\it nilpotent orbit of mixed Hodge structure} is a generalization 
(for one variable) of the definition of a nilpotent orbit of MHS given
in \cite{Schmid} and \cite{Cattani-Kaplan} (Definition 3.1); specifically,
\cite{Cattani-Kaplan}: (3.1) (iv) is not satisfied. A nilpotent orbit of
MHS will be defined properly in Proposition \& Definition \ref{Dagmar}.
\end{remark}

First of all, we define a nilpotent chain morphism
\[
N:A^\bullet\rightarrow A^\bullet
\]
by defining it in the following way on $B^\bullet$. 
Let $N$ be the zero
map on $\bigoplus_{i\le r-1}{\bigwedge}^\bullet_i (\frac{dp}{p})\oplus
\bigoplus_{i>r-1} 
E^\bullet (D_i \log P_i)$ and define it on 
$E^\bullet(\Delta^1)\otimes{\Lambda}^\bullet$ as 
\[
\begin{array}{cccc}
N: & E^\bullet(\Delta^1)\otimes{\Lambda}^\bullet & \rightarrow & 
E^\bullet(\Delta^1)\otimes{\Lambda}^\bullet \\*
   & \Xi  & \mapsto & \frac{d}{d(-u)} \Xi = - \frac{d}{du} \Xi \;.
\end{array}
\]
Note that $N$ respects the compatibility relations and hence defines 
a chain map $N:A^\bullet\rightarrow A^\bullet$. Concerning
the filtrations, $N$ induces maps
\[
N:W_{l+1} A^\bullet \rightarrow W_{l-1} A^\bullet \quad
\text{ and }\quad N:F^p A^\bullet \rightarrow F^{p-1} A^\bullet.
\]

In particular this chain map
induces a map $N:H^1(A^\bullet)\rightarrow H^1(A^\bullet)$ (note that 
we are committing `abuse of notation' here).

\begin{proposition} \label{Bitters"usser Nachtschatten}
The filtration $W_\bullet$ on $H^1(A^\bullet)$ and 
$N:H^1(A^\bullet)\rightarrow H^1(A^\bullet)$ satis\-fy $N^2=0$ and
\begin{enumerate}
\item
$W_{-1} H^1(A^\bullet) = 0$
\item
$W_{0}\;\, H^1(A^\bullet) = \im 
\left\{N:H^1(A^\bullet)\rightarrow H^1(A^\bullet) \right\}$
\item
$W_{1}\;\, H^1(A^\bullet) = \ker 
\left\{N:H^1(A^\bullet)\rightarrow H^1(A^\bullet) \right\}$
\item
$W_{2}\;\, H^1(A^\bullet) = H^1(A^\bullet)$.
\end{enumerate}
\end{proposition}

\noindent
{\bf Proof:} 
First, recall that we can always write a 
$[\varphi]\in H^1(A^\bullet)$ as 
\[
[\varphi]=[\sum_{i\ge 0} \omega_i 
+\sum_{[k<l]} \frac{1}{2\pi i} \rho_{kl} \,\Theta  + c_{kl} d\xi],
\]
where the $\omega_i$ are closed forms in $E^1(D_i \log P)$, the $c_{kl}$ 
are complex numbers and 
$\Theta:= (1-\xi)\frac{dx}{x}-\xi \frac{dy}{y} - \log t \, d\xi \in
W_1\left(E^\bullet\otimes\Lambda^\bullet\right)^1$.
This shows $N^2=0$ and $W_{-1} H^1(A^\bullet) = 0$ as well as
$W_{2}\; H^1(A^\bullet) = H^1(A^\bullet)$. 

\noindent 
{\bf ad(ii):} Since the sum of the residues for any closed form $\omega_i$
is zero, we have
\begin{align*}
W_{0}\; H^1(A^\bullet) &= \left\{ [\sum_{[k<l]} c_{kl} d\xi]\right\},\\
\im N &= \left\{ [\sum_{[k<l]} \frac{1}{2\pi i} \rho_{kl} d\xi]\bigg|
\sum_{[i<l]} \frac{1}{2\pi i} \rho_{il} - \sum_{[k<i]} 
\frac{1}{2\pi i} \rho_{ki} =0 \text{ for all }i\right\}.
\end{align*}
Obviously: $\im \left\{N:H^1(A^\bullet)\rightarrow H^1(A^\bullet) \right\}
\subseteq W_{0}\; H^1(A^\bullet)$. Now let $[\sum_{[k<l]} c_{kl} d\xi]$
be in $W_{0}\; H^1(A^\bullet)$ and let $m:=\sum_{[k<l]} 1$.  
Define complex numbers
\[
g_i:= \frac{1}{m}\left( \sum_{[k<i]} c_{ki} - \sum_{[i<l]} c_{il} \right).
\]
Note that $\sum_{i\ge 0} g_i =0$. It is easy to check that
\[
[\sum_{[k<l]} c_{kl} d\xi]=[\sum_{[k<l]} (c_{kl}-(g_l-g_k)\big) d\xi] 
\in \im N.
\]

\noindent 
{\bf ad(iii):} Observe that on the one hand
\begin{align*}
W_{1}\; H^1(A^\bullet) & = \left\{ [\sum_{i> 0} \omega_i 
+\sum_{[k<l]} c_{kl} d\xi]\right\}\\
\intertext{and on the other hand}
\ker N & = \left\{ [\sum_{i\ge 0} \omega_i 
+\sum_{[k<l]} \frac{1}{2\pi i} \rho_{kl} \,\Theta  + c_{kl} d\xi]\;\bigg|\;
[\sum_{[k<l]} \frac{1}{2\pi i} \rho_{kl} \; d\xi]=0 \right\}.
\end{align*}

Hence, we have to show: $[\sum_{[k<l]} \frac{1}{2\pi i} \rho_{kl} \; d\xi]=0 
\in H^1(A^\bullet) \Rightarrow \text{ all } \rho_{kl}=0$. For 
simplicity define $R_{kl}:=\frac{1}{2\pi i} \rho_{kl}$.
Let $g_i\in\cpxm$ be complex numbers such that 
$d \left(\sum_{i\ge 0} g_i +\sum_{[k<l]} \{g_k+\xi(g_l-g_k)\}\right)=
\sum_{[k<l]} R_{kl} \; d\xi$, 
i~e.~$(g_l-g_k)=R_{kl}$  $\forall [k<l]$. We show that the real parts,
$\Re(R_{kl})$, of all $R_{kl}$ are zero and leave it to the reader to show
that the imaginary parts are zero too. 

Let $i_0$ be an index such that $\Re(g_{i_0})=\max_{i} \Re(g_{i})$. 
As a consequence, all the real parts of residues at $\omega_{i_0}$,
that is the numbers $\Re (\Res_{p_{kl}}\omega_{i_0})$ for $\{k,l\}\ni i_0$
and $[k<l]$, are non negative ($\Res_{p_{i_0l}}\omega_{i_0}=R_{i_0l}
=g_{i_0}-g_l$ and $\Res_{p_{ki_0}}\omega_{i_0} = -R_{ki_0}
=g_{i_0}-g_k$). Now
the sum of the residues of the form $\omega_{i_0}$ has to be zero. 
This can only happen if for all the indices
$i$ with $D_i\cap D_{i_0}\neq \emptyset$ holds: $g_i=g_{i_0}$. 
Since the system od curves $D_i$ is connected, this shows 
$\Re(g_i)=\Re(g_{i_0})$ for all $i$ and hence $\Re(R_{kl})=0$ for all
$[k<l]$. \qed

\begin{remark} \rm
Proposition \ref{Bitters"usser Nachtschatten} shows that $W_\bullet$ is
{\it the weight filtration of $N$ relative to $W_\bullet$} as defined
in \cite{Steenbrink-Zucker} (2.5). The proof of Proposition 
\ref{Bitters"usser Nachtschatten} is trivial, when one uses 
Proposition (2.14) in \cite{Steenbrink-Zucker}.
\end{remark}

Let $T\in\Aut\left(H^1(A^\bullet)\right)$ be the monodromy of the local system
\[
\left\{\left(J_{\vec{v},\vec{w}}\big/
J^2_{\vec{v},\vec{w}}\right)^*_\itgm\right\}_{\vec{v}\in (T_0\Delta)^*}.
\]
Recall that by means of Theorem \ref{Dr. Yoghurt} the automorphism 
$T$ can be identified with
the monodromy of the local system $\{H^1(Z_t;\itgm)\}_{t\in \Delta^*}$.

A basic example of {\it a variation of mixed Hodge structure} is 
given by the following definition. For a complex number 
$\lambda\in\cpxm^*$ and a nilpotent vectorspace endomorphism $A$ define 
the endomorphism
$\lambda^{A}$ to be $e^{\log \lambda A}= {\pmb 1}+\log \lambda A+
\frac{1}{2}\log^2 \lambda A^2 + \cdots$. 
\begin{proposition+definition} \label{Dagmar} \rm
Let $H$ be a $\itgm$-MHS and $N:H\rightarrow H$ be a nilpotent
endomorphism such that 
$N(F^pH)\subset F^{p-1}H$ and $N(W_{l+1}H)\subset W_{l-1}H$.
Then for any $\lambda\in\cpxm^*$, the triple
\[
H_\lambda:=\left( \lambda^{-N} H_{\itgm}, (H_\ratm,W_\bullet),
(H_\cpxm,W_\bullet,F^\bullet)\right)
\]
is a MHS and the family of MHSs 
$\{H_\lambda\}_{\lambda\in\cpxm^*}$
is called {\it nilpotent orbit of mixed Hodge structure}. \qed
\end{proposition+definition}

In the next theorem we describe, how the lattice 
$(J_{\vec{v},\vec{w}}\big/J^2_{\vec{v},\vec{w}})^*_\itgm
= H^1(A^\bullet)_\itgm\subset H^1(A^\bullet)$
moves, when we move $\vec{v}$ and $\vec{w}$. 

\begin{theorem} \label{Noli me tangere} 
$e^{-2\pi i\; N}= T$ and for $\lambda$, $\mu\in\cpxm$ holds:
\[
\bigg(J_{\lambda\vec{v},\mu\vec{w}}\big/
J^2_{\lambda\vec{v},\mu\vec{w}}\bigg)^*_\itgm=\lambda^{-N}
\bigg(J_{\vec{v},\vec{w}}\big/J^2_{\vec{v},\vec{w}}\bigg)^*_\itgm
\subset H^1(A^\bullet).
\]
\end{theorem}

\noindent
{\bf Proof:} 
Note that $N^2=0$ implies $\lambda^{N}={\pmb 1}+\log \lambda N$.
Let $\eta=\lambda\times\mu:[0;1] \rightarrow 
(\cpxm\setminus\reellm^{\le 0})\times(\cpxm\setminus\reellm^{\le 0})$ be a 
path with $\eta(0)=(1,1)$ 
and let $H:[0;1]\times [0;1] \rightarrow Z_0$ be a homotopy such that 
$H(\cdot,s)$ is a path over $\lambda(s)\vec{v}$ based at $\mu(s)\vec{w}$.
Define $\gamma_{\vec{v},\vec{w}}:=H(\cdot,0)$.

Using the definitions (a), (b), (c) and (d) on page \pageref{Kormoran},
it is easy to compute for a closed $\varphi\in A^1$:
\[
\int_{H(\cdot,s)} \varphi = \int_{\gamma_{\vec{v},\vec{w}}} \varphi
+ \log \lambda(s) \int_{\gamma_{\vec{v},\vec{w}}} N \varphi
= \int_{\gamma_{\vec{v},\vec{w}}} \lambda(s)^N \varphi.
\]
Suppose $\int_{\gamma_{\vec{v},\vec{w}}} \varphi$ is an integer, then
$\int_{H(\cdot,s)} \lambda(s)^{-N} \varphi$ is an integer too. \qed

The following statement is a direct consequence of Theorem 
\ref{Bitters"usser Nachtschatten} and Theorem \ref{Noli me tangere}.
\begin{theorem}
The family of mixed Hodge structures
\[
\left\{\left(J_{\vec{v},\vec{w}}\big/J^2_{\vec{v},\vec{w}}\right)^* 
\right\}_{\vec{v}\in (T_0\Delta)^*}
\]
is a nilpotent orbit of mixed Hodge structure. \qed
\end{theorem}

Consider the graph $\Gamma$, whose vertices are the components $D_i$,
where two vertices $D_k$ and $D_l$ are connected by an edge if
and only if $D_k\cap D_l \neq \emptyset$. Then we have the following
corollary of Theorem \ref{Bitters"usser Nachtschatten} and
Theorem \ref{Noli me tangere}. \label{Herr Graf}

\begin{corollary} \label{Gr"unschnabel}
$\Gamma$ is a tree if and only if $T=id$. 
\end{corollary}

\noindent
{\bf Proof:}
By Theorem \ref{Bitters"usser Nachtschatten} and Theorem \ref{Noli me tangere}
$T=id$ is equivalent to $N=0$ or $\ker N = W_1H^1(A^\bullet) 
= W_2H^1(A^\bullet) = H^1(A^\bullet)$. We observed in Remark \ref{sleutel}
that holds:
\[
\dim_\cpxm\; Gr_2^W H^1(A^\bullet) = (\sum_{[k<l]} 1 )
- ( \sum_{i\ge 0} 1 ) +1 .
\]
Since $\Gamma$ is connected, we have: $\Gamma$ is a tree $\Leftrightarrow$
$\#\text{edges} = \#\text{vertices} -1 $. \qed

\subsection{Appendix} \label{Menno}

We still have to give the proof of Theorem \ref{Graugans}.

\noindent
{\bf Proof of \ref{Graugans}:} The only closed elements in $A^0$ are the constant compatible
functions. Therefore:
$H^0(A^\bullet)\cong H^0(Z_t;\cpxm)\cong\cpxm$.

We make use of the fact
that the map $\stackrel{\circ}{h}: \stackrel{\circ}{Z} \rightarrow
\stackrel{\circ}{\Delta}$ in \ref{Kolk} is a locally trivial 
fibration and prove:
$H^\bullet(A^\bullet)\cong H^\bullet(\stackrel{\circ}{Z}_{(0,1)};\cpxm)$.

Recall that we can build up $\stackrel{\circ}{Z}_{(0,1)}$ in the following
way. For each $p_{kl}$ with $[k<l]$ we have the open disks
$U_{kl}^k\subset D_k$ and $U_{kl}^l\subset D_l$ around the double
point $p_{kl}$. Now define
\[
B_{kl}:= \left(\stackrel{\circ}{U}_{kl}^k {\scriptscriptstyle\coprod}
\stackrel{\circ}{U}_{kl}^l\right) \bigg/ \sim\qquad ,
\]
where ``$\sim$" denotes the equivalence relation:
\[
\left(R_1, e^{2\pi i \phi_1}\right) \sim 
\left(R_2, e^{2\pi i \phi_2}\right)
\; :\Leftrightarrow \; R_1=R_2=0 \text{ and } \phi_1=-\phi_2.
\]

Then ${U_{kl}^k}^* {\scriptscriptstyle\coprod} {U_{kl}^l}^*$ can be
considered as subspace of
$B_{kl}$ and as subspace of $\coprod_{i\ge 0} D_i^*=
\coprod_{i\ge 0} D_i\setminus P_i$. We construct the fiber
$\stackrel{\circ}{Z}_{(0,1)}$ by glueing $\coprod_{[k<l]} B_{kl}$
and $\coprod_{i\ge 0} D_i^*$ together along the
${U_{kl}^k}^* {\scriptscriptstyle\coprod} {U_{kl}^l}^*$'s. Note that
$(\coprod_{i\ge 0} D_i^*,\, \coprod_{[k<l]} B_{kl})$
is an excisive couple (cf.~\cite{Spanier} p.~188) for the space
$\stackrel{\circ}{Z}_{(0,1)}$, which is to say that we may apply
Mayer-Vietoris. Hence, there is an exact sequence
\begin{multline} \label{Blaessgans}
0\rightarrow H^0(\stackrel{\circ}{Z}_{(0,1)})\rightarrow
\bigoplus\limits_{i\ge 0} H^0(D_i^*)\oplus\bigoplus\limits_{[k<l]}H^0(B_{kl})
\rightarrow\bigoplus\limits_{[k<l]}H^0({U_{kl}^k}^*) \oplus H^0({U_{kl}^k}^*)\\
\rightarrow H^1(\stackrel{\circ}{Z}_{(0,1)})\rightarrow
\bigoplus\limits_{i\ge 0} H^1(D_i^*)\oplus\bigoplus\limits_{[k<l]}H^1(B_{kl})
\rightarrow\bigoplus\limits_{[k<l]}H^1({U_{kl}^k}^*) \oplus H^1({U_{kl}^k}^*)
\rightarrow 0.
\end{multline}

Now we compare the last terms in (\ref{Blaessgans}) with the last terms in
(\ref{Saatgans}). That is, we define isomorphisms $\varphi$, $\psi$ such
that
\[
{\scriptstyle  \begin{array}{ccccccc}
 H^1(B^\bullet) & \!\!\!\!\!\stackrel{\Phi}{\rightarrow}&\!\!\!\!\! H^1(C^\bullet)
 & \!\!\!\!\!\rightarrow  & \!\!\!\!\! H^2(A^\bullet) & \!\!\!\!\!\rightarrow &
\!\!\!\!\! 0 \\
 & & & & & & \\
\cong\big\downarrow\varphi & &\cong\big\downarrow\psi & &\big\downarrow & & \\
 & & & & & & \\
\bigoplus\limits_{i\ge 0} H^1(D_i^*)\oplus\bigoplus\limits_{[k<l]}H^1(B_{kl}) &
\!\!\!\!\!\rightarrow & \!\!\!\!\!\bigoplus\limits_{[k<l]}H^1({U_{kl}^k}^*) \oplus
H^1({U_{kl}^k}^*) & \!\!\!\!\!\rightarrow & 0 \end{array}
}\]
commutes. It is then a consequence of the 5-lemma that $H^2(A^\bullet)=0$.

Denote by $\{\frac{dx}{x}\}$ (resp.~$\{\frac{dy}{y}\}$) the cohomology class
in $H^1(B_{kl})$, which has value $2\pi i$ on the positive generator of
$H_1({U_{kl}^k}^*)$ (resp.~$H_1({U_{kl}^l}^*)$) in $H_1(B_{kl})$.
Then $\{\frac{dx}{x}\}=-\{\frac{dy}{y}\}$.
We define an isomorphism
\[
{\scriptstyle
\begin{array}{rccc}
\varphi:& \!\!\! H^1(B^\bullet) & \!\!\!\!\!\rightarrow &
\!\!\!\!\!\!\!\bigoplus\limits_{i\ge 0} H^1(D_i^*)
\oplus\bigoplus\limits_{[k<l]}H^1(B_{kl}) \\
 & \!\!\! [\sum\limits_{i\ge 0} \omega_i +
\sum\limits_{[k<l]} K\frac{dx}{x}+L\frac{dy}{y}+H d\xi] & \!\!\!\!\!\mapsto &
\!\!\!\!\!\!\!\sum\limits_{i\ge 0}[\omega_i]
+\sum\limits_{[k<l]} K\{\frac{dx}{x}\}+L\{\frac{dy}{y}.\}
\end{array}}\]
(Note that for closed $\sum_{i\ge 0} \omega_i +
\sum_{[k<l]} K\frac{dx}{x}+L\frac{dy}{y}+H d\xi$, the polynomial
$K(\xi,u)-L(\xi,u)$ is constant, i.~e.~$K(\xi,u)-L(\xi,u)=K(0,0)-L(0,0)
=K(1,0)-L(1,0)$. Therefore,
$K\{\frac{dx}{x}\}+L\{\frac{dy}{y}\}=(K-L)\{\frac{dx}{x}\}$ makes sense.)

Define moreover the isomorphism
\[ {\scriptstyle
\begin{array}{rccc}
\psi: &\!\!\! H^1(C^\bullet) &\!\!\! \rightarrow & \!\!\!\!\!\!
\bigoplus\limits_{[k<l]}H^1({U_{kl}^k}^*) \oplus H^1({U_{kl}^k}^*)\\
 & \!\!\!\!\!\!\!\!\!\!\!\!\! \sum\limits_{[k<l]}(a_k\frac{dx}{x}
+b_k \frac{dy}{y})+(a_l\frac{dx}{x}+b_l\frac{dy}{y}) & \!\!\!\mapsto &
\!\!\!\!\!\!\sum\limits_{[k<l]} (a_k-b_k)\{\frac{dx}{x}\}
+(b_l-a_l)\{\frac{dy}{y}\}
\end{array}}\]
and note that the above diagram commutes.

Finally, if we use the fact that the Euler characteristics
in (\ref{Saatgans}) on page \pageref{Saatgans} and (\ref{Blaessgans}) on
page \pageref{Blaessgans} are zero, we find:
\begin{eqnarray*}
\dim_\cpxm H^1(A^\bullet)=\dim_\cpxm H^1(\stackrel{\circ}{Z}_{(0,1)})
&=& 1+\sum_{i\ge 0}\left(\dim_\cpxm H^1(D_i^*)-1 \right)\\
&=& r+1+\sum_{i\ge 0}\left(\dim_\cpxm H^1(D_i)-2 \right) +\sum_{[k<l]} 2. 
\quad\Box
\end{eqnarray*}


\newcommand{\eint}[1]{\ensuremath{\underset{#1}{\varepsilon\!\!\!\!\int}}}
\newcommand{\eintt}[1]{\ensuremath{\varepsilon\!\!\!\int_{#1}}}
\newcommand{\chen}{I=\sum_{|J|\le s} a_J\;
\varphi_{j_1}\otimes \cdots \otimes\varphi_{j_r}\in \bigoplus_{r=1}^s
\bigotimes^r A^1}
\newcommand{\chenkl}{I_{kl}=\sum_{|J|\le s} a_J\;
\varphi_{j_1}\otimes \cdots \otimes\varphi_{j_r}\in \bigoplus_{r=1}^s
\bigotimes^r A^1_{kl}}
\newcommand{\epsT}{{\textstyle {}^{\varepsilon}}\!\!\!
{\text{\rm \Large T}}}
\newcommand{\T}{\text{\rm \Large T}}

\section[Iterated Integrals on $\pi_1(Z_{\vec{v}},\vec{w})$]
{Iterated Integrals on the Nearby Fundamental Group} \label{Bl"omscher}

Let $I=\sum_{|J|\le s} a_J\; \varphi_{j_1}\otimes \cdots \otimes
\varphi_{j_r}\in \bigoplus_{r=1}^s \bigotimes^r A^1$.
As a generalization of the notion of closedness of an element in $A^1$
the notion of Chen-closedness for an element $I$ like above 
will be given. Under the assumption that such an $I$ is Chen-closed
we will define an {\it iterated integral of $I$ along a path $\gamma$ over
$\vec{v}$ (based at $\vec{w}$)}.

First, for a sufficiently small $\varepsilon> 0$ and
arbitrary $\varphi_1\otimes \cdots \otimes \varphi_r\in \bigotimes^r
A^1$, we define what we call the {\it $\varepsilon$-iterated integral}
of $\varphi_1\otimes \cdots \otimes \varphi_r$ along a path $\gamma$
over $\vec{v}$ (based at $\vec{w}$) and will denote it by:
\[
\eint{\gamma} \varphi_1 \cdots \varphi_r.
\]
This definition depends on the choice of coordinates around the double 
points (and on $\varepsilon$). Moreover, in general the expression diverges
when passing $\varepsilon\to 0$.
But we can prove that for a {\it Chen-closed element} $I$ 
and for a path $\gamma$ over $\vec{v}$ based at $\vec{w}$ the limit
\[
\lim_{\varepsilon\to 0}\eint{\gamma} I
=\lim_{\varepsilon\to 0} \sum_{|J|\le s} a_J\;
\eint{\gamma} \varphi_{j_1} \cdots \varphi_{j_r}
\]
exists and neither depends on the choice of coordinates around the double
points nor on the choice of $\gamma$ within a nearby homotopy class.
This will be done in Section \ref{Placenta}.
The main tool for the proofs is a dga $A^\bullet_{kl}$,
which we consider to be $A^\bullet$ localized in a sector at a double 
point $p_{kl}$.

These iterated integrals along paths over $\vec{v}$ based at $\vec{w}$
will give a way to describe
$\Hom_\itgm({\vec{J}}\big/ {\vec{J}^{s+1}},\cpxm)$
in terms of the dga $A^\bullet$, similar to Chen's $\pi_1$-De Rham theorem 
(cf.~(4.4) in \cite{Hain-the-geometry}).
Finally, this leads to the definition of a MHS on
$({\vec{J}}\big/ {\vec{J}^{s+1}})^*$ such that
the short exact sequence dual to \eqref{Hausgans} becomes a short exact 
sequence of MHSs. For instance the reference \cite{Hain-the-geometry} 
contains the necessary
background on iterated (line) integrals on manifolds.

\subsection{Review on Differential Graded Algebras}

\newcommand{\A}{{\mathcal A}^\bullet}
\newcommand{\Ae}{{\mathcal A}^1}
\newcommand{\bA}{{\bar{\mathcal A}}^\bullet}
\newcommand{\bAe}{{\bar{\mathcal A}}^1}
\newcommand{\R}{{\mathcal R}}
\newcommand{\K}{{\mathcal K}}

Here we are going to introduce a little of Chen's {\it Reduced Bar
Construction}; only so much as we will need in the sequel. For a broader introduction 
to this subject we refer to \cite{Chen-RBC} or to 
\cite{Hain-de-Rham-homotopy}.

Let $\A$ be a connected dga over $\korpm=\reellm$ or $\cpxm$ 
with augmentation map $\epsilon:\A\rightarrow 
\korpm$ and let $\R(\A,\epsilon)$ be the sub-vector space of
\[
\bigoplus\limits_{r=0}^\infty {\bigotimes}^r \A
\]
generated by elements $R_i(u,f)$, $i=1,\ldots,r$, which are defined for 
$f\in {\mathcal A}^0$ and 
$u=\varphi_1\otimes\cdots\otimes\varphi_r\in{\bigotimes}^r {\A}$ by:
\begin{align} \label{Johannes}
R_1(u,f) := & df\otimes\varphi_2\otimes\cdots\otimes\varphi_r  -
(f-\epsilon(f))\varphi_2\otimes\varphi_3\otimes\cdots\otimes\varphi_r, \\
R_i(u,f) := 
&\varphi_1\otimes\cdots\otimes\varphi_{i-1}\otimes df\otimes\varphi_{i+1}
\otimes\cdots\otimes\varphi_r \nonumber \\
 +&\varphi_1\otimes\cdots\otimes f\varphi_{i-1}\otimes\varphi_{i+1}
\otimes\cdots\otimes\varphi_r \nonumber \\
 - &\varphi_1\otimes\cdots
\otimes\varphi_{i-1}\otimes f\varphi_{i+1}
\otimes\cdots\otimes\varphi_r  \quad \text{ for } 1< i < r ,\label{Wilhelm}\\
\label{Sigrid}
R_r(u,f) :=& \varphi_1\otimes\cdots\otimes\varphi_{r-1}\otimes df +
\varphi_1\otimes\cdots\otimes (f-\epsilon(f))\varphi_{r-1}.
\end{align}

We define two maps\footnote{Define $\bigotimes^0 \A $ to be $\korpm$ and
$d_C:\Ae\rightarrow \korpm$ to be the zero map.}, the {\it internal} and 
the {\it combinatorial differential},
\[
d_I,\; d_C: \bigoplus\limits_{r=0}^\infty {\bigotimes}^r \Ae
\longrightarrow
\bigoplus\limits_{r=0}^\infty {\bigotimes}^r \A
\]
by giving their values on elements of the form
$\varphi_1\otimes\cdots\otimes\varphi_s$ and extending $\korpm$-linearly:
\begin{eqnarray*}
d_I(\varphi_1\otimes\cdots\otimes\varphi_s) & := &
\sum\limits_{i=1}^s\varphi_1\otimes\cdots\otimes\varphi_{i-1}
\otimes d\varphi_i\otimes\varphi_{i+1}\otimes\cdots\otimes\varphi_s \, ,\\
d_C(\varphi_1\otimes\cdots\otimes\varphi_s) & := &
\sum\limits_{i=1}^{s-1}\varphi_1\otimes\cdots\otimes\varphi_{i-1}
\otimes\varphi_{i}\wedge\varphi_{i+1}\otimes\varphi_{i+2}\otimes
\cdots\otimes\varphi_s.
\end{eqnarray*}
Here the maps $d_I$ and $d_C$ shift the degrees in the following way:
\[
d_I: {\bigotimes}^r \Ae \longrightarrow
\left({\bigotimes}^r \left[\Ae\oplus {\mathcal A}^2 \right]\right)^{r+1},
\;
d_C: {\bigotimes}^r \Ae \longrightarrow
\left({\bigotimes}^{r-1} \left[\Ae\oplus {\mathcal A}^2 \right]\right)^{r}_.
\]
Define 
\begin{align*}
\R^s_1(\A,\epsilon) &:= \left(\bigoplus\nolimits_{r=0}^s
\bigotimes\nolimits^r\Ae \right)\cap 
\R(\A,\epsilon)
\end{align*}
and note that $(d_I+d_C)(\R^s_1(\A,\epsilon))\subseteq\R^s(\A,\epsilon)$.

\newcommand{\dchen}{d_{\text{Chen}}}
\begin{definition}\label{Rotfuchs} \rm
The sum $\left(d_I+ d_C\right)$ induces a map
\[
\dchen: \bigoplus\limits_{r=0}^\infty {\bigotimes}^r \Ae \bigg/ 
\R^s_1(\A,\epsilon)
\longrightarrow
\left(\bigoplus\limits_{r=0}^\infty {\bigotimes}^r \A\right) \bigg/ 
\R(\A,\epsilon),
\]
the {\it Chen differential} and its kernel consists of what 
is called {\it Chen-closed elements}. We call an element $I\in
\bigoplus_{r=0}^\infty {\bigotimes}^r \Ae$ 
{\it Chen-closed} if $\dchen([I])=0$.
For each $s$ the kernel of $\dchen$ gets a name 
\[
H^0{B}_s(\A,\epsilon) :=
\ker\left\{\dchen:\frac{\left(\bigoplus_{r=0}^s{\bigotimes}^r\Ae \right)} 
{\R^s_1(\A,\epsilon)} \rightarrow
\frac{\left(\bigoplus_{r=0}^\infty {\bigotimes}^r \A\right)} 
{\R(\A,\epsilon)}\right\}.
\]
Let $H^0\bar{B}_s(\A,\epsilon)$ consist of those elements in 
$H^0 {B}_s(\A,\epsilon)$ with zero constant term, i.~e.~define
$H^0\bar{B}_s(\A,\epsilon):=\{ [I] \in H^0{B}_s(\A,\epsilon) |\;
I\in \bigoplus_{r=1}^s{\bigotimes}^r\Ae \}.$
An element of $H^0{B}_s(\A,\epsilon)$ respectively $H^0\bar{B}_s(\A,\epsilon)$
represented by $\sum_{J}a_J\;\varphi_{j_1}\otimes\cdots\otimes
\varphi_{j_r}$ is classically written as 
\[
\sum_{J}a_J\;[\varphi_{j_1}|\cdots |\varphi_{j_r}].
\]
\end{definition}

Observe that we have inclusions for all $s\ge 1$ 
\[
H^0{B}_s(\A,\epsilon)\subset H^0{B}_{s+1}(\A,\epsilon)
\quad \text{   and   }\quad 
H^0\bar{B}_s(\A,\epsilon)\subset H^0\bar{B}_{s+1}(\A,\epsilon).
\]
\begin{remark}\rm
In the general version of the `Reduced Bar Construction' 
(cf.~\cite{Chen-RBC} or \cite{Hain-de-Rham-homotopy}) the vectorspaces $H^0{B}_s(\A,\epsilon)$ and $H^0\bar{B}_s(\A,\epsilon)$ appear as first 
cohomology groups of certain complexes, whence the notation. 
\end{remark}

The following four Lemmas will turn out to be useful later in
Section \ref{Placenta}. 
Let $\bA$ be a connected sub-dga of $\A$ such that 
${\bar{\mathcal A}}^p={\mathcal A}^p$ for $p>1$ and 
\begin{equation} \label{Vierdaagse}
\Ae = d {\mathcal A}^0 \oplus \bAe.
\end{equation}
The proof of the following two lemmas is similar to the proof of the
Theorem on p.23 in \cite{Chen-RBC}. 

\begin{lemma} \label{Claus Hertling}
Suppose $I\in \bigoplus_{r=0}^s {\bigotimes}^r\Ae$. Then 
there is an $R\in \R_1^s(\A,\epsilon)$ and an
$\bar{I}\in\bigoplus_{r=0}^s {\bigotimes}^r\bAe$ such that 
$I=\bar{I}+R$ and 
\begin{gather*}
\dchen [I] =0 \quad \Leftrightarrow
\quad (d_I+d_C)(\bar{I})=0.
\end{gather*}
\end{lemma}

\noindent
{\bf Proof:}
By finite induction on $\tilde{s}$ we prove that for every 
$\tilde{s}$ between $s$ and $0$ there is an 
$I^{\tilde{s}}\in\bigoplus_{r=0}^{\tilde{s}}{\bigotimes}^r\Ae$ and an 
$R_{\tilde{s}}\in \R_1^s(\A,\epsilon)$ as well as an
$\bar{I}_{\tilde{s}}\in\bigoplus_{r=0}^s {\bigotimes}^r\bAe$ such that
$I=\bar{I}_{\tilde{s}}+I^{\tilde{s}}+R_{\tilde{s}}$.
First, for
{\bf ${\bf \tilde{s}=s}$}
let $I^{\tilde{s}}=I$ and $\bar{I}_{\tilde{s}}=R_{\tilde{s}}=0$.

\noindent
{\bf ${\bf \tilde{s}<s}$:}
Assume now that there is an
$I^{\tilde{s}+1}\in\bigoplus_{r=0}^{\tilde{s}+1} {\bigotimes}^r\Ae$ and an 
$R_{\tilde{s}+1}\in \R_1^s(\A,\epsilon)$ as well as an
$\bar{I}_{\tilde{s}+1}\in\bigoplus_{r=0}^s {\bigotimes}^r\bAe$ such that
\[
I=\bar{I}_{\tilde{s}+1}+I^{\tilde{s}+1}+R_{\tilde{s}+1}.
\]
Write $I^{\tilde{s}+1}=\sum_{|J|\le \tilde{s}+1} a_J \;
\varphi_{j_1}\otimes\cdots\otimes\varphi_{j_r}$ and let $\varphi_{j_\nu}=
\bar{\varphi}_{j_\nu}+df_{j_\nu}$ be the decomposition according to 
\eqref{Vierdaagse}. Then observe that
\begin{align*}
& \begin{array}{rrl} \!\!\! I^{\tilde{s}} 
:=\sum\limits_{|J|\le \tilde{s}+1} a_J \; & \big\{ 
df_{j_1}\otimes\varphi_{j_2}\otimes \cdots \otimes \varphi_{j_r} 
 & -\; R_1(\varphi_{j_2}\otimes \cdots \otimes \varphi_{j_r};f_{j_1})\\
& + \;\bar{\varphi}_{j_1}\otimes df_{j_2} \otimes \cdots \otimes \varphi_{j_r} 
& - \;R_2(\bar{\varphi}_{j_1}\otimes \varphi_{j_3}\otimes \cdots 
\otimes \varphi_{j_r};f_{j_2})\\
& \vdots\qquad\qquad &  \qquad\qquad\vdots \\
 & +\; \bar{\varphi}_{j_1} \otimes \cdots \otimes 
\bar{\varphi}_{j_{r-1}}df_{j_r}
& -\; R_r(\bar{\varphi}_{j_1}\otimes \cdots\otimes\bar{\varphi}_{j_{r-1}};
f_{j_r})\big\} \end{array} \\
\intertext{and $ 
\bar{I}_{\tilde{s}}:= \sum\limits_{|J|\le \tilde{s}+1} a_J \; 
\bar{\varphi}_{j_1}\otimes \cdots \otimes \bar{\varphi}_{j_r}$ as well as}
& \begin{array}{rl} \!\!\! R_{\tilde{s}} 
:=\sum\limits_{|J|\le \tilde{s}+1} a_J \;  \big\{ 
   & R_1(\varphi_{j_2}\otimes \cdots \otimes \varphi_{j_r};f_{j_1})\\
 + & R_2(\bar{\varphi}_{j_1}\otimes \varphi_{j_3}\otimes \cdots 
\otimes \varphi_{j_r};f_{j_2})\\
 &  \qquad\qquad\vdots \\
 + & R_r(\bar{\varphi}_{j_1}\otimes \cdots\otimes\bar{\varphi}_{j_{r-1}};
f_{j_r})\big\} \end{array} 
\end{align*}
satisfy the induction hypothesis for $\tilde{s}$. This completes the induction.
Finally let $\bar{I}:=\bar{I}_0$ and note that
$\dchen\left(\left[{\bar{I}\;}\right]\right) = 0 \; 
\Leftrightarrow \; (d_I+d_C)\left({\bar{I}\;}\right) = 0$. \qed

\begin{lemma}	    \label{Mimikry}
Suppose $I\in \bigoplus_{r=0}^s \bigotimes^r \Ae$ is Chen-closed.
Then we can write 
\[
I=\sum\nolimits_{|J|=r\le s} 
a_J \;\varphi_{j_1}\otimes\cdots\otimes\varphi_{j_r},
\]
where for any $J$ with $|J|=s$ and $a_J \neq 0$ either all $\varphi_{j_1},
\ldots,\varphi_{j_s}$ are closed or at least one of them is exact.
\end{lemma}

\noindent
{\bf Proof:} Let $I\in \bigoplus_{r=0}^s \bigotimes^r \Ae$ be Chen-closed.
By Lemma \ref{Claus Hertling} we may assume without loss of generality
that $I\in \bigoplus_{r=0}^s \bigotimes^r \bAe$ and $(d_I+d_C)I=0$.
Write $I=\sum_{J} a_J' \;\varphi_{j_1}'\otimes\cdots\otimes\varphi_{j_r}'$ and
consider the part of $I$, which
is contained in the s-th tensor power of $\Ae$, say
$T^s:=\sum_{|J|=s} a_J' \;\varphi_{j_1}'\otimes\cdots\otimes\varphi_{j_r}' $,
as element of $\left[{\bigotimes}^s \A \right]^s$.
The condition $( d_C +d_I) (\sum_{J} a_J' \;\varphi_{j_1}'
\otimes\cdots\otimes\varphi_{j_r}' )=0$ implies that $T^s$ is a cocycle
with respect to the differential of
$\left({\bigotimes}^s \A,\; d\right)$. The K{\"u}nneth formula
(\cite{Spanier}, p.~228 and p.~247, Theorem 11) yields the short exact
sequence (i.~e.~the isomorphism)
\[
0\longrightarrow
\left[{\bigotimes}^s H^\bullet(\A)\right]^s\longrightarrow
H^s\left({\bigotimes}^s \A\right)\longrightarrow 0.
\]
This is the reason why we can write $T^s$ as
$\sum_{|J|=s} a_J \;\varphi_{j_1}\otimes\cdots\otimes\varphi_{j_r} $,
where for each $J$ with $|J|=s$ either all forms $\varphi_{j_1},
\ldots,\varphi_{j_s}$ are closed or at least one of them is exact.
\hfill $\Box$

Suppose that $H^2(\A)=0$. Given closed elements 
$\varphi_1,\ldots,\varphi_s\in\Ae$, the next proposition shows us, 
how we can find Chen-closed elements in
$\bigoplus_{r=0}^\infty {\bigotimes}^r \Ae$ whose part in
${\bigotimes}^s \Ae$ is the tensor product
$\varphi_1\otimes\cdots\otimes\varphi_s$.

\begin{lemma} \label{Ligeti}
Suppose $H^2(\A)=0$ and assume $\varphi_1,\ldots,\varphi_s$
are closed elements in $\Ae$.
Then there exist elements $\varphi_{i_1,\ldots,i_\nu}\in\Ae$ 
with $\nu>1$ 
such that for all $\nu \ge 1$ (i.~e.~the forms $\varphi_1,\ldots,\varphi_s$
included) and for any multi-index $(i_1,\ldots,i_\nu)\in\{1,\ldots,s\}^\nu$
holds:
\[
d \varphi_{i_1,\ldots,i_\nu} + \sum_{k=i}^{\nu-1}
\varphi_{i_1,\ldots,i_k}\wedge \varphi_{i_{k+1},\ldots,i_\nu} =0.
\]
The following elements in 
$\bigoplus_{r=1}^\nu {\bigotimes}^r \Ae$ are closed under $(d_I+d_C)$:
\[
I=\sum\limits_{r=1}^{\nu}\quad\sum\limits_{0<{\alpha_1}<\cdots 
<\alpha_{r-1}<\nu}
\varphi_{i_1,\ldots,i_{\alpha_1}}\otimes\varphi_{i_{\alpha_1+1},
\ldots,i_{\alpha_2}}\otimes\cdots\otimes\varphi_{i_{\alpha_{r-1}+1},
\ldots,i_s}.
\]
\end{lemma}

\noindent
{\bf Proof:}
The proof of the existence of $\varphi_{i_1,\ldots,i_\nu}\in\Ae$ with
$\nu>1$ goes by induction on $\nu$. The computations are straightforward.
\qed
 
\subsection{Definition of Iterated Integrals on the Nearby 
Fundamental Group} \label{Douwe Egberts}

Let us begin with the definition of an $\varepsilon$-iterated integral
of a $\varphi_1\otimes\cdots\otimes\varphi_r\in\bigotimes^r A^1$ along
a path over $\vec{v}$ (based at $\vec{w}$). Each of the $\varphi_\nu\in A^1$
can be written as
\[
\varphi_\nu=\sum_{i\ge 0} \omega_i^{(\nu)}+\sum_{[k<l]} K_{kl}^{(\nu)}
\frac{dx}{x} + L_{kl}^{(\nu)} \frac{dy}{y} + H_{kl}^{(\nu)} d\xi.
\]
Choose around every double point $p_{kl}$ coordinates
$(x,y):W_{kl}=U_{kl}^k\times U_{kl}^l\rightarrow \cpxm^2$ and a
coordinate $t:\Delta\rightarrow \cpxm$ such that $h_{\big| W_{kl}}(x,y)
=x\cdot y$ and $\vec{v}=\frac{\partial}{\partial t}$.
For the double point $p_0\in D_0$ we require moreover that
the coordinate $x$ (which we also called $p$) satisfies 
$\frac{\partial}{\partial x}=
\frac{\partial}{\partial p}=\vec{w}$.
Always, when we talk about $\varepsilon$-iterated integrals in this
section, we will refer to this coordinate system.

As for the definition of line integrals in \ref{Kormoran} we give
the definition of iterated integrals in several steps.

\begin{description}
\item[(a)] 
Let $\gamma:[a;b]\rightarrow Z_0$ be a path over $\vec{v}$, which meets the set
of double points only once with parameter value $\tau_0\in ]a;b[$, where
it changes from $D_k$ to $D_l$.
Define $\gamma_x(\tau)$ (resp.~$\gamma_y(\tau)$) to be $x(\gamma(\tau))$
(resp.~$y(\gamma(\tau))$) for all $\tau\in [a;b]$ with $\gamma(\tau)
\in U_{kl}^k$ (resp.~$\gamma(\tau)\in U_{kl}^l$).
To shorten the following formulas we abbreviate for $\varepsilon>0$ small
enough:
\[
\eta_\varepsilon^{(\nu)}:=H_{kl}^{(\nu)}\big(\xi,\log\left(
\gamma_x(\tau_0-\varepsilon)
\cdot \gamma_y (\tau_0+\varepsilon)\right)\big)\; d\xi.
\]
Then define
\footnote{Here and also later we make the convention
that all the expressions with an integral sign, where
from the upper indices the right one is smaller than the left one,
like $\int_{[0;1]} \eta_\varepsilon^{(5)}\cdots
\eta_\varepsilon^{(4)}$, take the value 1.}
\[
\eint{\gamma} \!\! \varphi_1\cdots \varphi_r
\!:= \sum_{0\le\alpha\le\beta\le r}
\;\;\int\limits_{\gamma^{\le \tau_0-\varepsilon}} \!\!\!\!\!\!\omega_k^{(1)}
\cdots\omega_k^{(\alpha)}
\!\!\!\int\limits_{[0;1]} \eta_\varepsilon^{(\alpha+1)}\cdots
\eta_\varepsilon^{(\beta)}
\!\!\!\!\!\!\!\!\int\limits_{\gamma_{\ge \tau_0+\varepsilon}} \!\!\!\!\!\!
\omega_l^{(\beta+1)}\cdots \omega_l^{(r)}.
\]
\item[(b)]
Let $\gamma:[a;b]\rightarrow Z_0$ be a path over $\vec{v}$, which meets the set
of double points only once with parameter value $\tau_0\in ]a;b[$, where
it stays in one component $D_k$. Then define:
\[
\eint{\gamma} \varphi_1\cdots \varphi_r:=\sum_{0\le \alpha \le r}\;\;
\int\limits_{\gamma^{\le \tau_0-\varepsilon}} \omega_k^{(1)}\cdots
\omega_k^{(\alpha)}
\int\limits_{\gamma_{\ge \tau_0+\varepsilon}} \omega_k^{(\alpha+1)}
\cdots\omega_k^{(r)}.
\]
\item[(c)]
Now let $\gamma:[a;b]\rightarrow Z_0$ be a path over $\vec{v}$
such that $\gamma(a)$ and $\gamma(b)$ are no double points.
Then there is a finite number of parameter values $\tau_1,\ldots,\tau_N$
with $a<\tau_1<\cdots<\tau_N<b$, which are mapped onto double points.
Choose a $\tau_i^*\in ]\tau_i;\tau_{i+1}[$ for $i=1,\ldots,N-1$ and let
$\tau_0^*:=a$ and $\tau_N^*:=b$. Then define for $i=1,\ldots,N$ the paths
$\gamma_i:= \gamma_{\big| [\tau_{i-1}^*,\tau_i^*]}$ with which we define
$\eintt{\gamma} \varphi_1\cdots \varphi_r$ as 
\[
\sum_{0\le\alpha_1\le\cdots\le\alpha_{N-1}\le r}\;\;
\eint{\gamma_{1}}
\varphi_1\cdots \varphi_{\alpha_1}
\eint{\gamma_{2}}
\varphi_{\alpha_1+1}\cdots \varphi_{\alpha_2} \cdots
\eint{\gamma_{N}}
\varphi_{\alpha_{N-1}+1}\cdots \varphi_{\alpha_r}.
\]
\item[(d)] 
Finally, let $\gamma:[a;b]\rightarrow Z_0$ be a path in $Z_0$ over $\vec{v}$
based at $\vec{w}$. The coordinate $p:D_0\rightarrow\cpxm$ allows us to 
consider $D_0$ as part of the complex plane. There, in the complex plane,
we have the differential form (and not the symbol) $\frac{dp}{p}$. 
Let like before $\sigma$ be the 
straight path $\tau\mapsto (1-\tau)$ in $\cpxm$, then we define 
\[
\eint{\gamma} \varphi_1\cdots \varphi_r
:=\eint{\sigma\star\gamma\star \sigma^{-1}}
\varphi_1\cdots \varphi_r.
\]
\end{description}

We leave it to the reader to prove the following proposition. It is a
consequence of the corresponding formula for iterated integrals on manifolds
(see \cite{Hain-the-geometry}, Prop.~2.9).
\begin{proposition} \label{Badings}
The definition (c) does not depend on the choice of the $\tau_i^*$.
Let $\bf p$, $\bf q$, ${\bf r}\in Z_0\setminus\bigcup_{[k<l]}\{p_{kl}\}$ and let
$\alpha,\;\beta:[0;1]\rightarrow Z_0$ be paths over $\vec{v}$,
where $\alpha(0)={\bf p}$, $\alpha(1)=\beta(0)={\bf q}$ and $\beta(1)={\bf r}$.
Then for forms $\varphi_1,\ldots,\varphi_r\in A^1$ and $\varepsilon>0$
small enough holds:
\[
\eint{\alpha\star\beta}\varphi_{1}\cdots \varphi_r=
\sum_{0\le m\le r}\quad\eint{\alpha} \varphi_1\cdots \varphi_m\;
\eint{\beta} \varphi_{m+1}\cdots \varphi_r
\]
\qed
\end{proposition}

This behaviour of iterated integrals under composition of paths
will lead to the definition of $\epsT$ in Section \ref{Placenta}, 
motivated by Chen's theory of formal connections.
The main result of Section \ref{Placenta} will be the proof the 
following theorem.
 
\bigskip
\noindent
{\bf Theorem \ref{Reinhard May}} \; {\it
Let $\chen$ be a Chen-closed element and let $\gamma:[a;b]\rightarrow
Z_0$ be a path over $\vec{v}$ based at $\vec{w}$. Then
\[
\lim_{\varepsilon\to 0} \eint{\gamma} 
\sum_J a_J\; \varphi_{j_1}\cdots\varphi_{j_r}
\]
converges (\ref{Babymassage}) and does not depend on the choice of 
coordinates (\ref{Vormilch}).
Moreover it does not depend upon the choice of $\gamma$ within
a nearby homotopy class of paths over $\vec{v}$ based at $\vec{w}$
(\ref{Blaek Foes}), 
i.~e.~we have a nearby homotopy functional:
\[
\int \sum_J a_J\; \varphi_{j_1}\cdots\varphi_{j_r}:
\pi_1(Z_{\vec{v}},\vec{w})\rightarrow \cpxm.
\]
Finally, this homotopy functional is zero if $I\in \R_1^s (A^\bullet, a)$.}
\bigskip

Theorem \ref{Reinhard May} assures that the 
following definition of iterated integrals of Chen-closed elements along 
paths over $\vec{v}$ makes sense.  

\begin{definition} \label{Taetaerae}\rm
Let $\chen$ be a Chen-closed element. Then {\it the iterated
integral of $I$ along a path $\gamma$ over $\vec{v}$ based at
$\vec{w}$} is defined by:
\[
\int\limits_\gamma \sum_J a_J\; \varphi_{j_1}\cdots\varphi_{j_r}:=
\lim_{\varepsilon\to 0}\left\{
\eint{\gamma} \sum_J a_J\;\varphi_{j_1}\cdots\varphi_{j_r}
\right\}.
\]
\end{definition}

\subsection{The $\pi_1$-De Rham Theorem for the Nearby Fundamental Group}

Denote by $P(Z_{\vec{v}},\vec{w})$
the set of all paths over $\vec{v}$ based at $\vec{w}$. Note that we can
consider an $\varepsilon$-iterated integral not just as a function on 
$P(Z_{\vec{v}},\vec{w})$ but also on $\cpxm P(Z_{\vec{v}},\vec{w})$,
the free abelian group generated by the elements of $P(Z_{\vec{v}},\vec{w})$,
by extending it linearly. Note that elements in $P(Z_{\vec{v}},\vec{w})$
can be composed and this operation extends to $\cpxm P(Z_{\vec{v}},\vec{w})$.

Also the following proposition corresponds to a rule  
for iterated integrals on manifolds (see (2.10), (2.13) 
in \cite{Hain-the-geometry}). 
\begin{proposition} \label{Hoffmannstal}
Let $\alpha_1,\ldots,\alpha_s$ be paths over $\vec{v}$ based at $\vec{w}$
and let $\varphi_1,\ldots,\varphi_r\in A^1$ with $s\ge r$. Then
\[
\eint{(\alpha_1-1)\cdots(\alpha_s-1)}\!\!\!\!\!\!
\varphi_1\cdots\varphi_r =
\left\{\begin{array}{ll}
\eintt{\alpha_1} \varphi_1\cdots \eintt{\alpha_s} \varphi_s
+ V(\varepsilon), & \text{  if } r=s \\
W(\varepsilon), & \text{  if } r<s. \end{array}
\right.
\]
where $V(\varepsilon)$, $W(\varepsilon)$ are functions 
$V,W:\reellm^{>0}\rightarrow\cpxm$ 
converging to $0$ for $\varepsilon\to 0$. 
\end{proposition}

\noindent
{\bf Proof:}
We write $A(\varepsilon)\equiv B(\varepsilon)$
if $A(\varepsilon)-B(\varepsilon)=V(\varepsilon)$ for a function
$V(\varepsilon)$ with $\lim_{\varepsilon\to 0} V(\varepsilon)=0$.
One can prove (and this is easy) the following formula by
induction on $k$:
\begin{eqnarray*}
& & \quad\quad\quad
\eint{(\alpha_1-1)
\cdots(\alpha_k-1)\alpha_{k+1}\cdots\alpha_s}
\varphi_1\cdots \varphi_r\\
& \equiv &
\!\!\!\!\!\!\!\!\!\!\!\!\!\!\!\!\!\!\!\!
\sum\limits_{\topa{k\le N\le s}{0<\nu_1<\cdots<\nu_k\le\nu_{k+1}
\le\cdots\le\nu_{N-1}\le\nu_{N}=r}}
\eint{\alpha_1}\varphi_1\cdots\varphi_{\nu_1}
\eint{\alpha_2}\varphi_{\nu_1+1}\cdots\varphi_{\nu_2}
\cdots
\eint{\alpha_N}\varphi_{\nu_{N-1}+1}\cdots\varphi_{\nu_N}.
\end{eqnarray*}
For $k>r$, the sum is zero and for $k=r=s$ the sum is:
$\eintt{\alpha_1} \varphi_1\cdots \eintt{\alpha_s} \varphi_s$.
\qed

Recall that we defined an augmentation $a:A^\bullet\rightarrow \cpxm$
(see \ref{Weisser Schwan}).
Therefore Definition \ref{Rotfuchs} makes
sense for the pair $(A^\bullet,a)$, whence we have the complex vector spaces
\[
H^0\bar{B}_s(A^\bullet,a)\subset H^0{B}_s(A^\bullet,a).
\]
Note in particular: $H^0\bar{B}_1(A^\bullet,a)=H^1(A^\bullet)$.
Theorem \ref{Reinhard May} and Proposition \ref{Hoffmannstal} yield 
a well-defined integration map
\[
\int: H^0 B_s(A^\bullet,a) \longrightarrow \Hom_\itgm (\itgm
\pi_1(Z_{\vec{v}},\vec{w})\big/{\vec{J}}^{s+1};\cpxm).
\]
The analogue for the nearby fundamental group of Chen's 
$\pi_1$-De Rham Theorem states that this 
map is an isomorphism.

Also the group ring $\itgm\pi_1(Z_{\vec{v}},\vec{w})$ comes with an
augmentation map $\itgm\pi_1(Z_{\vec{v}},\vec{w})\rightarrow
\itgm$. This augmentation has a natural section $\sigma:n\mapsto n{\pmb 1}$
which makes the short exact sequence
\[
0\rightarrow \vec{J}\big/{\vec{J}}^{s+1} \rightarrow
\itgm\pi_1(Z_{\vec{v}},\vec{w}) \big/{\vec{J}}^{s+1}
\xrightarrow{\topa{\sigma}{\dashleftarrow}} \itgm \rightarrow 0
\]
split. This splitting gives a natural isomorphism
\[
\Hom_\itgm\left(\itgm\pi_1(Z_{\vec{v}},\vec{w}) \big/{\vec{J}}^{s+1};
\itgm\right)
\cong
\Hom_\itgm\left(\vec{J}\big/{\vec{J}}^{s+1};\itgm\right)\oplus\itgm
\]
Via the integration map this decomposition corresponds with
\[
H^0{B}_s(A^\bullet,a)\cong H^0\bar{B}_s(A^\bullet,a)\oplus \cpxm.
\]
Therefore it suffices to consider
$ \Hom_\itgm(\vec{J}\big/{\vec{J}}^{s+1},\itgm)$ in order to
describe the dual of 
$\itgm\pi_1(Z_{\vec{v}},\vec{w}) \big/{\vec{J}}^{s+1}$.

\begin{theorem}    \label{Novalis}
Integration of Chen-closed elements in
$\bigoplus_{r=1}^s\bigotimes^r A^1$
along paths over $\vec{v}$ based at $\vec{w}$
defines an isomorphism of complex vector spaces
\[
H^0{B}_s(A^\bullet,a)\xrightarrow{\cong}
\Hom_\itgm \left(\itgm\pi_1(Z_{\vec{v}},\vec{w}) \big/{\vec{J}}^{s+1} ; 
\cpxm\right).
\]
\end{theorem}

\noindent
{\bf Proof:}
We first indicate that for each $s\ge 1$ there is a short exact sequence
\[
0\rightarrow H^0\bar{B}_{s-1}(A^\bullet,a)\rightarrow H^0\bar{B}_s(A^\bullet,a)
\rightarrow \bigotimes^s H^1(A^1)\rightarrow 0.
\]
In particular, the map
\[
\begin{matrix}
H^0\bar{B}_s(A^\bullet,a)& \rightarrow & {\bigotimes}^s H^1(A^1) \\
\sum_J a_J [\varphi_{j_1}|\cdots|\varphi_{j_r}] & \mapsto &
\sum_J a_J [\varphi_{j_1}]\otimes\cdots\otimes[\varphi_{j_r}]
\end{matrix}
\]
is well-defined and is surjective by Proposition \ref{Ligeti}. 
Its kernel is $H^0\bar{B}_{s-1}(A^\bullet,a)$ by Lemma \ref{Mimikry}.
Then by Proposition \ref{Hoffmannstal} we have 
the commutative diagram:\footnote{
The map $(\vec{J}^s\big/\vec{J}^{s+1})^*\rightarrow
\bigotimes^s(\vec{J}\big/\vec{J}^{2})^*$, which is given by 
the ring structure on $\itgm\pi_1(Z_{\vec{v}},\vec{w})$, is  
surjective since we have the
isomorphism $\pi_1(Z_{\vec{v}},\vec{w})\cong \pi_1(Z_t,\sigma_p(t))$
for some $t\in\Delta^*$ and $p\in D_0^*$. In $\pi_1(Z_t,\sigma_p(t))$
it is true that $({J}^s\big/{J}^{s+1})^*  \rightarrow 
\bigotimes^s (J\big/{J}^{\;2} )^*$ is surjective.}
\[
\begin{array}{ccccccccc}
0 &\rightarrow & H^0\bar{B}_{s-1}(A^\bullet,a) & \rightarrow &
H^0\bar{B}_s(A^\bullet,a) & \rightarrow & \bigotimes^s H^1(A^1)&\rightarrow & 0 \\
 & & \downarrow & & \downarrow & & \downarrow & &  \\
0 &\rightarrow & (\vec{J}\big/\vec{J}^s)^* & \rightarrow &
(\vec{J}\big/\vec{J}^{s+1})^* & \rightarrow &
\bigotimes^s (\vec{J}\big/\vec{J}^{\;2} )^* &\rightarrow & 0
\end{array}
\]
The theorem can now be proved by induction on $s$ and the use of the
5-lemma. The beginning of the induction is given by Theorem \ref{Steinkauz}.
\qed

\begin{warning} \label{Schwarzbunte} \rm
Assume that we are given two
closed 1-forms $\varphi,\,\psi\in A^1$,
\[
\varphi=\sum_{i\ge 0} 0 + \sum_{[k<l]} c_{kl} d\xi 
\qquad\text{ and }\qquad
\psi=\sum_{i\ge 0} \omega_i + \sum_{[k<l]} c_{kl}' d\xi,
\]
where all $\omega_i\in E^1(D_i)$. Then $\varphi\wedge\psi = 0$ and
hence is $\int \varphi\psi $ a nearby homotopy functional. 
Define $\psi_0=\sum_{i\ge 0} 0 + \sum_{[k<l]} c_{kl}' d\xi$. When we 
multiply $\varphi$ and $\psi$ (componentwise), we get the same
as if we multiply $\varphi$ and $\psi_0$. But nonetheless, in general holds:
$\int \varphi\psi \neq \int \varphi\psi_0 =0$.
In general, there exists such $\varphi$ and $\psi$ for which there is
a closed path $\alpha$ in one of the $D_i$'s such that 
$\int_{\alpha} \omega_i\neq 0$ and a path $u$ over $\vec{v}$ starting
at $\vec{w}$ and ending at $\alpha(0)=\alpha(1)$ such 
that $\int_u\varphi \neq 0$. Then
\[
\begin{split}
\int_{u\star\alpha\star u^{-1}} \varphi\psi &= 
\int_{u} \varphi \int_{\alpha } \psi 
-\int_{\alpha} \varphi \int_u \psi +\int_{\alpha} \varphi\psi \\
&= \int_{u} \varphi \int_{\alpha } \psi 
\qquad\qquad\qquad\qquad\qquad\qquad\qquad \diamond
\end{split}
\]
\end{warning}

\subsection[The MHS on $\vec{J}\big/\vec{J}^{s+1}$]{ The Mixed Hodge Structure
on $\vec{J}\big/\vec{J}^{s+1}$}

After having developed this theory of iterated integrals along paths
over $\vec{v}$, we are ready now to define a mixed Hodge Structure
on $\vec{J}\big/\vec{J}^{s+1}$. With all what we know about
iterated integrals along paths over $\vec{v}$ yet, this is just
an application of standard techniques in mixed Hodge theory.

\subsubsection{The Hodge and the Weight Filtration}

Let us use Theorem \ref{Novalis} to define a weight - and a Hodge
filtration on $(\vec{J}\big/\vec{J}^{s+1})^*$.

On $\bigoplus_{r=1}^s{\bigotimes}^r A^1$
the {\it weight filtration} $W_\bullet$ is defined by
\[
W_l\left( \bigoplus_{r=1}^s{\bigotimes}^r A^1\right):=
\bigoplus_{r=1}^s\;\bigoplus_{l_1+\cdots+l_r +r\le l}
W_{l_1}A^1\otimes\cdots\otimes W_{l_r}A^1
\]
and the {\it Hodge filtration} $F^\bullet$ by
\[
F^p\left( \bigoplus_{r=1}^s{\bigotimes}^r A^1\right):=
\bigoplus_{r=1}^s\;\bigoplus_{p_1+\cdots+p_r \ge p}
F^{p_1}A^1\otimes\cdots\otimes F^{p_r}A^1.
\]
These filtrations $W_\bullet$ and $F^\bullet$ on
$\bigoplus_{r=1}^s\bigotimes^r A^1$ induce filtrations 
$W_\bullet$ and $F^\bullet$ on the subspace $\ker(d_I+d_C)$ and on
$H^0\bar{B}_s(A^\bullet,a)\subset H^0{B}_s(A^\bullet,a)$ by:
\begin{align*}
W_l H^0{B}_s(A^\bullet,a) &:= \im\left\{ W_l \ker(d_I+d_C) 
\rightarrow H^0\bar{B}_s(A^\bullet,a)\right\} \\
F^p H^0{B}_s(A^\bullet,a) &:= \im\left\{ F^p \ker(d_I+d_C) 
\rightarrow H^0\bar{B}_s(A^\bullet,a)\right\}. 
\end{align*}

\begin{proposition} \label{Eckermann}
The short exact sequences (with $s\ge 2$)
\[
0\rightarrow H^0\bar{B}_{s-1}(A^\bullet,a) \xrightarrow{\iota}
H^0\bar{B}_s(A^\bullet,a) \xrightarrow{pr}
\bigotimes^s H^1(A^1)\rightarrow 0
\]
are strict with respect to $W_\bullet$ and $F^\bullet$.
\end{proposition}

\noindent
{\bf Proof:}
It is easy to see that $pr$ is strict. We shall prove that
$\iota$ is strict with respect to both filtrations.

\noindent
{\bf Hodge filtration:}
Let $\bar{I}=\sum_J a_J \;[\varphi_{j_1}|\cdots|\varphi_{j_r}] 
\in F^p H^0\bar{B}_s(A^\bullet,a)\cap
H^0\bar{B}_{s-1}(A^\bullet,a)$. Then there is a Chen-closed
$I=\sum_{|J|\le s} a_J \varphi_{j_1}\otimes\cdots\otimes\varphi_{j_r}
\in F^p \bigoplus_{r=1}^s\bigotimes^r A^1$ such that
$\int I=\int \bar{I}$. Since $\bar{I} \in H^0\bar{B}_{s-1}(A^\bullet,a)$ 
we may assume
by Lemma \ref{Mimikry}, Theorem \ref{Reinhard May} and Proposition
\ref{Hoffmannstal} that for $J$ with $|J|=s$ all
the $\varphi_{j_1},\ldots,\varphi_{j_s}$ are exact:
$\varphi_{j_m} = d f_{j_m}$ and $f_{j_m}(p_{01})=0$.
We have by definition:
$\varphi_{j_1}\otimes\cdots\otimes\varphi_{j_s}\in F^p {\bigotimes}^s A^1$
$:\Leftrightarrow$ $\exists (p_1,\ldots,p_s)\in \itgm^s$ with 
$p_1+\cdots+p_s\ge p:$ $\varphi_{j_1}\in F^{p_1},\ldots,
\varphi_{j_s}\in F^{p_s}$.
Now there is either one $p_m\ge 1$  (which implies $\varphi_{j_m}=0$
as $A^0\stackrel{d}{\rightarrow} A^1$ is strict with respect to $F^\bullet$,
i.~e.~$d A^0\cap F^1 A^1=0$) or $p_1=\cdots=p_s=0$. 
Then
\[
\left(f_{j_1}-a(f_{j_1})\right) \varphi_{j_2}\otimes\varphi_{j_3}
\otimes\cdots\otimes\varphi_{j_s}
\in F^p {\bigotimes}^{s-1} A^1.
\]
Therefore, there is a Chen-closed $\tilde{I}\in F^p
\bigoplus_{r=1}^{s-1}{\bigotimes}^{r} A^1$ such that $\int \bar{I}
=\int\tilde{I}$.

\noindent
{\bf Weight filtration:}
Let $\bar{I}\in W_l H^0\bar{B}_s(A^\bullet,a)\cap
H^0\bar{B}_{s-1}(A^\bullet,a)$. Then there is a Chen-closed
$I=\sum_{|J|\le s} a_J \varphi_{j_1}\otimes\cdots\otimes\varphi_{j_r}
\in W_l\bigoplus_{r=1}^s\bigotimes^r A^1$ such that
$\int I=\int \bar{I}$. Because of 
$\bar{I} \in H^0\bar{B}_{s-1}(A^\bullet,a)$ we may again assume
by Lemma \ref{Mimikry}, Theorem \ref{Reinhard May} and Proposition
\ref{Hoffmannstal} that for $J$ with $|J|=s$ all
the $\varphi_{j_1},\ldots,\varphi_{j_s}$ are exact:
$\varphi_{j_m} = d f_{j_m}$ and $f_{j_m}(p_{01})=0$.

Note that if $df\in W_\ell A^1$ then $f\in W_{\ell+1} A^0$. This can be seen 
as follows. Write $f=\sum_{i\ge 0} g_i + 
\sum_{[k<l]} P_0+ P_1 u + \cdots P_m u^m$. Now $f\in W_{\ell+1}\Leftrightarrow
2m\le \ell+1$ and $df\in W_{\ell}\Leftrightarrow 2m-1\le \ell$.
We have by definition: 
$\varphi_{j_1}\otimes\cdots\otimes\varphi_{j_s}\in W_l {\bigotimes}^s A^1$
$:\Leftrightarrow$ $\exists (l_1,\ldots,l_s)\in \itgm^s$ with 
$l_1+\cdots+l_s+s\le l:$ $\varphi_{j_1}\in W_{l_1},\ldots,
\varphi_{j_s}\in W_{l_s}$.
Then
\[
\left(f_{j_1}-a(f_{j_1})\right)\varphi_{j_2}\otimes
\varphi_{j_3}\otimes\cdots\otimes\varphi_{j_s}
\in W_l {\bigotimes}^{s-1} A^1.
\]
Hence, there is a Chen-closed $\tilde{I}\in W_l
\bigoplus_{r=1}^{s-1}{\bigotimes}^{r} A^1$ such that $\int\tilde{I}=
\int\bar{I}$. \qed

The main theorem of this section is the following.
\begin{theorem}  \label{Blaumeise}
For any $s\ge 1$ is 
\[
( \vec{J}\big/\vec{J}^{s+1})^*:=
\left(( \vec{J}\big/\vec{J}^{s+1})^*_\itgm,\,
(( \vec{J}\big/\vec{J}^{s+1})^*_\ratm, W_\bullet)\,
(( \vec{J}\big/\vec{J}^{s+1})^*_\cpxm, W_\bullet, F^\bullet)\right)
\]
a mixed Hodge structure and the short exact sequence 
\begin{equation} 
0\rightarrow ( \vec{J}\big/\vec{J}^{s})^*\rightarrow
( \vec{J}\big/\vec{J}^{s+1})^*\rightarrow 
H^1(A^\bullet)^{\otimes s} \rightarrow 0
\end{equation}
is a short exact sequence of MHSs.
\end{theorem}

\noindent
{\bf Proof:}
Using Proposition \ref{Eckermann} and 1.16 in \cite{Griffiths-Schmid}, the
proof of the theorem can be done by induction on $s$.
\qed

\subsection{The Monodromy on $\vec{J}\big/\vec{J}^{s+1}$} \label{Blindfisch}

Recall that we defined in \ref{Maibaum} the nilpotent chain-morphism 
\[
\begin{array}{cccc}
N: & A^\bullet & \rightarrow & A^\bullet \\*
   & \phi  & \mapsto & \frac{d}{d(-u)} \phi = - \frac{d}{du} \phi \;.
\end{array}
\]
with the properties $N(W_{l+1}A^\bullet)\subset W_{l-1}A^\bullet$ and 
$N(F^{p}A^\bullet)\subset F^{p-1}A^\bullet$. 

Now we define a second 
nilpotent chain-morphism:
\[
M:A^\bullet\rightarrow A^\bullet
\]
by defining it first on $B^\bullet$ in the following way: Let $M$ be the zero
map on 
\[
\bigoplus_{i\le r-1}{\bigwedge}^\bullet_i (\frac{dp}{p})\oplus 
\bigoplus_{i>r-1} E^\bullet (D_i \log P_i)
\oplus\bigoplus_{[0<k<l]} E^\bullet(\Delta^1)\otimes{\Lambda}^\bullet.
\]
Only on the one component $E^\bullet(\Delta^1)\otimes{\Lambda}^\bullet$ for 
$(k,l)=(0,1)$ it is, like $N$, defined by
\[
\begin{array}{cccc}
M: & E^\bullet(\Delta^1)\otimes{\Lambda}^\bullet & \rightarrow & 
E^\bullet(\Delta^1)\otimes{\Lambda}^\bullet \\*
   & \Xi  & \mapsto & \frac{d}{d(-u)} \Xi = - \frac{d}{du} \Xi \;.
\end{array}
\]
Also $M$ satisfies $M(W_{l+1}A^\bullet)\subset W_{l-1}A^\bullet$ and 
$M(F^{p}A^\bullet)\subset F^{p-1}A^\bullet$.

\begin{remark} \rm
Note that for a $Q(u)\in E^\bullet(\Delta^1)\otimes{\bigwedge}^\bullet
(\frac{dx}{x},\frac{dy}{y})[u]$ and for a $\lambda\in\cpxm^*$ always holds:
\[
\lambda^N Q(u)= Q(u+\log \lambda).
\]
Similar for $M$. Moreover, $N$ and $M$ define linear maps
\begin{equation} \label{Rich Hornung}
N,\; M:\bigoplus_{r=1}^s\bigotimes^r A^\bullet \rightarrow
\bigoplus_{r=1}^s\bigotimes^r A^\bullet, 
\end{equation}
which are defined on a generator 
$\varphi_{j_1}\otimes\cdots\otimes\varphi_{j_r}$ by the Leibniz rule:
\begin{align*}
N\left(\varphi_{j_1}\otimes\cdots\otimes\varphi_{j_r} \right)&=
\sum_{\nu=1}^r \varphi_{j_1}\otimes\cdots\otimes N\varphi_{j_\nu}\otimes
\cdots\otimes \varphi_{j_r}
\quad\text{ resp. }\\
M\left(\varphi_{j_1}\otimes\cdots\otimes\varphi_{j_r}\right) &=
\sum_{\nu=1}^r \varphi_{j_1}\otimes\cdots\otimes M\varphi_{j_\nu}\otimes
\cdots\otimes \varphi_{j_r}.
\end{align*}
\end{remark}

The proof of the following proposition is a standard computation:
\begin{proposition} For $\lambda\in\cpxm^*$ and 
$\varphi_{j_1}\otimes\cdots\otimes\varphi_{j_r}\in\bigotimes^r A^\bullet$
holds:
\begin{align*}
\lambda^N\left(\varphi_{j_1}\otimes\cdots\otimes\varphi_{j_r}\right)
&=(\lambda^N\varphi_{j_1})\otimes\cdots\otimes(\lambda^N\varphi_{j_r})
\quad\text{and} \\
\lambda^M\left(\varphi_{j_1}\otimes\cdots\otimes\varphi_{j_r}\right)
&=(\lambda^M\varphi_{j_1})\otimes\cdots\otimes(\lambda^M\varphi_{j_r}). 
\end{align*}
\qed
\end{proposition}

Both, $N$ and $M$ in (\ref{Rich Hornung}) induce maps 
$N,\,M:H^0\bar{B}_s(A^\bullet,a)\rightarrow H^0\bar{B}_s(A^\bullet,a)$.
The complex vector space $H^0\bar{B}_s(A^\bullet,a)$ does neither 
depend on $\vec{v}$ nor on $\vec{w}$. Therefore, we can study, how
the lattice 
\[
({J}_{\vec{v},\vec{w}}\big/{J}_{\vec{v},\vec{w}}^{s+1})^*_\itgm
\subset H^0\bar{B}_s(A^\bullet,a)
\]
behaves, when we move $\vec{v}\in(T_0\Delta)^*$ and $\vec{w}\in (T_{p_0}D_0)^*$.
Let $T$ and $S$ be the monodromies of the local systems 
\[
\left\{\left({J}_{\vec{v},\vec{w}}\big/
{J}_{\vec{v},\vec{w}}^{s+1}\right)^*_\itgm\right\}_{\vec{v}\in(T_0\Delta)^*}
\text{ and }
\left\{\left({J}_{\vec{v},\vec{w}}\big/
{J}_{\vec{v},\vec{w}}^{s+1}\right)^*_\itgm\right\}_{\vec{w}\in (T_{p_0}D_0)^*}
\]
respectively. The main result of this subsection is the following theorem.
\begin{theorem} \label{Gr"unspan}
$e^{-2\pi i N}=T$ and $e^{2\pi i M}=S$. For $\lambda,\,\mu\in\cpxm^*$ holds:
\[
\left({J}_{\lambda\vec{v},\mu\vec{w}}\big/
{J}_{\lambda\vec{v},\mu\vec{w}}^{s+1}\right)^*_\itgm=\lambda^{-N}\mu^{M}
\left({J}_{\vec{v},\vec{w}}\big/{J}_{\vec{v},\vec{w}}^{s+1}\right)^*_\itgm
\subset H^0\bar{B}_s(A^\bullet,a).
\]
\end{theorem}
The proof of this theorem is given in 
\ref{Hokuspokus}. As a consequence of Theorem \ref{Gr"unspan} 
we obtain the following.

\begin{theorem}
For any $s\ge 1$ the family of mixed Hodge structures
\[
\left\{\left(J_{\vec{v},\vec{w}}\big/J^{s+1}_{\vec{v},\vec{w}}\right)^* 
\right\}_{\vec{v}\in (T_0\Delta)^*}
\]
is a nilpotent orbit of mixed Hodge structure. \qed
\end{theorem}

\subsection[The MHS on the Fundamental Group of $Z_0$]{The MHS 
on the Fundamental Group of the Central Fiber}

Let $\pi_1(Z_0,p_0)$ be the fundamental group of the central fiber and let 
$\itgm\pi_1(Z_0,p_0)$ be its group ring with augmentation ideal $J_0$. 
In this subsection we define a MHS on the fundamental group of the central fiber.
Apart from the contribution coming from the non compact disks $D_i$ with
$i=1,\ldots,r-1$, the 
construction of this MHS is a special case of the general construction
of a MHS on the fundamental group of a complex algebraic varietiy as 
introduced by Hain \cite{Hain-de-Rham-homotopy}. The main result here in
this subsection is that the obvious group homomorphism 
\[
c:\pi_1(Z_{\vec{v}},\vec{w})\rightarrow \pi_1(Z_{0},p_0)
\]
induces inclusions of MHSs:
\[
c^*:\big(J_0\big/J_0^{s+1}\big)^*\rightarrow 
\big(\vec{J}\big/\vec{J}^{s+1}\big)^*.
\]
Most proofs are very similar to the proofs of the corresponding statements
for $\pi_1(Z_{\vec{v}},\vec{w})$ and are therefore left to the reader.

\subsubsection[A DGA of Differential Forms on $Z_0$]{A DGA of 
Differential Forms on the Central Fiber}

Now let us construct a dga $A_0^\bullet$ which computes the cohomology 
of $Z_0$ and allows to define the MHS on $\pi_1(Z_0,p_0)$.
Define
\[
B_0^\bullet:= \cpxm\oplus \bigoplus_{i>0} E^\bullet(D_i)\oplus
\bigoplus_{[k<l]} E^\bullet(p_{kl})\otimes E^\bullet(\Delta^1).
\]

By setting $\cpxm={\bigwedge}^\bullet(\frac{dp}{p})$ and 
$E^\bullet(p_{kl})=\Lambda^0$
we may consider $B_0^\bullet$ as sub-dga of the complex $B^\bullet$
as it was defined in \ref{Ultraschall}.
Now define the dga:\footnote{We avoid to introduce 
the language of simplicial and cosimplicial
objects here. However except for the contribution from
$D_0,\dots,D_{r-1}$, the complex $A_0^\bullet$ is isomorphic to the de 
Rham complex of the (complex part of the) cosimplicial mixed Hodge complex 
on the simplicial variety $Z_0$.}
\[
A_0^\bullet:= B_0^\bullet\cap A^\bullet.
\]
Note that: $A_0^\bullet=W_0 A^\bullet$, i.~e.~$A_0^\bullet$ consists 
of exactly those elements in $A^\bullet$,
which have no poles and no ``$\log t=u$''. The dga 
$A_0^\bullet$ inherits a Hodge- and a weight filtration $F^\bullet$ and 
$W_\bullet$ as well as an augmentation map $a:A_0^\bullet\rightarrow\cpxm$ 
from $A^\bullet$.

Similarly to Theorem \ref{Graugans} one can prove:
\begin{theorem}
$H^\bullet(A_0^\bullet)\cong H^\bullet(Z_0;\cpxm)$. \qed
\end{theorem}

\subsubsection{Collapsing the Vanishing Cycles}

A path $\gamma:[0;1]\rightarrow Z_0$ over $\vec{v}$ based at $\vec{w}$ is 
in particular a path in $D^+$ based at $p_0$. 
We leave it to the reader to prove:
\begin{proposition}
The obvious map
\[
c: \pi_1(Z_{\vec{v}},\vec{w})\rightarrow \pi_1(Z_{0},p_0).
\]
is a well-defined surjective group homomorphism. 
\end{proposition}

\subsubsection{Iterated Integrals on the Central Fiber}

Since $A_0^\bullet$ is a sub-dga of $A^\bullet$ we may define:
\begin{definition}
Let $I=\sum_{|J|\le s} a_J\;
\varphi_{j_1}\otimes \cdots \otimes\varphi_{j_r}\in \bigoplus_{r=1}^s
\bigotimes^r A_0^1$ be Chen-closed and let $[\gamma]\in\pi_1(Z_0,p_0)$
be such that it is represented by a path $\gamma$ over $\vec{v}$ based at
$\vec{w}$. Then define
\[
\int_{[c(\gamma)]} I := \int_{\gamma} I.
\]
\end{definition}

Now let $J_0$ be the augmentation ideal in the group ring 
$\itgm\pi_1(Z_0,p_0)$. Similarly to Theorem \ref{Steinkauz} 
on page \pageref{Steinkauz} one proves:
\begin{theorem}
The integration of closed forms in $A^1_0$ along elements of $\pi_1(Z_0,p_0)$
defines an isomorphism
\[
H^1(A_0^\bullet)\xrightarrow{\cong} \Hom_\itgm\left(J_0\big/ J_0^2;\cpxm\right).
\]
\end{theorem}

Since by definition the integration maps are compatible with 
$c: \pi_1(Z_{\vec{v}},\vec{w})\rightarrow \pi_1(Z_{0},p_0)$, we have an 
inclusion
\[
H^1(A_0^\bullet)\hookrightarrow H^1(A^\bullet).
\]
Define for $R=\itgm$, $\ratm$, $\reellm$ or $\cpxm$: 
\[
H^1(A_0^\bullet)_R:=H^1(A^\bullet)_R\cap H^1(A_0^\bullet) 
\text{ and }
\left(J_0\big/ J_0^{s+1}\right)^*_R:=
\Hom_\itgm\left(J_0\big/ J_0^{s+1};R \right).
\]

It follows from the proof of Theorem \ref{Rothalsgans} on page 
\pageref{Rothalsgans} that holds:
\begin{theorem} 
\[
(J_{0}/ J_{0}^2)^*
:= \left(H^1(A_0^\bullet)_\itgm,\;(H^1(A_0^\bullet)_\ratm, W_\bullet),\;
(H^1(A_0^\bullet), W_\bullet, F^\bullet)\right) =
W_1(\vec{J}\big/\vec{J}^2)^*
\]
In particular, $(J_{0}/ J_{0}^2)^*$ is a $\itgm$-MHS of possible weights 
0 and 1.
\end{theorem}

Moreover, one can prove (similarly to Theorem \ref{Novalis}):
\begin{theorem}
Integration of Chen-closed elements
in $\bigoplus_{r=1}^s{\bigotimes}^s A_0^1$ along elements of 
$\pi_1(Z_0,p_0)$ defines an isomorphism of complex vector spaces
\[
H^0\bar{B}_s(A^\bullet_0,a)\xrightarrow{\cong} 
\Hom_\itgm\left(J_0\big/ J_0^{s+1};\cpxm\right). 
\] 
\qed
\end{theorem}

Since the integration maps are compatible with $c$ and since 
$\vec{J}\big/ \vec{J}^{s+1}\rightarrow J_0\big/ J_0^{s+1}$ is surjective, 
we find (the path $\sigma$ in the definition {\bf (d)} in 
\ref{Douwe Egberts}
does not contribute anything to the integral since there are no residues):
\begin{corollary}
$H^0\bar{B}_s(A^\bullet_0,a)\hookrightarrow H^0\bar{B}_s(A^\bullet,a)$. \qed
\end{corollary}

Similarly to Theorem \ref{Blaumeise} one proves:
\begin{theorem}
For any $s\ge 1$ is 
\[
\left( J_0\big/J_0^{s+1}\right)^*:=
\left(\left( J_0\big/J_0^{s+1}\right)^*_\itgm,\,
\left(\left( J_0\big/J_0^{s+1}\right)^*_\ratm, W_\bullet\right)\,
\left(\left( J_0\big/J_0^{s+1}\right)^*_\cpxm, W_\bullet, F^\bullet\right)
\right)
\]
a mixed Hodge structure and the short exact sequence 
\begin{equation} 
0\rightarrow \left( J_0\big/J_0^{s}\right)^*\rightarrow
\left( J_0\big/J_0^{s+1}\right)^*\rightarrow 
H^1(A_0^\bullet)^{\otimes s} \rightarrow 0
\end{equation}
is a short exact sequence of MHSs.
\end{theorem}

Finally, let us state the following easy to prove theorem.
\begin{theorem}
The map $c: \pi_1(Z_{\vec{v}},\vec{w})\rightarrow \pi_1(Z_{0},p_0)$
induces an inclusion of MHSs
\[
\big(J_0\big/ J_0^{s+1} \big)^*\hookrightarrow 
\big(\vec{J}\big/ \vec{J}^{s+1} \big)^*.
\] 
\qed
\end{theorem}


\section[Well-Definedness and Formal Connections]{Well-Definedness and Formal 
Connections}  \label{Placenta}
The main results of this section can be summarized in the following theorem. 
\footnote{The numbers in brackets refer to the place, where the 
respective assertions are proved.}
 
\begin{theorem} \label{Reinhard May}
Let $\chen$ be a Chen-closed element and let $\gamma:[a;b]\rightarrow
Z_0$ be a path over $\vec{v}$ based at $\vec{w}$. Then
\[
\lim_{\varepsilon\to 0} \eint{\gamma} 
\sum_J a_J\; \varphi_{j_1}\cdots\varphi_{j_r}
\]
converges (\ref{Babymassage}) and does not depend on the choice of 
coordinates (\ref{Vormilch}).
Moreover it does not depend upon the choice of $\gamma$ within
a nearby homotopy class of paths over $\vec{v}$ based at $\vec{w}$
(\ref{Blaek Foes}), 
i.~e.~we have a nearby homotopy functional:
\[
\int \sum_J a_J\; \varphi_{j_1}\cdots\varphi_{j_r}:
\pi_1(Z_{\vec{v}},\vec{w})\rightarrow \cpxm.
\]
Finally, this homotopy functional is zero if $I\in \R_1^s (A^\bullet, a)$.
\end{theorem}

In general on a manifold one can express iterated integrals
over Chen-closed linear combinations of tensor products of smooth
1-forms in three different ways as line integrals: First, one can pull
back the iterated integral to the parameter interval according to
the definition of iterated integrals.
Second, one can express the iterated integrals
as line integrals on the universal covering. And third, one can 
split up the iterated integral into pieces over parts of the path,
which lie in simply connected open sets. 
In all three cases one transforms the general iterated integral into
an iterated integral on a space, whose dga of differential forms
is acyclic. 

The idea to prove \ref{Reinhard May} is to express iterated integrals
along paths over $\vec{v}$ in a similar way locally as line integrals
along paths over $\vec{v}$. Away from the double points in $Z_0$ this
can be done like described above. But around a double point the
problem arises that if we restrict our dga $A^\bullet$ in the
naive way to a dga living on an open neighbourhood of the double point,
still the result will not be an acyclic complex. After all the dga 
$A^\bullet$ restricted to a double point computes the cohomology
of the vanishing cycle at that double point.

Since a path over $\vec{v}$ locally at a double point is forced to approach
the double point in a particular sector we can find 
simply connected open sets in $Z_0$, which contain
all of this (local) path apart from the double point itself. If we are
given a closed form $\varphi$ in $A^1$ then we can find primitives of this
form at least on those open sets. But since $\varphi$ may have a simple pole
in the double point,
these primitives can have logarithmic singularities there.

Recall that reducing the length of an iterated integral
over a Chen-closed element with differential forms in an acyclic complex
can be done as follows: Choose a primitive of a closed form and 
multiply this function
in an appropriate way with the other forms (cf.~\cite{Hain-the-geometry},
Prop.~(1.3), p.~252).
When one takes an iterated integral over a composition of two 
paths, one can take
primitives of this closed form over each of the paths and has to take
care that their values on the junction point match.

For instance in the case of a path which passes only through
one double point $p_{kl}$ and changes from one component to another,
we would like to
decompose this path into three parts: two parts where it stays in one of
the components and a third part which is represented by
the 1-simplex $\Delta^1$.
(Observe that also in the definition (a) on page \pageref{Kormoran} we
thought of $\Delta^1$ as being part of the path $\gamma$.)
But, when we now choose primitives on each of the components of a
closed compatible 1-form, in general these primitives
have logarithmic singularities at the junction
points of the three parts
of the path. Hence, we cannot make their values match. It will
turn out that although we cannot make their values match, we can
make their formal shapes or the type of their singularity match.
This leads to a kind of {\it higher compatibility} and finally
to the definition of the dga $A^\bullet_{kl}$.
In this way we can then indeed locally express an iterated integral
along a path $\gamma$ over $\vec{v}$ as line integral of a closed
form in $A^1_{kl}$ along $\gamma$. This construction 
allows us to prove in the local case 
convergence, independence of local coordinates around the double points
and invariance under nearby homotopies between paths over $\vec{v}$
of iterated integrals of Chen-closed elements.

A crucial role in all convergence considerations is played by
the following classical fact (cf.~\cite{Courant} III.2, p.~192).

\begin{fact} \label{Schwan}
For $\varepsilon_0>0$ let $h:]-\varepsilon_0,\varepsilon_0[\rightarrow
\reellm$ be a smooth function that vanishes at 0.
Then for any $m\in\natm$ holds:
$\lim_{\varepsilon\searrow 0} h(\varepsilon) \log^m \varepsilon=0.$
\end{fact}

Since iterated integrals are not additive on paths in the naive sense,
these  results do not carry over immediately to the corresponding
assertions in the global situation. Here we first have to introduce 
an analogue to Chen's theory of {\it formal power series connections}
in order to describe the relation between iterated integrals on a
composition of two paths over $\vec{v}$ and iterated integrals over 
these paths seperately.

\subsection{Local Considerations}

We start out by considering iterated integrals along paths which are 
defined locally in a neighbourhood of a double point.

\subsubsection{A Localization of $A^\bullet$ 
in a Sector at a Double Point} \label{Tommi Engel}

Consider a double point $p_{kl}$ in $Z_0$ with $[k<l]$.
Recall that there are coordinates
$(x,y):W_{kl}=U_{kl}^k\times U_{kl}^l\rightarrow\cpxm^2$ around
this double point. Using these coordinates, we want to construct a 
dga $A^\bullet_{kl}$ for this double point $p_{kl}$.

Define $U_{kl}^{k-}:=U_{kl}^{k}\setminus \reellm^{\le 0}$
and $U_{kl}^{l-}:=U_{kl}^{l}\setminus \reellm^{\le 0}$. Note that
the function $\log x$ (resp.~$\log y$) can be defined univalently on
$U_{kl}^{k-}$ (resp.~$U_{kl}^{l-}$). Now let $E^\bullet(U_{kl}^{k-}\log p_{kl})$
(resp.~$E^\bullet(U_{kl}^{l-}\log p_{kl})$) be the sub-dga of
$E^\bullet(U_{kl}^{k-})$ (resp.~$E^\bullet(U_{kl}^{l-})$), which is generated
as $E^\bullet(U_{kl}^{k})$ (resp.~$E^\bullet(U_{kl}^{l})$)--modules by powers of
$\log x$ (resp.~$\log y$) and their derivatives. That is, we have
embeddings of dga's
\begin{align*}
E^\bullet(U_{kl}^{k})\otimes_\cpxm {\bigwedge}^\bullet (\frac{dx}{x})[\log x]
&\hookrightarrow E^\bullet(U_{kl}^{k-})\\
\intertext{and}
E^\bullet(U_{kl}^{l})\otimes_\cpxm {\bigwedge}^\bullet (\frac{dy}{y})[\log y]
&\hookrightarrow E^\bullet(U_{kl}^{l-}),
\end{align*}
whose images are $E^\bullet(U_{kl}^{k-}\log p_{kl})$ and
$E^\bullet(U_{kl}^{l-}\log p_{kl})$. \footnote{These maps are embeddings
for the following reason: If $\sum_{\nu=0}^m g_\nu(x) \log^\nu x \equiv 0$
with $g_\nu\in E^0(U_{kl}^k)$ it follows that all $g_\nu\equiv 0$
for $\nu=0,\ldots,m$, because for any $x_0\in U_{kl}^k$ the polynomial
$\sum_{\nu=0}^m g_\nu(x_0) w^\nu \in\cpxm[w]$ has infinitely many zeroes:
$\{\log x_0 +2\pi i n | n\in\itgm\}$.}

Since $E^\bullet(U_{kl}^{k})$ (resp.~$E^\bullet(U_{kl}^{l})$) and
${\bigwedge}^\bullet (\frac{dx}{x})[\log x]$
(resp.~${\bigwedge}^\bullet (\frac{dy}{y})[\log y]$) have cohomology only
in degree 0, the K{\"u}nneth formula tells us that also
$E^\bullet(U_{kl}^{k-}\log p_{kl})$ (resp.~$E^\bullet(U_{kl}^{l-}\log p_{kl})$)
have no cohomology except in degree 0.
If $\{k,\,l\}=\{0,\,l\}$ then define
\[
E^\bullet(U_{0l}^{0-}\log p_{01}):={\bigwedge}^\bullet (\frac{dp}{p})[\log p],
\]
which are the polynomials in the variable $\log p$ and coefficients
in ${\bigwedge}^\bullet (\frac{dp}{p})$, where the differential
is given by the Leibniz rule, $d\,\frac{dp}{p}=0$ and
$d\,\log p=\frac{dp}{p}$.

Now let $\Lambda_{kl}^\bullet:=
{\bigwedge}^\bullet(\frac{dx}{x},\frac{dy}{y})[v,w]$	
be the dga, which is defined as follows: 
Let $\Lambda^\bullet:={\bigwedge}^\bullet(\frac{dx}{x},\frac{dy}{y})[v,w]$
be the dga ${\bigwedge}^\bullet(\frac{dx}{x},\frac{dy}{y})[v,w]$
of polynomials in two variables $v$ and $w$ with coefficients
in ${\bigwedge}^\bullet(\frac{dx}{x},\frac{dy}{y})$, where the differential
$d$ is given by the Leibniz rule,
$d\left(\frac{dx}{x}\right)=d\left(\frac{dy}{y}\right)=0$ and
\[
d\,v:=\frac{dx}{x}\qquad\text{	 and	}\qquad d\, w := \frac{dy}{y}.
\]
Observe that there is an embedding
\[
{\bigwedge}^\bullet(\frac{dx}{x},\frac{dy}{y})[u]\hookrightarrow
{\bigwedge}^\bullet(\frac{dx}{x},\frac{dy}{y})[v,w]
\]
by setting $u=v+w$.
Note that $(\Lambda_{kl}^\bullet,\,d)$ only has cohomology
in degree 0.

Similar like in the definition of $A^\bullet$ we define the
dga $A_{kl}^\bullet$ as sub-dga of
\[
B_{kl}^\bullet := E^\bullet(U_{kl}^{k-}\log p_{kl})\oplus
E^\bullet(U_{kl}^{l-}\log p_{kl}) \oplus
E^\bullet(\Delta^1)\otimes_\cpxm\Lambda_{kl}^\bullet,
\]
which contains all those elements that are subject to the following
compatibility conditions generalizing the compatibility conditions
in the definition of $A^\bullet$. Observe that $B_{kl}^n=0$ for $n\ge 4$.

We shall use $P_{kl}$, $K_{kl}$, $L_{kl}$, $H_{kl}$,
$R_{kl}$, $S_{kl}$ and $T_{kl}$ to denote elements in
$\cpxm[\xi,v,w]$. If the context is clear, we omit the indices $k,l$.

\noindent
$\bf A_{kl}^0 :$
An element
\[
f=\sum_{\nu=0}^m g_{k,\nu} (x) \log^\nu x
+\sum_{\nu=0}^m g_{l,\nu} (y) \log^\nu y  + P_{kl}(\xi,v,w)
\in B_{kl}^0
\]
is called {\it compatible}, iff holds:
\[
P_{kl}(0,v,w)=\sum_{\nu=0}^m g_{k,\nu} (0) v^\nu \text{ and }
P_{kl}(1,v,w)=\sum_{\nu=0}^m g_{l,\nu} (0) w^\nu.
\]

\noindent
$\bf A_{kl}^1 :$
We call an element in $B_{kl}^1$,
\[{
\varphi=\sum\limits_{\nu= 0}^m \omega_{k,\nu}\log^\nu x +
\sum\limits_{\nu= 0}^m \omega_{l,\nu}\log^\nu y
+ K\,\frac{dx}{x}
+L\,\frac{dy}{y} + H \, d\xi
}\]
a {\it compatible} element iff it satifies:
\begin{eqnarray*}
K(0,v,w) = &\sum_{\nu= 0}^m \Res_{p_{kl}}\omega_{k,\nu}\,v^\nu
&,\quad L(0,v,w)=0, \\
K(1,v,w) = & 0 &,\quad L(1,v,w)=\sum_{\nu= 0}^m
\Res_{p_{kl}}\omega_{l,\nu}\,w^\nu .
\end{eqnarray*}

\noindent
$\bf A_{kl}^2 :$
Let us call an element
\[
\phi = \Omega^{(k)} + \Omega^{(l)}
+ R \,d\xi\wedge\frac{dx}{x}
+ S\, d\xi\wedge\frac{dy}{y} + T\,
\frac{dx}{x}\wedge\frac{dy}{y} \in B_{kl}^2
\]
{\it compatible} iff holds:
\[
T(0,v,w)=T(1,v,w)=0.
\]

\noindent
$\bf A_{kl}^3 :$
And finally let $A_{kl}^3:=B_{kl}^3$.

Note that $dA_{kl}^i\subset A_{kl}^{i+1}$ for $i\ge 0$ and that
$A_{kl}^\bullet$ is a dga.
Again, similar like in the definition of $A^\bullet$, we can give an
alternative definition of $A_{kl}^\bullet$ as follows.
Consider the surjective map of complexes
\[
\Phi_{kl}: B_{kl}^\bullet \longrightarrow
\Lambda_{kl}^\bullet\oplus\Lambda_{kl}^\bullet,
\]
which sends  the element $\sum_{\nu=0}^m g_{k,\nu} (x) \log^\nu x
+\sum_{\nu=0}^m g_{l,\nu} (y) \log^\nu y  + P_{kl}(\xi,v,w) \in B_{kl}^0$
to $\left(P_{kl}(0,v,w)-\sum_{\nu=0}^m g_{k,\nu} (0) v^\nu\right)\oplus
\left(P_{kl}(1,v,w)-\sum_{\nu=0}^m g_{l,\nu} (0) w^\nu\right)$ and
\begin{gather*} 
\sum_{\nu= 0}^m \omega_{k,\nu}\log^\nu x +
\sum_{\nu= 0}^m \omega_{l,\nu}\log^\nu y
+ K(\xi,v,w)\frac{dx}{x}
+L(\xi,v,w)\frac{dy}{y} + H(\xi,v,w) d\xi 
\end{gather*}
\begin{align*}
\text{in $B_{kl}^1$ to} \qquad\qquad\qquad 
 &\left[\left(K(0,v,w)-\sum_{\nu= 0}^m
\Res_{p_{kl}}\omega_{k,\nu}\,v^\nu\right)\frac{dx}{x}
+L(0,v,w)\frac{dy}{y} \right] \\
\oplus
&\left[K(1,v,w)\frac{dx}{x}+\left(L(1,v,w)-\sum_{\nu= 0}^m
\Res_{p_{kl}}\omega_{l,\nu}\,w^\nu
\right)\frac{dy}{y} \right]. 
\end{align*}
Moreover, $\Phi_{kl}$ maps
\[
\Omega^{(k)} + \Omega^{(l)} + R(\xi,v,w) \,d\xi\wedge\frac{dx}{x}
+S(\xi,v,w)\, d\xi\wedge\frac{dy}{y} + T(\xi,v,w)\,
\frac{dx}{x}\wedge\frac{dy}{y} \in B_{kl}^2
\]
to $ T(0,v,w) \frac{dx}{x}\wedge\frac{dy}{y} \oplus
T(1,v,w) \frac{dx}{x}\wedge\frac{dy}{y}$. It is easy to check that $\Phi$ is
indeed a surjective map of complexes. Then we find:
\[
A_{kl}^\bullet = \ker \Phi_{kl}
\]
and with $C_{kl}^\bullet:= \Lambda_{kl}^\bullet\oplus\Lambda_{kl}^\bullet$
we get the short exact sequence
\[
0\longrightarrow A_{kl}^\bullet\longrightarrow B_{kl}^\bullet\longrightarrow
C_{kl}^\bullet \longrightarrow 0
\]
and may conclude from the long exact cohomology sequence using
$H^0(A^\bullet_{kl})=\cpxm$ that $H^i(A^\bullet_{kl})=0$ for $i\ge 1$.

Note that there is an obvious dga-morphism: 
$A^\bullet\rightarrow A^\bullet_{kl}$.

\subsubsection{A Path Locally at a Double Point}

Consider a path over $\vec{v}$, which lies
entirely in $U_{kl}^{k-}\cup\{p_{kl}\}\cup U_{kl}^{l-}$ for a double point
$p_{kl}$. We want to extend the definitions (a), (b), (c) and (d) in 
\ref{Douwe Egberts} to the definition of
$\varepsilon$-iterated integrals
with forms in $A^1_{kl}$ along this path. Again like before we
define these $\varepsilon$-iterated integrals in several steps.
At the end we let the iterated integral of a Chen-closed element
with these forms along such a path be the limit of the sum of
$\varepsilon$-iterated integrals.

Let $\varphi_1\otimes\cdots\otimes\varphi_r$ in $\bigotimes^r A^1_{kl}$,
where we write each
\begin{gather*}
\varphi_j=\tilde{\omega}_k^{(j)}+\tilde{\omega}_l^{(j)} + \Xi^{(j)} \\
\intertext{with}
\tilde{\omega}_k^{(j)}=\sum_{\nu=0}^m \omega_{k,\nu}^{(j)}\,
\log^\nu x \qquad \text{  and	}\qquad\tilde{\omega}_l^{(j)}=\sum_{\nu=0}^m
\omega_{l,\nu}^{(j)}\,\log^\nu y \\
\intertext{and}
\Xi^{(j)}=K^{(j)}\frac{dx}{x} +L^{(j)}\frac{dy}{y}
+H^{(j)}d\xi.
\end{gather*}
We define now the iterated integral of
$\varphi_1\otimes\cdots\otimes\varphi_r$ along a path over $\vec{v}$
in $U_{kl}^{k-}\cup\{p_{kl}\}\cup U_{kl}^{l-}$ as follows.

\begin{description}
\item[(a)] \label{Pekingente}
Let $\gamma:[a;b]\rightarrow U_{kl}^{k-}\cup\{p_{kl}\}\cup U_{kl}^{l-}$ be
a path over $\vec{v}$, which meets the set
of double points only once with parameter value $\tau_0\in ]a;b[$, where
it changes from $D_k$ to $D_l$.
Define $\gamma_x(\tau)$ (resp.~$\gamma_y(\tau)$) to be $x(\gamma(\tau))$
(resp.~$y(\gamma(\tau))$).
We abbreviate for $\varepsilon>0$ small enough:
\[
\eta_\varepsilon^{(j)}:=H^{(j)}\big(\xi,\log\left(
\gamma_x(\tau_0-\varepsilon)
\cdot \gamma_y (\tau_0+\varepsilon)\right)\big)\; d\xi.
\]
Then define
\[
\eint{\gamma} \!\! \varphi_1\cdots \varphi_r
\!:=\!\!\!\!\!\!\! \sum_{0\le\alpha\le\beta\le r}\;
\int\limits_{\gamma^{\le \tau_0-\varepsilon}} \!\!\!\!\!\!
\tilde{\omega}_k^{(1)}\cdots\tilde{\omega}_k^{(\alpha)}
\!\!\!\int\limits_{[0;1]} \!\eta_\varepsilon^{(\alpha+1)}\cdots
\eta_\varepsilon^{(\beta)}
\!\!\!\!\!\!\!\int\limits_{\gamma_{\ge \tau_0+\varepsilon}} \!\!\!\!\!
\tilde{\omega}_l^{(\beta+1)}\cdots \tilde{\omega}_l^{(r)}.
\]
\item[(b)]
Let $\gamma:[a;b]\rightarrow U_{kl}^{k-}\cup\{p_{kl}\}\cup U_{kl}^{l-}$
be a path over $\vec{v}$, which meets the set
of double points only once with parameter value $\tau_0\in ]a;b[$, where
it stays in one component $D_k$. Then define:
\[
\eint{\gamma} \varphi_1\cdots \varphi_r:=\sum_{0\le \alpha \le r}\;
\int\limits_{\gamma^{\le \tau_0-\varepsilon}}
\tilde{\omega}_k^{(1)}\cdots\tilde{\omega}_k^{(\alpha)}
\int\limits_{\gamma_{\ge \tau_0+\varepsilon}}
\tilde{\omega}_k^{(\alpha+1)}\cdots\tilde{\omega}_k^{(r)}.
\]
\item[(c)]
Now let $\gamma:[a;b]\rightarrow U_{kl}^{k-}\cup\{p_{kl}\}\cup U_{kl}^{l-}$
be a path over $\vec{v}$
such that $\gamma(a)$ and $\gamma(b)$ are different from $p_{kl}$.
Then there is a finite number of parameter values $\tau_1,\ldots,\tau_N$
with $a<\tau_1<\cdots<\tau_N<b$, which are mapped onto $p_{kl}$.
Choose a $\tau_i^*\in ]\tau_i;\tau_{i+1}[$ for $i=1,\ldots,N-1$ and let
$\tau_0^*:=a$ and $\tau_N^*:=b$. Then we define $\gamma_i:=
\gamma_{\big|[\tau_{i-1}^*,\tau_i^*]}$ for $i=1,\ldots,N$ and let:
$\eintt{\gamma} \varphi_1\cdots \varphi_r:=$
\[
\sum_{0\le\alpha_1\le\cdots\le\alpha_{N-1}\le r}\;
\eint{\gamma_{1}}
\varphi_1\cdots \varphi_{\alpha_1}
\eint{\gamma_{2}}
\varphi_{\alpha_1+1}\cdots \varphi_{\alpha_2}\;\cdots\;
\eint{\gamma_{N}}
\!\varphi_{\alpha_{N-1}+1}\cdots \varphi_{\alpha_r}.
\]
\item[(d)]
Finally assume that $p_{kl}=p_0$. Let 
$\gamma:[a;b]\rightarrow \{p_{01}\}\cup U_{01}^{1-}$
be a path over $\vec{v}$ starting (or ending) at $\vec{w}$. 
Like before, the coordinate $p:D_0\rightarrow\cpxm$ allows us to 
consider $D_0$ as part of the complex plane. There, in $\cpxm$,
we have the differential form (and not the symbol) $\frac{dp}{p}$.
Let $\sigma$ be the straight
path $\tau\mapsto (1-\tau)$ in $\cpxm$, then we define 
$\varepsilon>0$ small enough 
\[
\eint{\gamma} \varphi_1\cdots \varphi_r
:=\eint{\sigma\star\gamma}
\varphi_1\cdots \varphi_r\quad \left(\text{ or }
\eint{\gamma\star\sigma^{-1}}
\varphi_1\cdots \varphi_r\right).
\]
\end{description}

Note that if all the $\varphi_j$ are even in the image of the obvious
dga morphism $A^\bullet\rightarrow A^\bullet_{kl}$, then clearly
these definitions coincide with the definitions (a), (b), (c) and
(d) in \ref{Douwe Egberts}.

\begin{proposition} \label{Kafka}
Let $f\in A^0_{kl}$. Then for any path $\gamma:
[a;b]\rightarrow U_{kl}^{k-}\cup\{p_{kl}\}\cup U_{kl}^{l-}$ over
$\vec{v}$ holds:
\[
\lim_{\varepsilon\to 0}\eint{\gamma} df = f(\gamma(b))- f(\gamma(a)).
\]
\end{proposition}

\noindent
{\bf Proof:}
Let $\varphi=\tilde{\omega}_k+\tilde{\omega}_l +
K\frac{dx}{x} +L\frac{dy}{y}+Hd\xi$ $=df=d\left(G_{k}(x) + G_l(y)
+P(\xi,v,w)\right)$ and write
\[
\tilde{\omega}_k=\sum_{\nu=0}^m \omega_{k,\nu}\,
\log^\nu x \qquad\text{  and  } \qquad\tilde{\omega}_l=\sum_{\nu=0}^m
\omega_{l,\nu}\,\log^\nu y
\]
as well as
\[
G_k(x)=\sum_{\nu=0}^m g_{k,\nu}(x)\,
\log^\nu x \qquad\text{  and  }\qquad G_l=\sum_{\nu=0}^m
g_{l,\nu}(y)\,\log^\nu y.
\]

We can decompose $\gamma$ into paths $\gamma_i$, for $i=1,\ldots,N$,
i.~e.~$\gamma=\gamma_1\star\cdots\star\gamma_N$ such that each
of the paths $\gamma_i$ meets the double point $p_{kl}$ only once.
Because line integrals are additive we have
\[
\int_\gamma \varphi = \sum_{i=1}^N \int_{\gamma_i} \varphi,
\]
provided that all the integrals on the right hand side converge.
Therefore we may assume without loss of generality that $\gamma$
meets the double point $p_{kl}$ only once with parameter value
$\tau_0\in]a;b[$. We distinguish two cases.

\noindent
{\bf A path traversing the double point:}
Suppose that $\gamma$ changes components, say, it runs from $D_k$ to $D_l$.
Observe first that holds:
\[
\int_{[0;1]} H(\xi,v,w) d\xi = P(1,v,w) - P(0,v,w).
\]
Then we may compute -- using the compatibility:
\begin{align*}
&\int_{\gamma^{\le \tau_0-\varepsilon}} d G_k -	P(0,
\log \gamma_x(\tau_0-\varepsilon),\log \gamma_y(\tau_0+\varepsilon)) \\
= &G_k(\gamma_x(\tau_0-\varepsilon))
-\sum_{\nu=0}^m g_{k,\nu}(0)\log^\nu \gamma_x(\tau_0-\varepsilon)
-G_k(\gamma_x(a)) \\
= &\sum_{\nu=0}^m
\left(g_{k,\nu}(\gamma_x(\tau_0-\varepsilon))-
g_{k,\nu}(0)\right)\log^\nu \gamma_x(\tau_0-\varepsilon)
-G_k(\gamma_x(a)).
\end{align*}
With Fact \ref{Schwan} we derive
\begin{align*}
&\lim\limits_{\varepsilon\to 0}
\int\limits_{\gamma^{\le \tau_0-\varepsilon}}  \sum\limits_{\nu=0}^m
\omega_{k,\nu}\log^\nu  x - P(0,
\log \gamma_x(\tau_0-\varepsilon),\log \gamma_y(\tau_0+\varepsilon))
=-G_k(\gamma_x(a)) \\
\intertext{and similarly we find:}
&\lim\limits_{\varepsilon\to 0}
\int\limits_{\gamma^{\ge \tau_0+\varepsilon}}  \sum\limits_{\nu=0}^m
\omega_{l,\nu}\log^\nu y - P(1,
\log \gamma_x(\tau_0-\varepsilon),\log \gamma_y(\tau_0+\varepsilon))
= G_l(\gamma_y(b)).
\end{align*}

\noindent
{\bf A path colliding with the double point:}
Suppose now that $\gamma$ stays in one component $D_k$.
Here we compute:
\[
\int\limits_{\gamma^{\le \tau_0-\varepsilon}} \!\!\!\tilde{\omega}_k
+\!\!\!\int\limits_{\gamma^{\ge \tau_0+\varepsilon}}\!\!\!\tilde{\omega}_k
=\left(G_k(\gamma_x(\tau_0-\varepsilon))-G_k(\gamma_y(\tau_0+\varepsilon))
\right)+\left(G_k(\gamma_x(b))-G_k(\gamma_x(a))\right)
\]
and the first summand tends to 0 as $\varepsilon$ becomes small.
\qed

\bigskip
Before we formulate and prove a generalization of the preceeding 
Proposition to iterated integrals of arbitrary length we 
need some preparation. 
To perform all the convergence considerations in the proof of
Theorem \ref{Brentano} we need the following technical but
elementary Lemma.

\begin{lemma} \label{Flamingo}
Let $\gamma:[a;b]\rightarrow U_{kl}^{k-}\cup \{p_{kl}\}\cup
U_{kl}^{l-}$ be a path over $\vec{v}$ that meets the double point
once with parameter value $\tau_0$ and changes from $D_k$ to $D_l$.
Let $\tilde{\omega}_k^{(1)},\ldots,\tilde{\omega}_k^{(r)}\in
E^1(U_{kl}^{k-}\log p_{kl})$ and let $\tilde{\omega}_l^{(1)},\ldots,
\tilde{\omega}_l^{(r)}\in E^1(U_{kl}^{l-}\log p_{kl})$. Define for
the polynomials $h^{(1)},\ldots,h^{(r)}\in \cpxm[\xi,v,w]$:
\[
\eta^{(j)}_\varepsilon:= h^{(j)}\big(\xi,\log \gamma_x(\tau_0-\varepsilon),
\log \gamma_y(\tau_0+\varepsilon)\big)d\xi.
\]
Then there is a positive number $C$ and a natural number $N$ such that for
all $\varepsilon>0$ small enough holds:
\begin{eqnarray}
\left|\int_{\gamma^{\le \tau_0-\varepsilon}}
\tilde{\omega}_k^{(1)}\cdots\tilde{\omega}_k^{(r)}\right|
&\le& C \left| \log^N \varepsilon \right|  \label{Wim} \\
\left|\int_{[0;1]}
\eta_\varepsilon^{(1)}\cdots\eta_\varepsilon^{(r)}\right|
& \le & C \left| \log^N \varepsilon \right| \label{Wam}\\
\left|\int_{\gamma_{\ge\tau_0-\varepsilon}}
\tilde{\omega}_l^{(1)}\cdots\tilde{\omega}_l^{(r)}\right|
&\le& C \left| \log^N \varepsilon \right|. \label{Wum}
\end{eqnarray}
Also the corresponding assertion holds for a
path $\gamma$ over $\vec{v}$, which
meets the double point once but stays in one component.
\end{lemma}

\noindent
{\bf Proof:} For the proof of \eqref{Wam} in Lemma \ref{Flamingo},
observe that
\[
\int_{[0;1]}\left(h^{(1)}(\xi,v,w)d\xi\right)\cdots
\left(h^{(r)}(\xi,v,w)d\xi\right)
\]
is a polynomial, say $Q(v,w)$, in $\cpxm[v,w]$ such that
\[
\int_{[0;1]}\eta_\varepsilon^{(1)}\cdots\eta_\varepsilon^{(r)}=
Q(\log \gamma_x(\tau_0-\varepsilon),\log \gamma_y(\tau_0+\varepsilon)).
\]
Write $Q(v,w)=\sum_{m,n} a_{mn} v^m w^n$. Then we find
$Q(\log \gamma_x(\tau_0-\varepsilon),\log \gamma_y(\tau_0+\varepsilon))$
\[
=\sum_{m,n} a_{mn}\left(\log (\frac{\gamma_x(\tau_0-\varepsilon)}{\varepsilon})
+\log\varepsilon\right)^m
\left(\log (\frac{\gamma_y(\tau_0+\varepsilon)}{\varepsilon})
+\log\varepsilon\right)^n.
\]
Assertion \eqref{Wam} is true since:
\[
\lim_{\varepsilon\to 0}\log
\left(\frac{\gamma_x(\tau_0-\varepsilon)}{\varepsilon}\right)
=-\log \dot{\gamma}_x(\tau_0) \quad \text{  and  } \quad
\lim_{\varepsilon\to 0}\log
\left(\frac{\gamma_y(\tau_0+\varepsilon)}{\varepsilon}\right)
=\log \dot{\gamma}_y(\tau_0).
\]
\eqref{Wum} and \eqref{Wim} are quite similar to prove.
We show \eqref{Wim} and leave it to the reader to prove \eqref{Wum}.
Let $f_j(\tau) d\tau=\gamma_x^*\tilde{\omega}_k^{(i)}$
for $i=1,\ldots,r$. Then
\begin{eqnarray*}
 & &\left|\;\int\limits_{\gamma^{\le \tau_0-\varepsilon}}
\tilde{\omega}_k^{(1)}\cdots\tilde{\omega}_k^{(r)}\right| \\
&\le&\left| \!\!\!\int\limits_{\qquad
a\le\tau_1\le \cdots\le\tau_r\le\tau_0-\varepsilon}
\!\!\!\!\!\!\!\!\!\!
\!\!\!\!\!\!\!\!\!\!\cdots\int f_1(\tau_1)\cdots f_r(\tau_r)
d\tau_1\cdots d\tau_r\right|  \\
&\le& \int_a^{\tau_0-\varepsilon}
|f_1(\tau_1)|d\tau_1 \cdot \cdots\cdot
\int_a^{\tau_0-\varepsilon} |f_r(\tau_r)| d\tau_r
\end{eqnarray*}
and if	$\tilde{\omega}_k{(j)}=\sum_{\nu=0}^m A_\nu(x)\log^\nu x
\frac{dx}{x} + B_\nu(x)\log^\nu x d\bar{x}$, then we find
\[
f_j(\tau)=\sum_{\nu=0}^m\left\{
\frac{A_\nu(\gamma_x(\tau))\log^\nu\gamma_x(\tau)}{\gamma_x(\tau)}
\dot{\gamma}_x(\tau)
+{B_\nu(\gamma_x(\tau))
\log^\nu\gamma_x(\tau)}\overline{\dot{\gamma}_x(\tau_)}\right\}.
\]
There are $A_\nu$, $B_\nu>0$ such that $\sup_x |A_\nu(x)|\le A_\nu$
and $\sup_x |B_\nu(x)|\le b_\nu$. Moreover let $C_\nu:=\int_a^{\tau_0}
\left|\overline{\dot{\gamma}_x(\tau)}\right| d\tau$. Then
\begin{eqnarray*}
&&\int_a^{\tau_0-\varepsilon}|f_j(\tau)|d\tau \\
& \le  &\sum\limits_{\nu=0}^m\left\{
A_\nu\int\limits_a^{\tau_0-\varepsilon}
\left|\frac{\log^\nu\gamma_x(\tau)}{\gamma_x(\tau)}
\dot{\gamma}_x(\tau_)\right|d\tau
+B_\nu\int_a^{\tau_0-\varepsilon}
\left|\log^\nu\gamma_x(\tau)\overline{\dot{\gamma}_x(\tau_)}\right| \right\}\\
& \le & \sum\limits_{\nu=0}^m\left\{
A_\nu \left|\frac{1}{\nu+1}\left(\log^{\nu+1}\gamma_x(\tau_0-\varepsilon)
-\log^{\nu+1}\gamma_x(a)\right)\right| \right.
+\left. B_\nu\left|\log^\nu\gamma_x(\tau_0-\varepsilon)\right| C_\nu\right\}
\end{eqnarray*}
The assertion follows now from the differentiability of $\gamma_x$.
\qed

\bigskip
The following proposition is an analogue to the corresponding laws 
for iterated integrals on manifolds. There it is an easy
consequence of partial integration. Here we have to take the
limit $\varepsilon\to 0$ into account. 

\begin{proposition} \label{Erstlingsausstattung}
Assume now that $\gamma:[a;b]\rightarrow
U_{kl}^{k-}\cup \{p_{kl}\} \cup U_{kl}^{l-}$ is a path over $\vec{v}$.
Let $\varphi_1,\ldots,\varphi_s\in A^1_{kl}$ be such that at least one of them
is exact.
\begin{align*}
&\text{If $\varphi_1=df_1$, then} \\ 
&{\textstyle\lim\limits_{\varepsilon\to 0}\left\{
\eintt{\gamma}df_1{\varphi}_2\ldots \varphi_s -\eintt{\gamma}
\big((f_1-f_1(\gamma(a))\varphi_2\big)\varphi_3\cdots\varphi_s
\right\} =0,} \\ 
&\text{if $\varphi_\lambda =df_\lambda$ for $2\le \lambda\le s$, then}
\\
&{\textstyle\lim\limits_{\varepsilon\to 0}\left\{
\eintt{\gamma} \varphi_1\!\cdots\! df_\lambda\cdots\varphi_s -
\eintt{\gamma} \varphi_1\!\cdots\!(f_\lambda\varphi_{\lambda+1})\!\cdots\!
\varphi_s + \eintt{\gamma} \varphi_1
\!\cdots\!(f_\lambda\varphi_{\lambda-1})\!\cdots\!
\varphi_s\!\right\}=0,} \\
&\text{if $\varphi_s=df_s$, then} \\ 
&{\textstyle\lim\limits_{\varepsilon\to 0}\left\{
\eintt{\gamma} \varphi_1\cdots\varphi_{s-1} df_s -
\eintt{\gamma} \varphi_1\cdots\varphi_{s-2}\big(f_s
-f_s(\gamma(b))\varphi_{s-1}\big)\right\}=0.} 
\end{align*}
\end{proposition}

\noindent
{\bf Proof:}
The computations for the proof of these laws are straightforward but tedious.
The idea is to express the iterated integrals of the proposition in terms of
iterated integrals along paths lying completely in either one of the 
components or on the 1-simplex $\Delta^1$. For the latter iterated integrals   
the corresponding laws hold. The complete calculations can be found in
\cite{Doktorarbeit}, pp.~121-129.

Let us here just prove the first of the three laws. The proof of the others
is comparable. Write 
\begin{gather*}
\varphi_j=\tilde{\omega}_k^{(j)}+\tilde{\omega}_l^{(j)} +
K^{(j)}\frac{dx}{x} +L^{(j)}\frac{dy}{y}+H^{(j)}d\xi \\
\intertext{and} 
f_1=f=G_{k}(x) + G_l(y)+P(\xi,v,w), \\
\intertext{where} 
\tilde{\omega}_k=\sum_{\nu=0}^m \omega_{k,\nu}\,
\log^\nu x \qquad\text{  and  } \qquad\tilde{\omega}_l=\sum_{\nu=0}^m
\omega_{l,\nu}\,\log^\nu y  \\
\intertext{and}
G_k(x)=\sum_{\nu=0}^m g_{k,\nu}(x)\,
\log^\nu x \qquad\text{  and  }\qquad G_l(y)=\sum_{\nu=0}^m
g_{l,\nu}(y)\,\log^\nu y.
\end{gather*}
Before we show the first law 
for general paths $\gamma$ over $\vec{v}$
in $U_{kl}^{k-}\cup \{p_{kl}\} \cup U_{kl}^{l-}$ we consider
the local situation at a double point first.

\bigskip
\noindent
{\bf A path traversing the double point once:}
Suppose that $\gamma:[a;b]\rightarrow
U_{kl}^{k-}\cup \{p_{kl}\} \cup U_{kl}^{l-}$ meets the double point $p_{kl}$
once with parameter value$\tau_0\in]a;b[$, where it changes from $D_k$ to
$D_l$.
Again, we make use of the notation:
\begin{gather*}\
P_\varepsilon(\xi):= P(\xi,\log\gamma_x(\tau_0-\varepsilon),
\log\gamma_y(\tau_0+\varepsilon)) \\
\intertext{ and }
\eta^{(j)}_\varepsilon:= H^{(j)}(\xi,\log\gamma_x(\tau_0-\varepsilon),
\log\gamma_y(\tau_0+\varepsilon))d\xi.
\end{gather*}
Note that: $d P_\varepsilon(\xi)=\eta^{(\lambda)}_\varepsilon$ in
$E^\bullet(\Delta^1)$. Assume without loss
of generality that $f(p)=f(\gamma(a))=0$. Then we compute
with $\varphi_1=df$:

\newcommand{\termk}[1]{\int\limits_{\gamma^{\le \tau_0-\varepsilon}}
\!\!\!\!\tilde{\omega}_k^{(#1)}\cdots\tilde{\omega}_k^{(\alpha)}}
\newcommand{\termx}[1]{\!\!\!\! \int\limits_{[0;1]}
\!\!\!\!\eta_\varepsilon^{(#1)}\!\!\!\!\cdots\eta_\varepsilon^{(\beta)}}
\newcommand{\terml}[1]{\!\!\!\! \int\limits_{\gamma_{\ge \tau_0+\varepsilon}}
\!\!\!\!\tilde{\omega}_l^{(#1)}\!\!\!\!\cdots\tilde{\omega}_l^{(s)}}

\begin{align*}
  & \eint{\gamma} (f\varphi_2)\varphi_3\cdots\varphi_s \displaybreak[0] \\
= & \int\limits_{\gamma_{\ge \tau_0+\varepsilon}}
  (G_l\tilde{\omega}_l^{(2)})\tilde{\omega}_l^{(3)}
  \cdots\tilde{\omega}_l^{(s)}\\
+ & \sum\limits_{1=\alpha<\beta\le s} \;\;
  \int\limits_{[0;1]}
  (P_\varepsilon\eta_\varepsilon^{(2)})\eta_\varepsilon^{(3)}
  \cdots\eta_\varepsilon^{(\beta)}
  \int_{\gamma^{\ge \tau_0+\varepsilon}}
  \tilde{\omega}_l^{(\beta+1)}\cdots\tilde{\omega}_l^{(s)}\\
+ &\sum\limits_{1<\alpha\le\beta\le s} \;\;
  \int\limits_{\gamma_{\le \tau_0-\varepsilon}}
  (G_k\tilde{\omega}_k^{(2)})\tilde{\omega}_k^{(3)}
  \cdots\tilde{\omega}_k^{(\alpha)} \! \termx{\alpha+1}\!\! \terml{\beta+1}
  \displaybreak[0]\\
= &\int\limits_{\gamma_{\ge \tau_0+\varepsilon}}
  \tilde{\omega}_l^{(1)}\cdots\tilde{\omega}_l^{(s)} \\
+ &\sum\limits_{1=\alpha<\beta\le s} \quad \termx{1}\terml{\beta+1}\\
+ &\sum\limits_{1<\alpha\le\beta\le s}\;\;\termk{1}
  \termx{\alpha+1}\terml{\beta+1} \\
+ &\;G_l(\gamma_y(\tau_0+\varepsilon))
  \int\limits_{\gamma_{\ge \tau_0+\varepsilon}}
  \tilde{\omega}_l^{(2)}\cdots\tilde{\omega}_l^{(s)}\\ 
+ &\sum\limits_{1=\alpha<\beta\le s} P_\varepsilon(0) \int\limits_{[0;1]}
  \eta_\varepsilon^{(2)})\eta_\varepsilon^{(3)}
  \cdots\eta_\varepsilon^{(\beta)}
  \int\limits_{\gamma^{\ge \tau_0+\varepsilon}}
  \tilde{\omega}_l^{(\beta+1)}\cdots\tilde{\omega}_l^{(s)}\displaybreak[0]\\
= &\int\limits_{\gamma_{\ge \tau_0+\varepsilon}}
  \tilde{\omega}_l^{(1)}\cdots\tilde{\omega}_l^{(s)}\\
+ &\sum\limits_{1=\alpha<\beta\le s} \quad \termx{1}\terml{\beta+1}\\
+ &\sum\limits_{1<\alpha\le\beta\le s}\;\;
  \termk{1}\termx{\alpha+1}\terml{\beta+1} \\
+ &\int\limits_{[0,1]}\eta_\varepsilon^{(1)}
  \int\limits_{\gamma_{\ge \tau_0+\varepsilon}}
  \tilde{\omega}_l^{(2)}\cdots\tilde{\omega}_l^{(s)}
+ \int\limits_{\gamma^{\le \tau_0-\varepsilon}}\tilde{\omega}_k^{(1)}
  \int\limits_{\gamma_{\ge \tau_0+\varepsilon}}
  \tilde{\omega}_l^{(2)}\cdots\tilde{\omega}_l^{(s)}\\
+ &\sum\limits_{1=\alpha<\beta\le s}
  \int\limits_{\gamma^{\le \tau_0-\varepsilon}}\tilde{\omega}_k^{(1)}
  \int\limits_{[0;1]} \eta_\varepsilon^{(2)}\cdots
  \eta_\varepsilon^{(\beta)}\terml{\beta+1} \\
+ &\left[\left( G_l(\gamma_y(\tau_0+\varepsilon))-P_\varepsilon(1)\right)
  + \left( P_\varepsilon(0)- G_k(\gamma_x(\tau_0-\varepsilon))\right)\right]\\
\cdot & \int\limits_{\gamma_{\ge \tau_0+\varepsilon}}
  \tilde{\omega}_l^{(2)}\cdots\tilde{\omega}_l^{(s)}\\
+ & \sum\limits_{1=\alpha<\beta\le s}
  \left( P_\varepsilon(0) - G_k(\log \gamma_x(\tau_0-\varepsilon)\right)\\
\cdot & \int\limits_{[0;1]} \eta_\varepsilon^{(2)}\cdots
  \eta_\varepsilon^{(\beta)}\terml{\beta+1} \displaybreak[0]\\
= &\quad\eint{\gamma} \varphi_1\cdots\varphi_s\\ 
+ & \left[\left( G_l(\gamma_y(\tau_0+\varepsilon))-P_\varepsilon(1)\right)
  + \left( P_\varepsilon(0)-
  G_k(\gamma_x(\tau_0-\varepsilon))\right)\right]
  \cdot\!\!\! \int\limits_{\gamma_{\ge \tau_0+\varepsilon}}
  \tilde{\omega}_l^{(2)}\cdots\tilde{\omega}_l^{(s)}\\ 
+ &\sum\limits_{1=\alpha<\beta\le s}
  \left( P_\varepsilon(0) - G_k(\log \gamma_x(\tau_0-\varepsilon)\right)
  \cdot \int\limits_{[0;1]} \eta_\varepsilon^{(2)}\cdots
  \eta_\varepsilon^{(\beta)}\terml{\beta+1}
\end{align*}

With Fact \ref{Schwan} and Lemma \ref{Flamingo} we find that
\[
\lim\limits_{\varepsilon\to 0}\left\{
\eint{\gamma} (f\varphi_2)\varphi_3\cdots\varphi_s - \eint{\gamma}
{\varphi}_1\ldots \varphi_s \right\} =0 .
\]

\noindent
{\bf A path colliding once with the double point.}
Now let $\gamma:[a;b]\rightarrow
U_{kl}^{k-}\cup \{p_{kl}\} \cup U_{kl}^{l-}$ be a path over $\vec{v}$
that meets the double point $p_{kl}$ once with parameter value
$\tau_0\in]a;b[$ and stays in one component $D_k$.
Assume again $f(p)=f(\gamma(a))=0$. Then we compute as follows.

\newcommand{\termkk}{ \int\limits_{\gamma^{\ge \tau_0+\varepsilon}}
\!\!\!\!\tilde{\omega}_k^{(\alpha+1)}\cdots\tilde{\omega}_k^{(s)}}

\begin{align*}
& \quad \eint{\gamma} (f\varphi_2)\varphi_3\cdots\varphi_s \displaybreak[0]\\
&= \int\limits_{\gamma_{\ge \tau_0+\varepsilon}}
(G_k\tilde{\omega}_k^{(2)})\tilde{\omega}_k^{(3)}
\cdots\tilde{\omega}_k^{(s)}\\ \displaybreak[0]
&+ \sum\limits_{1\le\alpha\le s}
\;\;\int\limits_{\gamma_{\le \tau_0-\varepsilon}}
(G_k\tilde{\omega}_k^{(2)})\tilde{\omega}_k^{(3)}
\cdots\tilde{\omega}_k^{(\alpha)} \!\! \termkk \displaybreak[0]\\
&= \sum\limits_{0\le\alpha\le s}\;\;\termk{1}\termkk \\
&+ \left( G_k(\gamma_x(\tau_0+\varepsilon))-
\int\limits_{\gamma^{\le \tau_0-\varepsilon}}\tilde{\omega}_k^{(1)}
\right) \displaybreak[0] \\
&= \eint{\gamma} (f\varphi_2)\varphi_3\cdots\varphi_s
+ \left( G_k(\gamma_x(\tau_0+\varepsilon))-
G_k(\gamma(\tau_0-\varepsilon)) \right).
\end{align*}

This shows by Fact \ref{Schwan} and Lemma \ref{Flamingo}:
\[
\lim\limits_{\varepsilon\to 0}\left\{
\eint{\gamma} (f\varphi_2)\varphi_3\cdots\varphi_s - \eint{\gamma}
{\varphi}_1\ldots \varphi_s \right\} =0 .
\]

\noindent
{\bf A path locally at a double point}
Assume now that $\gamma:[a;b]\rightarrow
U_{kl}^{k-}\cup \{p_{kl}\} \cup U_{kl}^{l-}$ is a path over $\vec{v}$.
Note that we can decompose $\gamma$ into pieces each of which
meets the double point once:
$\gamma=\gamma_1\star\cdots\star\gamma_N$.
We saw already that for each of the paths $\gamma_i$ with $i=1,\ldots,N$
the assertion of the Proposition holds.

It remains to show that if for a path $\alpha:[0;1]\rightarrow
U_{kl}^{k-}\cup \{p_{kl}\} \cup U_{kl}^{l-}$ over $\vec{v}$ and for a
path $\beta:[0;1]\rightarrow U_{kl}^{k-}\cup \{p_{kl}\} \cup U_{kl}^{l-}$
over $\vec{v}$ with $\alpha(1)=\beta(0)$ the assertion of the proposition 
holds, that then the assertion is true for their composition 
$\alpha\star\beta$. This is a straightforward calculation similar 
to those above,
which finally makes again use of \ref{Schwan} and \ref{Flamingo}.
\qed

\begin{theorem} \label{Brentano}
Let $p$, $q\in U_{kl}^{k-}\cup U_{kl}^{l-}$ and let
$I_{kl}\in \bigoplus_{r=1}^s {\bigotimes}^s A^1_{kl}$
be a Chen-closed element. Then there is a ``function'' $f_{pq}\in
A^0_{kl}$ such that for any path $\gamma:[a;b]\rightarrow Z_0$
over $\vec{v}$, which lies in $U_{kl}^{k-}\cup \{p_{kl}\}\cup
U_{kl}^{l-}$ and which starts in $p$ and ends in $q$, holds:
\[
\lim_{\varepsilon\to 0} \eint\gamma I_{kl} = f_{pq}(q)- f_{pq}(p).
\]
\end{theorem}
Theorem \ref{Brentano} has the following immediate consequence.

\begin{corollary}
$\lim_{\varepsilon\to 0}\eintt{\gamma}I_{kl}$ converges 
under the assumptions of Theorem \ref{Brentano}.
\end{corollary}

\noindent
{\bf Proof of \ref{Brentano}:}
We will prove the theorem by induction on the maximal length $s$ of
tensor products in $I_{kl}$.
For $s=1$ this is the assertion of Proposition \ref{Kafka}.

\noindent
{\bf $s\ge 2$:}
Let $I_{kl}=\sum_J a_J \;\varphi_{j_1}\otimes \cdots\otimes \varphi_{j_r}$
be Chen-closed.
According to Lemma \ref{Mimikry} there is for each $J$ with $|J|=s$
at least one closed and hence exact form $\varphi_{j_\lambda}=
d f_{j_\lambda}$.

We distinguish three cases:
\begin{enumerate}
\item
If $\lambda =1$, then define
$\tilde{I}_{j_1,\ldots, j_s} := \left(f_{j_1}-f_{j_1}(p)\right) \varphi_{j_2}
\otimes\cdots\otimes\varphi_{j_s}\in{\bigotimes}^{s-1} A^1_{kl}.$
\item
If $1< \lambda<s$ define $\tilde{I}_{j_1,\ldots, j_s}
\in{\bigotimes}^{s-1} A^1_{kl}$ to be:
\[
\varphi_{j_1}\otimes\cdots\otimes
f_{j_\lambda}\varphi_{j_{\lambda+1}}\otimes\cdots\otimes\varphi_{j_s}
-\varphi_{j_1}\otimes\cdots\otimes
f_{j_\lambda}\varphi_{j_{\lambda-1}}\otimes\cdots\otimes\varphi_{j_s}.
\]
\item
Set $\tilde{I}_{j_1,\ldots, j_s}
:=-\varphi_{j_1}\otimes\cdots\otimes
\left(f_{j_s}-f_{j_s}(q)\right) \varphi_{j_{s-1}}
\in{\bigotimes}^{s-1} A^1_{kl}$,
if $\lambda=s$.
\end{enumerate}
Note that the definition of the
$\tilde{I}_{j_1,\ldots, j_s}$ only depends on $p$ and $q$. Now
\[
\tilde{I}_{kl}:=\sum_{|J|= s} a_J\; \tilde{I}_{j_1,\ldots, j_s}
+\sum_{|J|\le s-1} a_J\;
\varphi_{j_1}\otimes \cdots \otimes\varphi_{j_r}\in \bigoplus_{r=1}^s
{\bigotimes}^r A^1_{kl}
\]
is Chen-closed and by induction hypothesis there is a $f_{kl}
\in A^0_{kl}$ such that $\int_{\gamma} \tilde{I}_{kl} =
f_{pq}(q)-f_{pq}(p)$. Proposition \ref{Erstlingsausstattung} tells us:
$\lim_{\varepsilon\to 0}\{\eintt{\gamma}
{\varphi}_{j_1}\ldots\varphi_{j_s} -
\eintt{\gamma} \tilde{I}_{j_1,\ldots, j_s} \} =0$.
And this implies $\int_{\gamma} {I}_{kl} = \int_{\gamma} \tilde{I}_{kl}$.
\qed 

\subsection{Formal Power Series Connections}

Iterated integrals are not additive on paths in the naive sense
as we saw already in Proposition \ref{Badings}.
However, once we know they satisfy rules (see
\ref{Zweitlingsausstattung}) similar to those in Proposition 
\ref{Erstlingsausstattung}, it is possible to 
give Proposition \ref{Badings} an interpretation of a 
homomorphism.
The appropriate language to do this is Chen's theory
of  {\it formal power series connections} 
(cf.~\cite{Chen-Advances}, \cite{Chen-loop} and \cite{Chen-connections}).
Here we introduce an anologous language for integration with
forms in $A^\bullet_{kl}$ and $A^\bullet$. 

But first let us state the global version of Proposition 
\ref{Erstlingsausstattung}. This is necessary in order to 
apply the reduced bar construction.
 
\begin{proposition} \label{Zweitlingsausstattung}
Assume that $\gamma:[a;b]\rightarrow Z_0 $ is a path over $\vec{v}$
(based at $\vec{w}$).
Let $\varphi_1,\ldots,\varphi_s\in A^1$ be such that at least one of them
is exact.
\begin{align*}
&\text{If $\varphi_1=df_1$, then} \\
&{\textstyle\lim\limits_{\varepsilon\to 0}\left\{
\eintt{\gamma}df_1{\varphi}_2\ldots \varphi_s -\eintt{\gamma}
\big((f_1-f_1(\gamma(a))\varphi_2\big)\varphi_3\cdots\varphi_s
\right\} =0,} \\
&\text{if $\varphi_\lambda =df_\lambda$ for $2\le \lambda\le s$, then}
\\
&{\textstyle\lim\limits_{\varepsilon\to 0}\left\{
\eintt{\gamma} \varphi_1\!\cdots\! df_\lambda\cdots\varphi_s -
\eintt{\gamma} \varphi_1\!\cdots\!(f_\lambda\varphi_{\lambda+1})\!\cdots\!
\varphi_s + \eintt{\gamma} \varphi_1
\!\cdots\!(f_\lambda\varphi_{\lambda-1})\!\cdots\!
\varphi_s\!\right\}=0,} \\
&\text{if $\varphi_s=df_s$, then} \\
&{\textstyle\lim\limits_{\varepsilon\to 0}\left\{
\eintt{\gamma} \varphi_1\cdots\varphi_{s-1} df_s -
\eintt{\gamma} \varphi_1\cdots\varphi_{s-2}\big(f_s
-f_s(\gamma(b))\varphi_{s-1}\big)\right\}=0.}
\end{align*}
\end{proposition}

\noindent
{\bf Proof:}
Assume now that $\gamma:[a;b]\rightarrow
Z_0$ is a path over $\vec{v}$ (based at $\vec{w}$).
Note that according to Proposition \ref{Erstlingsausstattung} 
we can decompose $\gamma$ into piece such that for each of them 
the assertion of Proposition \ref{Zweitlingsausstattung} holds:
$\gamma=\gamma_1\star\cdots\star\gamma_N$.
One has to show that if for a path $\alpha:[0;1]\rightarrow
Z_0$ over $\vec{v}$ (based at $\vec{w}$) and for a
path $\beta:[0;1]\rightarrow Z_0$ 
over $\vec{v}$ (based at $\vec{w}$) 
with $\alpha(1)=\beta(0)$ the assertion of 
Proposition \ref{Zweitlingsausstattung}
holds, that then the assertion is true for their composition
$\alpha\star\beta$. This is a straightforward calculation similar
to those in the proof of Proposition \ref{Erstlingsausstattung}.
\qed

\subsubsection{Formal Power Series Connections}
Let $\varphi_{i}$, $\varphi_{ij}$, $\varphi_{ijk}$, ...;
$i,\,j,\,k=1,2,\ldots,s$ be elements in $A^1$. Let 
$X_1,\ldots,X_s$ be noncommutative indeterminates.

The formal power series
\[
\omega= \sum \varphi_{i}X_i+\sum \varphi_{ij}X_i X_j
+\sum \varphi_{ijk} X_i X_j X_k + \cdots
\] 
will be called a {\it formal power series connection for $A^\bullet$}.
We refer to the number $s$ as to the {\it length} of $\omega$.
The {\it curvature} $\kappa$ of $\omega$ is defined by:
\[
\kappa=d\omega + \omega \wedge \omega.
\]
Associated to the given power series connection $\omega$ and a path 
$\gamma:[a;b]\rightarrow Z_0$ over $\vec{v}$ (based at $\vec{w}$) as well as
$\varepsilon>0$ small enough there is a formal power series with coefficients
in $\cpxm$, {\it the $\varepsilon$-transport of $\omega$ applied to 
$\gamma$}
\[
\epsT(\gamma)= 1 +\eint{\gamma} \omega+\eint{\gamma} \omega\omega
+\eint{\gamma} \omega\omega\omega+\cdots,
\] 
whose coefficients are given by 
\begin{gather*}
\epsT_i(\gamma) = \eint{\gamma} \varphi_i,\qquad
\epsT_{ij}(\gamma) = \eint{\gamma} \varphi_{i}\varphi_{j}+\varphi_{ij}   \\
\epsT_{ijk}(\gamma) = \eint{\gamma}\varphi_{i}\varphi_{j}\varphi_{k}
+\varphi_{ij}\varphi_{k}+\varphi_{i}\varphi_{jk}+ \varphi_{ijk}\\
\intertext{and in general}
\epsT_{i_1,\ldots,i_r}(\gamma) = \sum_{0<\alpha_1,\cdots<\alpha_{r-1}<r}
\;\eint{\gamma}\varphi_{i_1\cdots i_{\alpha_1}}\cdots
\eint{\gamma}\varphi_{i_{\alpha_{r-1}+1}\cdots i_{r}}.
\end{gather*}

The reason for us to make these definitions is the following Proposition.
\begin{proposition} \label{Babytrainer}
Let $\bf p$, $\bf q$, ${\bf r}\in Z_0\setminus\bigcup_{[k<l]}\{p_{kl}\}$ 
and let $\alpha,\;\beta:[0;1]\rightarrow Z_0$ be paths over $\vec{v}$
(based at $\vec{w}$), where $\alpha$ goes from $\bf p$ to $\bf q$ and
$\beta$ from $\bf q$ to $\bf r$. Let $\omega$ be a power series connection.
Then for $\varepsilon>0$ small enough:
\[
\epsT(\alpha\star\beta)=\epsT(\alpha)\cdot\epsT(\beta).
\]
\end{proposition}

\noindent
{\bf Proof:}
The proof follows from the following rule
\[
\eint{\alpha\star\beta}\underbrace{\omega\cdots\omega}_n=
\sum_{k=0}^n\quad\eint{\alpha}\underbrace{\omega\cdots\omega}_k\quad
            \eint{\beta}\underbrace{\omega\cdots\omega}_{n-k},
\]
which in turn can be seen to be a consequence of 
Proposition \ref{Badings}. \qed

Not difficult to prove is the following proposition when one uses
Lemma \ref{Ligeti}.
\begin{proposition}
Let $\omega= \sum \varphi_{i}X_i+\sum \varphi_{ij}X_i X_j
+\sum \varphi_{ijk} X_i X_j X_k + \cdots$ 
be a formal power series connection for $A^\bullet$. 
Then $\kappa=d\omega + \omega \wedge \omega =0$
if and only if the coefficients of the formal power series 
$1 + \omega+ \omega\otimes\omega+ \omega\otimes\omega\otimes\omega+\cdots$
are closed under $(d_I+d_C)$. \qed
\end{proposition}

The next proposition allows us to use formal power series connections
in order to prove things about Chen-closed elements in
$\bigoplus_{r=0}^s {\bigotimes}^r A^1$.

\begin{proposition} \label{dwarsligger}
Let $I\in\bigoplus_{r=0}^s {\bigotimes}^r A^1$  be Chen-closed. Then
there exists a finite number of
formal power series connections $\omega_1,\ldots,\omega_N$ of length $s$
with curvature 0 such that up to an element in $\R_1^s(A^\bullet,a)$
$I$ is the coefficient of $X_1\cdots X_s$ in: 
\[
\sum_{i=1}^N \left( 1+\omega_i +\omega_i\otimes\omega_i +
\omega_i\otimes\omega_i\otimes\omega_i+ \cdots \right).
\]
\end{proposition}

\noindent
{\bf Proof:}
Induction on $s$: $s=1$ is trivial. So $s>1$.
By Lemma \ref{Claus Hertling} and Lemma \ref{Ligeti} we find that
there are fomal power series connections (fpsc's)
$\omega_1^{(s)},\ldots,\omega_K^{(s)}$ of length $s$ such that the
difference $I^{(s)}$ between $I$ and the coefficient of $X_1\cdots X_s$ in
$\sum_{i=1}^K\left( 1+\omega_i^{(s)} +\omega_i^{(s)}\otimes\omega_i^{(s)} 
+ \cdots \right)$ is of length $s-1$ and closed under $(d_I+d_C)$. 
By induction hypothesis there exist fpsc's $\tilde{\omega}_1^{(s-1)},\ldots,
\tilde{\omega}_L^{(s-1)}$ of length $s-1$ such that $I^{(s-1)}$ is the 
coefficient of $X_1\cdots X_{s-1}$
in $\sum_{i=1}^L\left( 1+\tilde{\omega}_i^{(s-1)}
+\tilde{\omega}_i^{(s-1)}\otimes\tilde{\omega}_i^{(s-1)} + \cdots \right)$.
Define the
fpsc's ${\omega}_i^{(s-1)}$ to be $\tilde{\omega}_i^{(s-1)}$ except that the
indeterminate $X_{s-1}$ is replaced by $X_{s-1}X_s$. The
${\omega}_i^{(s-1)}$ still have zero curvature! and $I$ is the
coefficient of $X_1\cdots X_s$ in
\[
\sum_{i=1}^K \left( 1+\omega_i^{(s)} +\omega_i^{(s)}\otimes\omega_i^{(s)} +
\cdots \right)
+
\sum_{i=1}^L \left( 1+\omega_i^{(s-1)} +
\omega_i^{(s-1)}\otimes\omega_i^{(s-1)} +
\cdots \right).
\]
\qed

Proposition \ref{Babytrainer} enables us to prove that the 
iterated integrals along paths over $\vec{v}$ converge 
not just locally but also globally.
\begin{theorem} \label{Milcheinschuss}
Let $\omega$ be a formal power series connection for $A^\bullet$
with curvature $\kappa=0$, let
$\gamma:[0;1]\rightarrow Z_0$ be a path over $\vec{v}$
(based at $\vec{w}$) and let $\epsT(\gamma)$ be the associated 
$\varepsilon$-transport of $\omega$ applied to $\gamma$. Then 
\[
\T(\gamma):=\lim_{\varepsilon\to 0} \epsT(\gamma)
\]
converges and is called {\it the transport of $\omega$ applied to
$\gamma$}. Moreover if
$\alpha,\;\beta:[0;1]\rightarrow Z_0$ are paths over $\vec{v}$
(based at $\vec{w}$), which can be composed to a connected path over $\vec{v}$,
then 
\[
\T(\alpha\star\beta)=\T(\alpha)\cdot \T(\beta).
\]
\end{theorem}
 
\noindent
{\bf Proof:} 
Let $\gamma:[a;b]\rightarrow Z_0$ be a path over $\vec{v}$ (based at $\vec{w}$)
that meets the set of double points $\{p_{kl}\}$ with parameter values
$a=\tau_0<\tau_1<\cdots<\tau_m<\tau_{m+1}=b$. At the double point 
$p_{0}\in D_0$, we have the embedding $x=p:D_0\rightarrow \cpxm$.
Like before, let $\sigma$ be the straight path $t\mapsto (1-t)$ from 
$1$ to $0 \in D_0$. Decompose $\sigma\star\gamma\star\sigma^{-1}=
\gamma_0\star\cdots\star\gamma_{m+1}$ such that for each 
$\gamma_\kappa:[a_\kappa;b_\kappa]\rightarrow Z_0$, the interval
$[a_\kappa;b_\kappa]$ 
contains only the parameter value $\tau_\kappa$ in its interior. 
There is an $\varepsilon_0 >0$ with which we can decompose 
each of the paths $\gamma_\kappa$ into $\gamma_\kappa=\alpha_\kappa\star
\delta_\kappa\star\beta_\kappa$ such that $\alpha_\kappa$ resp.~$\beta_\kappa$
lie entirely in one component and each $\delta_\kappa$ is a map
\[
\delta_\kappa: [\tau_\kappa-\varepsilon_0;\tau_\kappa+\varepsilon_0]
\rightarrow U_{kl}^{k-}\cup \{p_{kl}\}\cup U_{kl}^{l-}
\]
for some $[k<l]$. 

By Theorem \ref{Brentano} we know that if $\vartheta$ is one of the paths $\alpha_\kappa$, $\delta_\kappa$, $\beta_\kappa$ then $\epsT(\vartheta)$
converges. The proof can now be accomplished by applying 
Proposition \ref{Babytrainer}.
\qed

Hence we may conclude the following corollary.
\begin{corollary} \label{Babymassage}
Let $I\in\bigoplus_{r=0}^s {\bigotimes}^r A^1$ 
be a Chen-closed element and let $\gamma:[0;1]\rightarrow Z_0$ 
be a path over $\vec{v}$ (based at $\vec{w}$). 
Then $\lim_{\varepsilon\to 0} \eintt{\gamma} I$ converges.
If $\gamma:[0;1]\rightarrow Z_0$ is a path over $\vec{v}$
based at $\vec{w}$, then $\int_{\gamma}(\R_1^s(A^\bullet,a))=0$.
\end{corollary}

\noindent
{\bf Proof:}
The assertions in the corollary are consequences of 
Proposition \ref{dwarsligger}, Proposition \ref{Zweitlingsausstattung} and Corollary \ref{Milcheinschuss}.
\qed 

\begin{theorem} \label{Vormilch}
Let $I\in\bigoplus_{r=0}^s {\bigotimes}^r A^1$ be Chen-closed and
let $\gamma:[0;1]\rightarrow Z_0$ be a path over $\vec{v}$ (based at
$\vec{w}$). Then $\int_\gamma I$ does not depend on the choice of
coordinates $(x,y)$ around the double points.
\end{theorem}

\begin{remark} \rm
Theorem \ref{Brentano} may not be applied directly to prove 
Theorem \ref{Vormilch} since a priori the definition of the dga 
$A_{kl}^\bullet$ depends on coordinates.
\end{remark}

\noindent
{\bf Proof of \ref{Vormilch}:}
By Lemma \ref{dwarsligger}, Proposition \ref{Zweitlingsausstattung}
and Corollary \ref{Milcheinschuss} it is enough to prove the theorem for 
the following local situation. Assume without loss of generality
that the path $\gamma:[0;1]\rightarrow Z_0$ over $\vec{v}$
lies entirely in a $U_{kl}^{k-}\cup \{p_{kl}\} \cup U_{kl}^{l-}$ for $[k<l]$
and meets the double point $p_{kl}$ once with parameter value $\tau_0$.
Like earlier we distinguish two cases.

\noindent
{\bf A path traversing the double point.}
We have
\[
\int_\gamma I =\lim_{\varepsilon\to 0}\left\{
\sum_J \eint{\gamma} \varphi_{j_1}\cdots\varphi_{j_r}\right\}
\]
Recall that for each $J$ there are polynomials $Q_{\alpha\beta}^J\in \cpxm[u]$
such that $\eintt{\gamma} \varphi_{j_1}\cdots\varphi_{j_r}$
\[
 = \sum\limits_{0\le\alpha\le\beta\le r}\termk{1}
Q_{\alpha\beta}^J\left(\log\gamma_x(\tau_0-\varepsilon)\cdot
\gamma_y(\tau_0+\varepsilon)\right)\terml{\beta+1}.
\]

Furthermore remark:
\[
\left(\log\left(\gamma_x(\tau_0-\varepsilon)\cdot
\gamma_y(\tau_0+\varepsilon)\right)\right)^m
=\left(\log\left(\frac{\gamma_x(\tau_0-\varepsilon)\cdot
\gamma_y(\tau_0+\varepsilon)}{\varepsilon^2}\right)
+\log\varepsilon^2\right)^m
\]
\begin{align*}
=\left(\log\varepsilon^2\right)^m 
&+\bigg\{\log\left({\small \frac{\gamma_x(\tau_0-\varepsilon)\cdot
\gamma_y(\tau_0+\varepsilon)}{\varepsilon^2}}\right) \\
&\cdot \left({\small \sum\limits_{k=1}^m \binom{m}{k}
\log^{k-1}\left(\frac{\gamma_x(\tau_0-\varepsilon)\cdot
\gamma_y(\tau_0+\varepsilon)}{\varepsilon^2}\right)
\left(\log{\varepsilon^2}\right)^{m-k}}\right)\bigg\}.
\end{align*}
The function $\log\left(\frac{\gamma_x(\tau_0-\varepsilon)\cdot
\gamma_y(\tau_0+\varepsilon)}{\varepsilon^2}\right)$ is a
differentiable function in $\varepsilon$ on an intervall
$(-\varepsilon_0;\,\varepsilon_0)$ for some $\varepsilon_0$,
which vanishes in $0$. Hence the vanishing of this function is stronger
than any logarithmic growth (Use Fact \ref{Schwan} on page \pageref{Schwan}).
This indicates that the following limit is zero:
\[
\lim_{\varepsilon\to 0}\left\{
\eint{\gamma} \varphi_{j_1}\cdots\varphi_{j_r}
-\sum\limits_{0\le\alpha\le\beta\le r}\termk{1}
Q_{\alpha\beta}^J\left(\log\varepsilon^2\right)\terml{\beta+1}\right\}.
\]
We find that $\int_\gamma I$ does not depend on coordinates.

\noindent
{\bf A path colliding with a double point}\\
Let $p:=\gamma(0)$ and $q:=\gamma(1)$. We know by \ref{Brentano}
that there exists a closed and hence exact
$\varphi_{pq}=d f_{pq}\in A^1_{kl}$ such that
\[
\int_\gamma I = \int_\gamma \varphi_{pq} = f_{pq}(q)- f_{pq}(p).
\]
Let $\breve{\gamma}:[a;b]\rightarrow U_{kl}^{k-}$ be a path
(over $\vec{v}$), which is nearby homotopic relative to the
end points with $\gamma$, but does not touch the double point.
According to Proposition \ref{Kafka} we observe:
\[
\int_\gamma \varphi_{pq} = f_{pq}(q)-f_{pq}(p) =
\int_{\breve{\gamma}} \varphi_{pq}
\]
The integral on the right is independent of the choice of coordinates.
\qed

\subsection{Nearby Homotopy Functionals}

In this subsection we want to prove that iterated
integrals of Chen-closed elements in $\bigoplus_{r=1}^s{\bigotimes}^r A^1$
are invariant under a homotopy over $\vec{v}$ based at $\vec{w}$. 
Hence, they are {\it nearby homotopy functionals}. There is no converse of 
this assertion since in general the iterated integrals only 
converge for Chen-closed elements in $\bigoplus_{r=1}^s{\bigotimes}^r A^1$.

\afterpage{\clearpage
\begin{figure}[H]
\qquad\qquad\qquad\qquad\qquad\epsfig{file=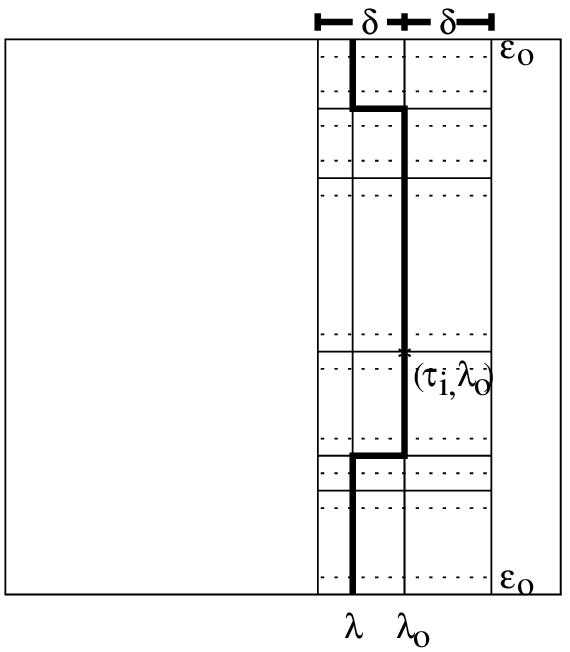,clip=}
\caption{Illustration of the domain of $H$, the path $\eta$
and the construction of the sequence of paths
$(\alpha_i)_{i=0,\ldots,M}$.}
\end{figure}}

\begin{theorem} \label{Blaek Foes}
Let $I\in\bigoplus_{r=0}^s {\bigotimes}^r A^1$ be Chen-closed. 
Then for two nearby homotopic paths $\gamma_0,\gamma_1:
[0;1]\rightarrow Z_0$ over $\vec{v}$ based at $\vec{w}$ holds:
\[
\int_{\gamma_0} I = \int_{\gamma_1} I.
\]
\end{theorem}

\noindent
{\bf Proof:}
By Lemma \ref{dwarsligger} and Proposition \ref{Zweitlingsausstattung}
it is enough to prove the corresponding theorem for a formal power 
series connection $\omega$ with curvature $\kappa=0$. Let $\T$ be the
associated transport. We want to prove: $\T({\gamma_0}) = \T({\gamma_1})$

Let $H:[0;1]\times[0;1]\rightarrow Z_0$ be the nearby homotopy.
That is, $H(\cdot,0)=\gamma_0$ and $H(\cdot,\gamma_1)=\gamma_1$
and for any $\lambda\in [0;1]$ is $H(\cdot,\lambda)$ is a path
over $\vec{v}$ based at $\vec{w}$. Consider the function:
\[
\begin{array}{cccc}
C: & [0;1] & \longrightarrow & \cpxm \\
   &\lambda& \longmapsto     & \T({H(\cdot,\lambda)}).
\end{array}
\]
We will show that $C$ is constant by showing that it is locally constant.

Let $\lambda_0\in [0;1]$. We show that $C$ is locally constant at
$\lambda_0$. Write $\eta:=H(\cdot,\lambda_0)$ and denote by
$0=\tau_0<\tau_1<\cdots<\tau_N<\tau_{n+1}=1$ the parameter values, where
$\eta$ meets the set of double points: $\eta$ defines a map
\[
\mu:\{0,\ldots,N+1\}\longrightarrow \left\{(k,l)\big|\; [k<l]\right\},
\]
by $\eta(\tau_i)=p_{\mu(i)}$ for $i=0,\ldots,N+1$.

For any $\varepsilon > 0$ we define $L_0^\varepsilon:=[0;\varepsilon]$,
$L_i^\varepsilon:=[\tau_i-\varepsilon;\tau_i+\varepsilon]$ for $i=1,\ldots,N$
and $L_{N+1}^\varepsilon:=[1-\varepsilon;1]$.
Choose a positive $\varepsilon_1< \min\limits_i \frac{|\tau_{i+1}-\tau_i|}{2}$
(then all the $L_i^{\varepsilon_1}$ are mutually disjoint) such that:
$\eta\left(L_i^{\varepsilon_1}\right) \subset U_{\mu(i)}^{k-}
\cup \{ p_{\mu(i)}\}\cup U_{\mu(i)}^{l-}$.

By the continuity of $H$ and the fact that all $H(\cdot,\lambda)$
are paths over $\vec{v}$ we can find a $\delta>0$ and an
$\varepsilon_0>0$ (with $\varepsilon_0<\varepsilon_1$) such that
for all $i=0,\ldots,N+1$ holds:
\[
H\left(L_i^{\varepsilon_0}\times[\lambda_0-\delta;\lambda_0+\delta]\right)
\subset U_{\mu(i)}^{k-}
\cup \{ p_{\mu(i)}\}\cup U_{\mu(i)}^{l-}.
\]
and the inverse image of the double points under
$H_{\big| [0;1]\times[\lambda_0-\delta,\lambda_0+\delta]}$
lies entirely in the interior of the set
$\bigcup_i L_i^{\varepsilon_0}\times[\lambda_0-\delta,\lambda_0+\delta]$.
(If $\lambda_0 =0$ resp.~$\lambda_0=1$ then replace
$[\lambda_0-\delta,\lambda_0+\delta]$ by $[0,\delta]$
resp.~$[1-\delta,1]$.)

Fix a $\lambda\in [\lambda_0-\delta,\lambda_0+\delta]$.
It is a consequence of Theorem \ref{Brentano}
that $\T$ takes the same value
on any two paths $H\circ\theta_i$ ($i=1,\,2$) over $\vec{v}$,
where the $\theta_i$ are both paths with the same end points
within $L_i^{\varepsilon_0}\times[\lambda_0-\delta;\lambda_0+\delta]$,
whose end points are two of the four points:
$\{ (\tau_i\pm \varepsilon_0,\lambda_0),\;
(\tau_i\pm \varepsilon_0,\lambda)\}$.

Moreover $\int I$ takes the same value
on any two paths $H\circ\theta_i$ ($i=1,\,2$) over $\vec{v}$,
where the $\theta_i$ are both paths with the same end points
within $[\tau_i+\varepsilon_0,\tau_{i+1}-\varepsilon_0]
\times[\lambda_0-\delta;\lambda_0+\delta]$,
whose end points are two of the four points:
$\{ (\tau_i+ \varepsilon_0,\lambda_0),\;
(\tau_{i+1}- \varepsilon_0,\lambda_0),\;
(\tau_i+ \varepsilon_0,\lambda),\;
(\tau_{i+1}-\varepsilon_0,\lambda)\}$,
likewise by Theorem \ref{Brentano}.

Now one easily constructs a finite sequence of paths
$\alpha_0$, $\alpha_1,\ldots,\alpha_M$ in $[0;1]
\times[\lambda_0-\delta;\lambda_0+\delta]$ with the following
properties:

\begin{itemize}
\item
$H\circ\alpha_0=H(\cdot,\lambda_0)=\eta$ and $H\circ\alpha_M=H(\cdot,\lambda)$,
\item
each $H\circ\alpha_i$ is a path over $\vec{v}$ based at $\vec{w}$,
\item
two successive paths $\alpha_i$, $\alpha_{i+1}$ 
coincide everywhere except either on an interval
$L_i^{\varepsilon_0}\times[\lambda_0-\delta;\lambda_0+\delta]$ or on
an interval $[\tau_i+\varepsilon_0,\tau_{i+1}-\varepsilon_0]
\times[\lambda_0-\delta;\lambda_0+\delta]$.
\end{itemize}
In other words we may write
$\alpha_i=\zeta\star\theta\star\vartheta$ and 
$\alpha_{i+1}=\zeta\star\theta'\star\vartheta$. The operator
$\T$ takes the same value on $\theta$ and $\theta'$.
Apply Corollary \ref{Milcheinschuss} to conclude
$\T({H\circ\alpha_i})=\T({H\circ\alpha_{i+1}})$
\qed

\subsection{Nilpotent Orbits} \label{Hokuspokus}

Here we will apply the formal power series connections to prove 
Theorem \ref{Gr"unspan}. Since the operators $M$ and $N$ commute, 
we may prove the theorem by proving the following 
two equations separately:
\begin{enumerate}
\item
$\left({J}_{\lambda\vec{v},\vec{w}}\big/
{J}_{\lambda\vec{v},\vec{w}}^{s+1}\right)^*_\itgm=\lambda^{-N}
\left({J}_{\vec{v},\vec{w}}\big/{J}_{\vec{v},\vec{w}}^{s+1}\right)^*_\itgm
\subset H^0\bar{B}_s(A^\bullet,a)$,
\item 
$\left({J}_{\vec{v},\mu\vec{w}}\big/
{J}_{\vec{v},\mu\vec{w}}^{s+1}\right)^*_\itgm=\mu^{M}
\left({J}_{\vec{v},\vec{w}}\big/{J}_{\vec{v},\vec{w}}^{s+1}\right)^*_\itgm
\subset H^0\bar{B}_s(A^\bullet,a)$.
\end{enumerate}

\noindent
{\bf Proof of Theorem \ref{Gr"unspan}:} 
We will prove (i) and leave (ii) for the reader. 

Let $\lambda:[0;1]\rightarrow \cpxm\setminus \reellm^{\le 0}$ be a path
with $\lambda(0)=1$. Let $\lambda_0$ be a positive real number such that 
$\lambda_0<|\lambda(\varsigma)|$ for all $\varsigma\in[0;1]$. 

Like always, we fix coordinates $(x,y):W_{kl}=U_{kl}^k\times U_{kl}^l
\rightarrow \cpxm^2$ with respect to $\vec{v}$, that means 
$\langle \frac{\partial}{\partial x},
\frac{\partial}{\partial y}\rangle=\vec{v}$.
For an appropriate $\rho>0$ let 
$\tilde{y}_\varsigma=\tilde{y}:U_{kl}^l\rightarrow \cpxm$ be another coordinate
that depends continuously on $\lambda(\varsigma)$ resp.~$\varsigma$ 
and that satisfies:
\begin{xalignat*}{2}
\tilde{y}_\varsigma &= y &&\text{ for } |y|\ge \frac{\rho}{\lambda_0}, \\
\tilde{y}_\varsigma &= \frac{1}{\lambda(\varsigma)}y &&\text{ for } |y|\le \rho
\text{ and }\\
\tilde{y}_0 &= y &&\text{ for all } y. 
\end{xalignat*}
For example, take a smoothing of the map that is given in polar 
coordinates by:
\[
\begin{pmatrix} \tilde{r} \\ \tilde{\varphi} \end{pmatrix} =
\begin{pmatrix} (r-\rho)\frac{|\lambda|-\lambda_0}{\lambda-|\lambda|\lambda_0}
+\frac{\rho}{|\lambda|} \\ 
\varphi + \arg (\lambda) \left(\frac{\lambda_0 r-\rho}{\rho-\rho\lambda_0}
\right) 
\end{pmatrix} \quad \text{for } \rho\le r \le \frac{\rho}{\lambda_0.}
\]
These coordinates $(x,\tilde{y}_\varsigma)$ are made such that holds: 
$\langle \frac{\partial}{\partial x},
\frac{\partial}{\partial \tilde{y}_\varsigma}\rangle
=\lambda(\varsigma)\vec{v}$. 

The coordinate changes $\tilde{y}_\varsigma^{-1}\circ y$ at any $U_{kl}^l$
can be put together to an isotopy
\[
\phi_\varsigma: Z_0\rightarrow Z_0,
\]
which is defined to be $\tilde{y}_\varsigma^{-1}\circ y$ on any $U_{kl}^l$
and the identity elsewhere. Note that $\phi_0=id$ and 
$\phi_\varsigma(y)=\lambda(\varsigma) y$ for $y\in U_{kl}^l$ with $|y|\le \rho$. 
Given a path $\gamma:[a;b]\rightarrow Z_0$ over $\vec{v}$ based at $\vec{w}$,
we can define a homotopy $H$ by $H(t,\varsigma):=\phi_\varsigma\circ \gamma(t)$.
This $H(\cdot,\varsigma)$ is a path over $\lambda(\varsigma)\vec{v}$ based at 
$\vec{w}$.

Now take a Chen-closed $I=\sum_{J} a_J\; \varphi_{j_1} \otimes\cdots\otimes
\varphi_{j_r}\in \bigoplus_{r=1}^s\bigotimes^r A^1$. We obtain by the definition
of iterated integrals along paths over $\lambda(\varsigma)\vec{v}$ based at $\vec{w}$:
\begin{equation} \label{Howdy}
\int_{H(\cdot,\varsigma)} I = \int_\gamma \phi_\varsigma^* I.
\end{equation}
Here the dga-automorphism $\phi_\varsigma^*:A^\bullet\rightarrow A^\bullet$
is defined by 
\begin{alignat*}{2}
({\phi_\varsigma}_{\big| D_i})^* &: E^\bullet(D_i\log P_i)
&\rightarrow &E^\bullet(D_i\log P_i) \text{     and } \\
id &: E^\bullet(\Delta^1)\otimes \Lambda^\bullet &\rightarrow 
&E^\bullet(\Delta^1)\otimes \Lambda^\bullet .
\end{alignat*}
$\phi_\varsigma^*$ induces an automorphism 
$\bigoplus_{r=1}^s\bigotimes^r \phi_\varsigma^*$ of 
$\bigoplus_{r=1}^s\bigotimes^r A^\bullet$, which we also denote by 
$\phi_\varsigma^*$. This is the $\phi_\varsigma^*$ in (\ref{Howdy}). 

We are done with the proof of the theorem if we can show for all paths
$\gamma\in\pi_1(Z_{\vec{v}},\vec{w})$:
\[
\int_\gamma \phi_\varsigma^* I = \int_\gamma \lambda(\varsigma)^N I.
\]
Also this statement will be proved by first looking at the local and then 
at the global situation.

\noindent
{\bf Local Situation:} 
Assume here that $\vartheta:[a;b]\rightarrow U_{kl}^{k-}\cup \{p_{kl}\}\cup
U_{kl}^{l-}$ is a path over $\vec{v}$ that meets the double point $p_{kl}$
once with parameter value $\tau_0$, where it changes from $D_k$ to $D_l$.
(Since $\phi_\varsigma^* I$ and $\lambda(\varsigma)^N I$ are nearby homotopy 
functionals, we do not have to worry about paths that {\it collide 
with a double point}.) Suppose the endpoint $\vartheta(b)$ of $\vartheta$
be such that $|y(\vartheta(b))|>\frac{\rho}{\lambda_0}$. 

Notice that the two coordinate systems $(x,y)$ and $(x,\tilde{y}_\varsigma)$
around $p_{kl}$ give two different {\it localizations in a sector 
at the double point $p_{kl}$}
\[
A^\bullet \rightarrow A^\bullet_{kl} \qquad \text{  and   } \qquad
A^\bullet \rightarrow \tilde{A}^\bullet_{kl}
\]
respectively. Denote the images of $I$ under these natural maps by $I_{kl}$
resp.~$\tilde{I}_{kl}$. Observe moreover that $\phi_\varsigma$ 
induces 
\[
\phi_\varsigma:U_{kl}^{k-}=
\{y\not\in\reellm^{\le 0}\}\xrightarrow{\approx} \tilde{U}_{kl}^{k-}=
\{\tilde{y}_\varsigma\not\in\reellm^{\le 0}\}. 
\]

Then we can extend the dga-automorphisms $\phi_\varsigma^*$ and 
$\lambda(\varsigma)^N$
of $A^\bullet$ to dga-isomorphisms with the same names,
$\phi_\varsigma^*$ and $\lambda(\varsigma)^N$, from $\tilde{A}^\bullet_{kl}$ 
to $A^\bullet_{kl}$ like indicated in the following two commutative
diagrams.
\[
\begin{CD}
A^\bullet @>>> \tilde{A}^\bullet_{kl} & \qquad & \log \tilde{y} &  
& \qquad & w \\
@V\phi_\varsigma^*VV @VV\phi_\varsigma^*V @VVV \qquad & @VVV  \\
A^\bullet @>>> {A}^\bullet_{kl} &\qquad & \log \tilde{y}\circ\phi_\varsigma(y) 
& & \qquad &  w,
\end{CD}
\]
I.~e.~$\phi_\varsigma$ on $\tilde{A}^\bullet_{kl}$ is the pull-back with 
$\phi_\varsigma^*$. Extend $\lambda(\varsigma)^N$ by
\[
\begin{CD}
A^\bullet @>>> \tilde{A}^\bullet_{kl} & \qquad & \log \tilde{y} &  
& \qquad & w \\
@V{\lambda(\varsigma)^N}VV @VV{\lambda(\varsigma)^N}V @VVV \qquad & @VVV  \\
A^\bullet @>>> {A}^\bullet_{kl} &\qquad & \log \tilde{y}(y) 
& & \qquad &  (w-\log\lambda).
\end{CD}
\]

We know by Theorem \ref{Brentano} that for $\tilde{I}_{kl}$ and the 
points $\vartheta(b)$, $\vartheta(a)$ there is a function 
$\tilde{f}_{kl}\in\tilde{A}^0_{kl}$ such that holds:
\[
\int_\vartheta {I}= \int_\vartheta \tilde{I}_{kl} =
\tilde{f}_{kl}\big(\tilde{y}\circ\vartheta(b)\big)-
\tilde{f}_{kl}(\tilde{y}\circ\vartheta(a)\big).
\]
Since both, $\phi_\varsigma^*$ and $\lambda(\varsigma)^N$, are isomorphisms of 
dga's, which preserve the values of functions on $\vartheta(b)$ and
$\vartheta(a)$, the pull-back $\phi_\varsigma^*\tilde{f}_{kl}$ is such a
function for $\phi_\varsigma^*\tilde{I}_{kl}$ and 
$\lambda(\varsigma)^N\tilde{f}_{kl}$ is such a
function for $\lambda(\varsigma)^N\tilde{I}_{kl}$. That is to say:
\begin{align*}
\int_\vartheta {\phi_\varsigma^* I}&= 
\int_\vartheta \phi_\varsigma^*\tilde{I}_{kl} =
\left(\phi_\varsigma^*\tilde{f}_{kl}\right)(\tilde{y}\circ\vartheta(b)\big)-
\left(\phi_\varsigma^*\tilde{f}_{kl}\right)(\tilde{y}\circ\vartheta(a)\big) \\
\intertext{and }
\int_\vartheta \lambda(\varsigma)^N {I} &= 
\int_\vartheta \lambda(\varsigma)^N\tilde{I}_{kl} =
\left(\lambda(\varsigma)^N\tilde{f}_{kl}\right)
(\tilde{y}\circ\vartheta(b)\big)-
\left(\lambda(\varsigma)^N\tilde{f}_{kl}\right)
(\tilde{y}\circ\vartheta(a)\big).
\end{align*}

If $\tilde{f}_{kl}$ is given by 
\begin{align*}
\tilde{f}_{kl} &= \sum_{\nu=0}^m g_{k,\nu}(x) \log^\nu x
+ \sum_{\nu=0}^m g_{l,\nu}(\tilde{y}) \log^\nu \tilde{y} + P(\xi,v,w) \\
\intertext{then}
\phi_\varsigma^*\tilde{f}_{kl} &= \sum_{\nu=0}^m g_{k,\nu}(x) \log^\nu x
+ \sum_{\nu=0}^m g_{l,\nu}({y}) \log^\nu {y} + P(\xi,v,w) \\
\intertext{and}
\lambda(\varsigma)^N\tilde{f}_{kl} &= \sum_{\nu=0}^m g_{k,\nu}(x) \log^\nu x
+ \sum_{\nu=0}^m g_{l,\nu}(\tilde{y}( y)) \log^\nu 
\big(\tilde{y}( y)\big) + P(\xi,v,w-\log \lambda).
\end{align*}

Because of $\tilde{y}\big( y\circ\vartheta(b)\big)=y\circ\vartheta(b)$ as
$|y\circ\vartheta(b)|>\frac{\rho}{\lambda_0}$, we finally find:
\[
\int_\vartheta {\phi_\varsigma^* I}=\int_\vartheta \lambda(\varsigma)^N {I}.
\] 

\noindent
{\bf Global Situation:} 
Let $\gamma:[a;b]\rightarrow Z_0$ be a path over $\vec{v}$ based at $\vec{w}$.
At the double point $p_{0}\in D_0$, we have the embedding 
$x=p:D_0\rightarrow \cpxm$.
Like before, let $\sigma$ be the straight path $t\mapsto (1-t)$ from 
$1$ to $0 \in D_0$. Decompose like in the proof of Theorem \ref{Milcheinschuss}
the path $\sigma\star\gamma\star\sigma^{-1}=
\gamma_0\star\cdots\star\gamma_{m+1}$ such that for a Chen-closed $I$ 
each of the paths $\gamma_\nu$ satisfies: $\int_{\gamma_\nu} 
{\phi_\varsigma^* I}=\int_{\gamma_\nu} \lambda(\varsigma)^N {I}$.

Let $\omega$ be a formal power series connection with curvature zero.
Then we can apply ${\phi_\varsigma^*}$ and $\lambda(\varsigma)^N$
to $\omega$ by applying them to each coefficient. We know that the
transport associated to ${\phi_\varsigma^*}\omega$ and the
transport associated to $\lambda(\varsigma)^N\omega$ take the same value on 
all paths $\gamma_\nu$. Apply Theorem \ref{Milcheinschuss} 
to prove the theorem.
\qed


\section{Plane Curve Singularities} \label{Kranenburg}

In this last section, we finally want to show, how the setting of 
\ref{Setzauf} arises in the case of plane curve 
singularities. Throughout the whole section let 
\[
f=f_0\cdots f_{r-1}:\left(\cpxm^2,0\right)\rightarrow\left(\cpxm,0\right)
\]
be a plane curve singularity with $r$ branches.

First we will construct 
a map $h:(Z,D^+)\rightarrow (\Delta,0)$ with the properties
described earlier in \ref{Setzauf} and which is related to $f$ in
a natural way. 
We obtain this map 
$h:(Z,D^+)\rightarrow (\Delta,0)$
by a construction that is described in \cite{Steenbrink-Oslo}.
We refer to this process as {\it semistable reduction}.

Then we construct orbits of tangent vectors in $T_0\Delta$ and in
$T_{p_0} D_0$ under the action of certain roots of unity.
All these objects are defined in such a way 
that any right-equivalence between two plane curve singularities
transforms these data into each other. That is, what the word {\it natural}
refers to in this context. 

Finally, we will define {\it the mixed Hodge structure 
on the nearby fundamental 
group of a plane curve singularity} and will finish with an
example.

\subsection{Semistable Reduction} \label{Semistable}

\label{Lol}
Consider the plane curve singularity 
$f=f_0\cdots f_{r-1}$
with $r$ branches and let $\pi:(\tilde{X},E^+)\rightarrow (\cpxm^2,0)$ be 
its embedded
(minimal) resolution. Let $E$ denote the divisor $(f\circ \pi)^{-1}(0)$ 
with components $E=\bigcup_{i\ge 0} E_i$, where $E_i$ has multiplicity
$e_i$. Let $E_0\cup\cdots\cup E_{r-1}$ be 
{\it the strict preimage} of $f^{-1}(0)$ 
in $\tilde{X}$; that 
is the closure of $(f\circ \pi)^{-1}(0)\setminus \pi^{-1}(0)$ in $\tilde{X}$.
Assume that the component $E_i$ corresponds to the factor $f_i$ in $f$.
Define $E^+$, {\it the exceptional divisor}, to be 
$\bigcup_{i>r-1} E_i$. 
The $E_i$'s for $i>r-1$ are all rational curves.
Call $q_0,\ldots,q_{r-1}$ the intersection points 
$E_i\cap E^+$ for $i=0,\ldots,r-1$. 

Let $d=\lcm_{i} e_i$ and let $\theta:\Delta\rightarrow \cpxm$ be the
map $t\mapsto t^d$. Now consider the normalization of the fibered product
of $f\circ \pi$ and $\theta$: 
\[
Z':=\widetilde{\tilde{X}\times_\cpxm\Delta}.
\]
There are projections $h':Z'\rightarrow \Delta'$ and
$P':Z'\rightarrow \tilde{X}$. Denote 
the closure of $(f\circ \pi\circ P')^{-1}(0)
\setminus (\pi\circ P')^{-1}(0)$ in $Z'$ by 
$D'_0\cup\cdots\cup D'_{r-1}$. 
It is not difficult to see that due to 
the fact that $e_0=\cdots =e_{r-1}=1$ the projection $P'$ maps 
$D'_0\cup\cdots\cup D'_{r-1}$ to $E_0\cup\cdots\cup E_{r-1}$. 
Assume the notation to be such that $D'_i$ is mapped to $E_i$.

\begin{proposition}  \label{Trees}
$Z'$ has possibly cyclic quotient singularities. More precisely, for
each singular point $p\in Z'$ there exists a natural number $n_p>0$ 
and an isomorphism
\[
\varphi_{p}: {\mathcal O}_{Z',p} \xrightarrow{\cong} 
\frac{\cpxm\{\alpha,\beta,\eta\}}{(\alpha\beta-\eta^{n_p})}
\]
that maps $h'$ to $\eta$.
\end{proposition}
\noindent
For a proof of this proposition we refer to \cite{Steenbrink-Oslo} 
or \cite{Doktorarbeit}.
 
Now we resolve these singularities of $Z'$ according to a well-known 
procedure (cf.~for instance \cite{Knorrer}). For a singularity of the type: $\alpha\beta-\eta^n$
we replace the singular point by a chain of $(n-1)$
copies of $\prom^1$ in the following way.

Consider $n$ copies of $\cpxm^2$ with coordinates $(u_i,v_i)$ for
$i=1,\ldots,n$. These spaces are glued together by the following
glueing maps:
\[
\begin{matrix}
\varphi_i: & \cpxm^2\setminus\{v_i=0\} & \rightarrow & 
\cpxm^2\setminus\{u_{i+1}=0\}  \\
 & (u_i,v_i) & \mapsto & (\frac{1}{v_i},\; u_iv_i^2)=
\left(u_{i+1}, v_{i+1}\right).
\end{matrix}
\]
Call the result $M$ and observe that $M$ is smooth. 
Now the resolution of the variety $\left\{\alpha\beta-\eta^n =0\right\}$ 
is given
by the map $M\rightarrow \left\{\alpha\beta-\eta^n =0\right\}$,
that is defined on the coordinates $(u_i,v_i)$ by
$(u_i,v_i) \mapsto \left(u_i^iv_i^{i-1};\; u_i^{n-i}v_i^{n+1-i};\;
u_i v_i \right)=\left(\alpha,\beta, \eta\right)$.
Notice that this map is well-defined. In this way we can resolve all
the singularities of $Z'$ and find a map $h:(Z,D^+)\rightarrow (\Delta,0)$ 
with the desired properties. 
We have a commutative diagram of space germs
\begin{equation} \label{Schwangerschaft}
\begin{CD}
(Z,D^+) @>P>> (\tilde{X},{D'}^+) \\
@VhVV @VVf\circ \pi V \\
(\Delta,0) @>>> (\cpxm,0) \\
t & \longmapsto & t^d \;,
\end{CD}
\end{equation}
such that $P$ induces an isomorphism over $\Delta^*$ between the
spaces $Z\setminus h^{-1}(0)$
and $(\cpxm^2\setminus f^{-1}(0))\times_{\cpxm^*} \Delta^*$.
With the notation of \ref{Setzauf} we are ready to prove:
\begin{proposition} \label{Reduction}
Let 
$f,\,\tilde{f}:\left(\cpxm^2,0\right)\rightarrow\left(\cpxm,0\right)$
be two right-equivalent plane curve singularities, i.~e.~there
is a biholomorphism $\phi:(\cpxm^2,0)\rightarrow (\cpxm^2,0)$
such that $\tilde{f}\circ\phi=f$. 
Let then $h:(Z,D^+)\rightarrow (\Delta,0)$ and 
$\tilde{h}:(\tilde{Z},\tilde{D}^+)\rightarrow (\Delta,0)$ be the maps 
constructed by the procdure above for $f$ and $\tilde{f}$ respectively.
Then $\phi$ lifts to a biholomorphism 
${\Phi}:(Z,D^+)\rightarrow (\tilde{Z},\tilde{D}^+)$
such that holds $\tilde{h}\circ \Phi  = h$.
\end{proposition}

\noindent
{\bf Proof:} 
The assertion follows from the fact that both copies of $(\cpxm^2,0)$
have undergone the same blowing-up procedures, which of course
do not depend on coordinates.  \qed

\subsection{Canonical Tangent Vectors}

In this subsection we first recall some well-known facts about the 
map $t\mapsto t^d$. These observations will be used to
construct an orbit under the action 
of the $d$-th roots of unity on $(T_0\Delta)^*$.

If we distinguish one branch $f_0$ of the singularity $f$ then
a careful study of the Puiseux expansion of $f_0$ 
will show subsequently that there
is a natural orbit under the action 
of certain roots of unity on $(T_{p_0}D_0)^*$, which we 
call {\it the monstrance} of $(f,f_0)$. 

\subsubsection{Pulling Back Tangent Vectors}

Let $\Delta$ and $\Delta'$ be two disks in $\cpxm$ around $0$. Always,
when we are given a holomorphic map 
$\varphi:\Delta\rightarrow \Delta'$ of multiplicity $m$,
then the differential of $\varphi$ in $0$ is zero, when $m>1$. However, 
a tangent vector ${\vec{v}\;}'\in (T_0\Delta')^*$ defines in a natural way
$m$ tangent vectors $\zeta_m^\nu \vec{v}\in (T_0\Delta)^*$ for 
$\nu=0,\ldots,m-1$. 

Let us refer to the set $\left\{\zeta_m^\nu \vec{v}\in (T_0\Delta)^*\big|
\;\nu=0,\ldots,m-1\right\}$ as to {\it the inverse star of ${\vec{v}\;}'$ under
$\varphi$}. We will give two descriptions of this phenomenon: 
a concrete one and an abstract one. 

The concrete description is as follows. Choose a coordinate $z'$ on 
$(\Delta',0)$ such that ${\vec{v}\;}'=\frac{\partial}{\partial z'}$. Then
there is up to a choice of an $m$-th root of unity only one coordinate
$z$ on $\Delta$ such that: $z'\circ\varphi=z^m$. 
Define the $m$ tangent vectors by
$\zeta_m^\nu \frac{\partial}{\partial z}\in (T_0\Delta)^*$, where 
$\nu=0,\ldots,m-1$. Notice that this definition does not depend
on the choice of $z'$.

The abstract description is the following. Let $\maxid$ resp.~$\maxid'$
be the maximal ideal in ${\mathcal O}_{\Delta,0}$ resp.~${\mathcal O}_{\Delta',0}$
and let $\maxid\big/\maxid^2$ resp.~$\maxid'\big/{\maxid'}^2$ be the
corresponding Zariski cotangent spaces. Observe that by the duality
\[
\begin{matrix}
T_0\Delta\otimes \left(\maxid\big/\maxid^2\right) & \rightarrow & \cpxm \\
(X\otimes g) & \mapsto & X(g)(0)
\end{matrix}
\]
a choice of a tangent vector in $T_0\Delta$ corresponds to the 
choice of an element in $\maxid\big/\maxid^2$. Now $\varphi$ induces a map
\[
\begin{matrix}
\varphi^*: & \maxid'\big/{\maxid'}^2 & 
\rightarrow & \maxid^m\big/{\maxid}^{m+1}\\
 & g' & \mapsto & g'\circ\varphi .
\end{matrix}
\]
On the other hand we have the composition of natural maps:
\[
\begin{matrix}
\maxid\big/{\maxid}^{2} & \rightarrow & 
\left(\maxid\big/{\maxid}^{2}\right)^{\otimes m} & \rightarrow &
\maxid^m\big/{\maxid}^{m+1} \\
g & \mapsto & g\otimes\cdots\otimes g & &  \\
& & a_1\otimes\cdots\otimes a_m & \mapsto & a_1\cdots a_m.
\end{matrix}
\]
Each non-zero element in $\maxid^m\big/{\maxid}^{m+1}$ has $m$ inverse images
in $\maxid\big/{\maxid}^{2}$ under this composition of maps. 
This shows that a tanget vector ${\vec{v}\;}'\in (T_0\Delta')^*$ determines
$m$ tangent vectors in $(T_0\Delta)^*$ in a natural way. 

For plane curve singularities 
$f:\left(\cpxm^2,0\right)\rightarrow\left(\cpxm,0\right)$, the
tangent vector $\frac{\partial}{\partial t}\in (T_0 \cpxm)^*$ is
preserved under right-equivalences. Denote by 
\[
{\mathcal S}(\frac{\partial}{\partial t})=\{\zeta_d^\nu \vec{v} \in 
(T_0\Delta)^*|\nu=0,\ldots,d-1\} 
\]
{\it the inverse star} of $\frac{\partial}{\partial t}$ under the map 
$\theta: t\mapsto t^d$.

\subsubsection{The Monstrance of $(f,f_0)$}

Let $f_0$ be a distinguished branch of $f$ with Puiseux pairs 
$(m_1,n_1)$, $\ldots,$ $(m_g,n_g)$.
With the definitions $m:=m_1\cdots m_g$
and $k:=n_1 m_2\cdots m_g$, we will define a natural orbit of $m\cdot k$
tangent
vectors in $(T_{p_0}E_0)^*$ under the induced action of the group
$\mu_{m\cdot k}$ of $m k$-th roots of unity. We will refer
to this orbit as to {\it the monstrance} of $(f,f_0)$ and 
$\frac{\partial}{\partial t}$. We denote
it by ${\mathcal M}(f_0;\frac{\partial}{\partial t})$. 
Its elements will serve as {\it base vectors} like 
the $\vec{w}$ in the discussion in the previous sections.

The idea for the construction is the following. Note that 
the map $P$ in (\ref{Schwangerschaft}) maps $E_0$ biholomorphically
to $D_0$. We know that if we 
have a parametrization $\tau:(\cpxm,0)\rightarrow (E_0,q_0)$ of
the strict preimage, then this defines a tangent vector in $T_{q_0}E_0$
respectively $T_{p_0}D_0$,
namely: $\frac{\partial}{\partial \tau}$. The Puiseux expansion is 
a parametrization $p:(\cpxm,0)\rightarrow \left(f_0^{-1}(0),0\right)$
with the property that it can be lifted to a smooth parametrization 
$\tilde{p}:(\cpxm,0)\rightarrow (E_0,q_0)$. In this way, a Puiseux expansion
yields a tangent vector in $T_{q_0}E_0$. Now for the irreducible 
plane curve singularity $f_0$, we 
can always find coordinates $(x,y)$ of $\cpxm^2$ such that 
\[
f_0(x,y)= 1\cdot x^m + \cdots + 1\cdot y^k +\cdots \;,
\]
where $x^m$ and $y^k$ are the smallest pure powers of $x$ and $y$ 
respectively. Up to first order, there are $m\cdot k$ such coordinate
systems. This requirement on the coordinates with respect to $f_0$ 
determines the Puiseux expansions up to first order. Up to first order,
we get $m\cdot k$ canonical parametrizations 
of $(E_0,q_0)$ and therefore $m\cdot k$ canonical tangent vectors. 

But let us give a geometric meaning to this construction. Since 
for the moment we do not want to fix coordinates on $\cpxm^2$,
we let $(X,x_0):=(\cpxm^2,0)$. 
\begin{definition} \rm
Let $\iota:(L,x_0)\hookrightarrow (X,x_0)$ be a 1-dimensional
smooth subspace of $(X,x_0)$. We will call a map
${\bf P}:(X,x_0)\rightarrow (L,x_0)$
{\it a projection} of $(X,x_0)$ onto $(L,x_0)$, if the induced map
$T_{x_0} L \xrightarrow{D\iota}  T_{x_0} X  \xrightarrow{D{\bf P}} T_{x_0}L$
is the identity. 
\end{definition}
 
Now let $\iota:(R,x_0)\hookrightarrow (X,x_0)$ be the tangent 
to $(f_0^{-1}(0),x_0)\subset (X,x_0)$. We formulate the 
following theorem. 
\begin{theorem} \label{Endspurt}
Let $f=f_0\cdots f_{r-1}:(X,x_0)\rightarrow (\cpxm,0)$ be a plane curve
singularity with one distinguished branch $f_0$. Let  
$(m_1,n_1),\ldots,(m_g,n_g)$ be the Puiseux pairs of $f_0$.
Set $m:=m_1\cdots m_g$ and $k:=n_1 m_2\cdots m_g$.

Let $\iota:(R,x_0)\hookrightarrow (X,x_0)$ be the tangent and
$\pi:(\tilde{X},E^+)\rightarrow (X,x_0)$ be the embedded resolution of $f$.
Then for any projection ${\bf P}:(X,x_0)\rightarrow (R,x_0)$,
the composition $\phi$ of maps
\[
(E_0,q_0)\xrightarrow{\pi_{|(E_0,q_0)}} (X,x_0) \xrightarrow{{\bf P}}
(R,x_0) \xrightarrow{{f_0}_{| (R,x_0)}} (\Delta,0)
\]
has multiplicity $m\cdot k$ and the inverse star of the tangent 
vector $\frac{\partial}{\partial t}$ under $\phi$ does 
not depend on the choice of the projection ${\bf P}$.  
\end{theorem}

\noindent
{\bf Proof:}
Choose coordinates $(x,y)$ on $(X,x_0)$ 
such that the following two conditions are satisfied:
\begin{itemize}
\item
The tangent line $(R,x_0)$ of $(f_0^{-1}(0),x_0)$ is given by $x=0$.
\item
$ f_0(x,y)= 1\cdot x^m + \cdots + 1\cdot y^k +\cdots $, where $x^m$ and $y^k$
are the smallest pure powers of $x$ and $y$ respectively.
\end{itemize}

Then the Puiseux expansion of $f_0$ with respect to $(x,y)$ has the
following form, where $m<k$ ($f_0$ is irreducible):
\begin{align*}
x(\tau) &= a_k \tau^k + a_{k+1} \tau^{k+1} + \cdots \\
y(\tau) &= \tau^m
\end{align*}
and $a_k^m=-1$. Moreover we can find $f_0$ back by:
\[
f(x,y)=\prod_{\nu=0}^{m-1} \left( x- x(\zeta_m^\nu y^{\frac{1}{m}})\right).
\]
$\tau$ can be used as a coordinate on the strict transform $(E,q_0)$. 
We may take $\tilde{y}:=y_{\big| (R,x_0)}$ as coordinate 
on the tangent line $(R,x_0)$. 

\enlargethispage{1cm}
Any projection ${\bf P}:(X,x_0)\rightarrow (R,x_0)$ is a priori given 
by $\tilde{y}\circ{\bf P}(x,y) = y({\bf P}(x,y))= 
Ax+By + {\bf h}(x,y)$ and ${\bf h}(x,y)\in \maxid_{(X,x_0)}^2$. 
The condition $D{\bf P}_{\big| T_{x_0} R}= id$ means here: $B=1$.
Hence, we find for a unit $u(\tau)\in \cpxm\{\tau\}\setminus \maxid$ with
$u(0)=\zeta_m^\nu$:
\[
\begin{split}
\tilde{y}\circ{\bf P}\left(x(\tau),y(\tau)\right) &= 
\tau^m + \tau^{m+1}\left(\cdots\right) \\
&= \tau^m\cdot u(\tau)^m. \\
\text{Therefore} \qquad\qquad & \\
\phi(\tau) &= f_0\left(0,\tilde{y}\circ{\bf P}\left(x(\tau),\;
y(\tau)\right)\right)\\
 &= (-1)^m\prod_{\nu=0}^{m-1} x\left(\zeta_m^\nu \tau\cdot u(\tau) \right)\\
 &= (-1)^m\prod_{\nu=0}^{m-1} \left(a_k(\zeta_m^\nu \tau\cdot u)^k +
a_{k+1}(\zeta_m^\nu \tau\cdot u)^{k+1}+\cdots \right) \\
 &= (-1)^{m+k}a_k^m\tau^{m\cdot k}+ \tau^{m\cdot k+1}\big(\cdots\big).
\end{split}
\]
This shows that $\phi(\tau)$ has multiplicity $m\cdot k$ and the
term of lowest order does not depend on the choice of ${\bf P}$. \qed

\begin{definition} \rm
Define {\it the monstrance of $(f,f_0)$} to be the inverse star of
$\frac{\partial}{\partial t}$ under a map $\phi$ like in Theorem 
\ref{Endspurt} and denote it by 
\[
{\mathcal M}(f_0;\frac{\partial}{\partial t})
=\{\zeta_{mk}^\nu \vec{w} \in 
(T_{q_0}E_0)^*|\;\nu=0,\ldots,mk-1\}. 
\]
\end{definition}

\begin{remark} \rm
$m\cdot k$ does not divide $d$ in general. For instance in the example 
of Section \ref{shakuhachi} is $m=4$, $k=6$ and $d=156$. 
However, if $f=f_0$ and $f_0$ has only one Puiseux pair 
$(m,k)=(m_1,n_1)$, then 
$m\cdot k$ divides $d$ as one easily sees by blowing up twice.
\end{remark}

\subsection{The Invariant}

Here finally we are going to associate to a plane curve 
singularity an invariant of its right-equivalence class. In fact it is an
invariant of the following seemingly coarser equivalence relation. 

\begin{definition} \rm
We say that two hypersurface singularities 
\[
f,\,g:(\cpxm^{n+1},0)\rightarrow (\cpxm,0) 
\]
are {\it right-equivalent to first order}
if there are isomorphisms $\phi:(\cpxm^{n+1},0)\rightarrow(\cpxm^{n+1},0)$
and $\varphi:(\cpxm,0)\rightarrow (\cpxm,0)$ such that 
$g\circ \phi = \varphi\circ f$ and 
$\varphi_*:T_0\cpxm\rightarrow T_0\cpxm$ is the identity.
\end{definition}

The following theorem\footnote{This theorem was pointed out to us by
Claus Hertling.} 
might be interesting in this context. 
\begin{theorem} \label{Geburtszange}
Two plane curve singularities are right-equivalent if they are  
right-equivalent to first order.
\end{theorem}

\noindent
{\bf Proof:}
Let $f:(\cpxm^{2},0)\rightarrow (\cpxm,0)$ be a plane curve singularity
and let $h(t)=t+a_2 t^2+a_3 t^3 + \cdots$ be in $\cpxm\{t\}$. We have to
show that $f$ and $h\circ f$ are right-equivalent. 
Let $\maxid$ be the maximal ideal in $\cpxm\{x,y\}$.
According to a theorem of Brian\c{c}on-Skoda \cite{Briancon-Skoda}
(see also \cite{Lipman-Teissier})
$f^2\in \maxid \cpxm\{x,y\}$. Using this result one
can show by a standard construction
(see for instance \cite{Lamotke} p.~189, V\S 4, Prop.~4) that all $f_t$ in 
$f_t(x,y):= f(x,y) + t\left( a_2 f^2(x,y) +a_3 f^3(x,y)+\cdots \right)$
are mutually right-equivalent. \qed 

\newcommand{\kaen}{\text{ ${\mathfrak P}$}}
\newcommand{\ddt}{{\scriptstyle \frac{\partial}{\partial t}}}

\begin{definition} \rm
Let $f=f_0\cdots f_{r-1}:
(\cpxm^2,0)\rightarrow (\cpxm,0)$ be a plane 
curve singularity. For all $s\ge 1$ we define
\begin{equation} \label{bunter Hund}
\kaen^s(f,\ddt)
:= \bigoplus_{i=0}^{r-1}\left\{\bigoplus_{\vec{v}\in{\mathcal S}(\ddt)}\quad
\bigoplus_{\vec{w}_i\in{\mathcal M}(f_i;\ddt)}
\left(J_{\vec{v},\vec{w}_i}\big/J_{\vec{v},\vec{w}_i}^{s+1}\right)^*\right\}
\end{equation}
and we call the direct limit
\[
\kaen(f,\ddt):=
\lim_{\xrightarrow[s]{}}\;\kaen^s(f,\ddt)
\]
{\it the mixed Hodge structure on the nearby fundamental group of
the plane curve singularity $f$}. 
\end{definition}

Note that $\mu_d$ acts on $\kaen(f,\ddt)$.
More precisely, $\zeta_d^\nu=e^{\frac{2\pi i}{d}\nu}\in \mu_d$ 
acts as the map: 
\begin{equation} \label{Pudel}
\left(\zeta_d^\nu\right)^{-N}=e^{2\pi i \nu \frac{N}{d}}: 
\left(J_{\vec{v},\vec{w}_i}\big/J_{\vec{v},\vec{w}_i}^{s+1}\right)^*_\itgm
\rightarrow 
\left(J_{\zeta_d^\nu\vec{v},\vec{w}_i}
\big/J_{\zeta_d^\nu\vec{v},\vec{w}_i}^{s+1}
\right)^*_\itgm.
\end{equation}

Let $(m^{(i)}_1,n^{(i)}_1),\ldots,(m^{(i)}_{g_i},n^{(i)}_{g_i})$ be the 
Puiseux pairs of the branch $f_i$ and define 
$m^{(i)}:= m^{(i)}_1\cdots m^{(i)}_{g_i}$ and 
$k^{(i)}:= n^{(i)}_1m^{(i)}_2\cdots m^{(i)}_{g_i}$. Then each of the groups
$\mu_{m^{(i)}k^{(i)}}$ acts on $\kaen(f,\ddt)$. Let $M^{(i)}$ be the map
$M$ of \ref{Gr"unspan} associated to the distinguished branch $f_i$.
Then $\zeta_{m^{(i)}k^{(i)}}^\nu$ induces the map 
\begin{equation} \label{Dackel}
\left(\zeta_{m^{(i)}k^{(i)}}^\nu\right)^{M^{(i)}}: 
\big(J_{\vec{v},\vec{w}_i}\big/J_{\vec{v},\vec{w}_i}^{s+1}\big)^*_\itgm
\rightarrow 
\big(J_{\zeta_d^\nu\vec{v},
\zeta_{m^{(i)}k^{(i)}}^\nu\vec{w}_i}\big/
J_{\zeta_d^\nu\vec{v},\zeta_{m^{(i)}k^{(i)}}^\nu\vec{w}_i}^{s+1}
\big)^*_\itgm.
\end{equation}

By Theorem \ref{Gr"unspan} all the maps in (\ref{Pudel}) and 
(\ref{Dackel}) are isomorphisms of MHSs over $\itgm$.
Let ${\mathcal T}$ be the monodromy of the local system 
\[
\kaen(f,\ddt)_\itgm=
\bigoplus_{i=0}^{r-1}\left\{\bigoplus_{\vec{v}\in{\mathcal S}(\ddt)}\quad
\bigoplus_{\vec{w}_i\in{\mathcal M}(f_i;\ddt)}
\left(J_{\vec{v},\vec{w}_i}
\big/J_{\vec{v},\vec{w}_i}^{s+1}\right)^*_\itgm\right\}.
\]
Define $L:=(\frac{1}{d}N-\sum_{i=0}^{r-1}\frac{1}{m^{(i)}k^{(i)}}M^{(i)}):
H^0\bar{B}(A^\bullet,a)
\rightarrow H^0\bar{B}(A^\bullet,a)$ and moreover
\[
{\mathcal L}:=\bigoplus_{i,\vec{v},\vec{w}_i} L:
\bigoplus_{i,\vec{v},\vec{w}_i}
 H^0\bar{B}(A^\bullet,a)\rightarrow
\bigoplus_{i,\zeta_d\vec{v},\zeta_{m^{(i)} k^{(i)}}\vec{w}_i}
 H^0\bar{B}(A^\bullet,a).
\]
The map ${\mathcal L}$ sends the summand with index $(i,\vec{v},\vec{w}_i)$ to 
the summand corresponding to 
$(i,\zeta_d\vec{v},\zeta_{m^{(i)}k^{(i)}}\vec{w}_i)$.
By virtue of Theorem \ref{Gr"unspan} we find for $\vartheta\in\cpxm^*$:
\[
{\mathcal T}=e^{-2\pi i {\mathcal L}} \quad\text{ and }\quad
\kaen^s(f,\vartheta \ddt)_\itgm 
= \vartheta^{-{\mathcal L}}
\kaen^s(f,\ddt)_\itgm,
\]
from where we may conclude that the family of MHSs
\[
\left\{\kaen(f,\ddt)
\right\}_{\ddt\in T_0\cpxm}
\]
is a nilpotent orbit of MHS. 
Note that on the components of the direct sum 
\eqref{bunter Hund} the monodromy ${\mathcal T}$ is given by
\[
e^{-2\pi i L}
: \left(J_{\vec{v},\vec{w}_i}\big/J_{\vec{v},\vec{w}_i}^{s+1}\right)^*_\itgm
\rightarrow
\left(J_{\zeta_d\vec{v},\zeta_{m^{(i)}k^{(i)}}\vec{w}_i}\big/
J_{\zeta_d\vec{v},\zeta_{m^{(i)}k^{(i)}}\vec{w}_i}^{s+1}\right)^*_\itgm.
\]
We consider both lattices as sub-lattices of $H^0\bar{B}(A^\bullet,a)$. 

\subsection{An Example} \label{shakuhachi}

\afterpage{\clearpage
\begin{figure}[H]
\centering\epsfig{file=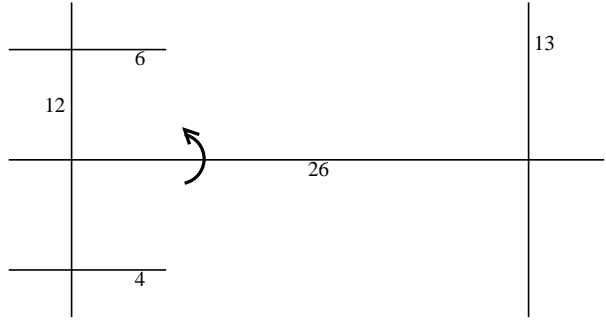,width=8cm,clip=}
\caption{The resolution of $f_\lambda(x,y)=
\left(y^2-x^3\right)^2-\left(4\lambda x^5 y
+\lambda^2 x^7\right)$. } \label{Aufloesung}
\end{figure}}

The only purpose of the example we give in this subsection is to indicate that 
{\it the mixed Hodge structure on the nearby fundamental group of
a plane curve singularity} detects moduli that are hidden
for the MHS on the vanishing cohomology. We focus especially on effects
coming from the monodromy on $\vec{J}\big/{\vec{J}\,}^4$. But in general,
we also expect interesting information from higher periods -- even 
on $\vec{J}\big/{\vec{J}\,}^3$. A systematic investigation of the invariant
on classes of singularities still remains to be done. 

Here we investigate the example \eqref{Liszt} that was mentioned 
in the introduction:
\begin{equation*}
f_\lambda(x,y)=\left(y^2-x^3\right)^2-\left(4\lambda x^5 y+
\lambda^2 x^7\right),\quad \lambda \neq 0.
\end{equation*}

All $f_\lambda$ have Milnor number 16 and the resolution of $f_\lambda$ 
is shown in Figure \ref{Aufloesung}.
From this resolution and from the Milnor number it is not difficult 
(cf.\cite{Brieskorn}) to derive that the divisor with normal crossings
$D\subset Z$ has the shape that is illustrated in Figure \ref{Weisswal}.
Moreover we see $d=156$. Notice also that $m=4$ and $k=6$.

Observe that the automorphism of $(\cpxm^2,0)$ which is given by
\begin{equation}
\phi: \begin{pmatrix} \tilde{x} \\ \tilde{y} \end{pmatrix}=
\begin{pmatrix} \vartheta^2 x \\ \vartheta^3 {y} \end{pmatrix}
\end{equation}
has the property that: $f_\lambda(\tilde{x}, \tilde{y}) = \vartheta^{12}
f_{\vartheta\lambda}(x,y)$. Or likewise, $\phi$ yields a 
commutative diagram:
\[
\begin{CD}
(\cpxm^2,0) @>\phi>> (\cpxm^2,0) \\
@V{f_{\vartheta\lambda}}
VV @VVf_{\lambda}V \\
(\cpxm,0) @>{\cdot\vartheta^{12}}>> (\cpxm,0) .
\end{CD}
\]
Hence, the MHS on the nearby fundamental group of the pair 
$({f_{\vartheta\lambda}},\ddt)$ is isomorphic
to the one of $({f_{\lambda}},\vartheta^{12}\ddt)$, i.~e.
$\kaen({f_{\vartheta\lambda}},\ddt)
\cong\kaen ({f_{\lambda}},\vartheta^{12}\ddt)$.

Therefore, we fix one $f_\lambda$ and consider how 
$\kaen ({f_{\lambda}},\vartheta^{12}\ddt)$ varies, 
when we change $\vartheta^{12}=:\theta$.
We know that for any $s\ge 1$, the variation of MHS
$\left\{\kaen^s({f_{\lambda}},\theta\ddt)\right\}_{
\theta\in\cpxm^*}$
is the nilpotent orbit of MHS associated to the pair 
$\left(\kaen^s({f_{\lambda}},\theta\ddt),
{\mathcal L}\right)$, where ${\mathcal L}=\bigoplus_{\vec{v},\vec{w}} L$
is the nilpotent endomorphism given by $L=\frac{1}{156} N - \frac{1}{24} M$.

Denote by $h_\lambda:(Z_\lambda,D_\lambda^+)\rightarrow (\Delta,0)$
the map that is obtained by semistable reduction from $f_{\lambda}$. 
By Theorem \ref{Reduction} we have a commutative diagram 
\[
\begin{CD}
(Z_\lambda,D_\lambda^+) @>>> (\cpxm,0) \\
@V{h_\lambda}VV @VV{f_\lambda}V \\
(\Delta,0) @>>> (\cpxm,0) \\
t & \longmapsto & t^{156}.
\end{CD}
\]

\afterpage{\clearpage
\begin{figure}[H]
\epsfig{file=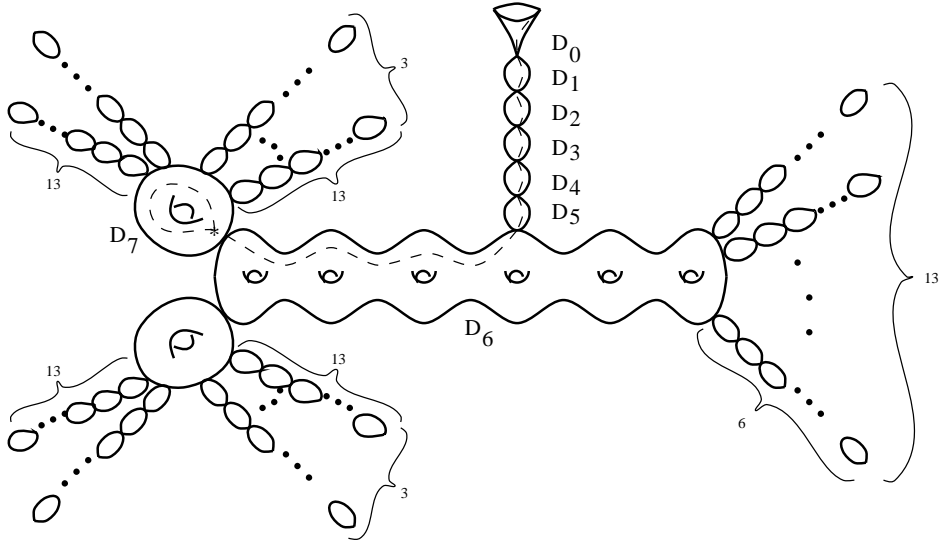,clip=}
\caption{The divisor $D$ for
$f_\lambda(x,y)=\left(y^2-x^3\right)^2-\left(4\lambda x^5 y
+\lambda^2 x^7\right)$.} \label{Weisswal}
\end{figure}}

To $h_\lambda$ we associate the dga $A^\bullet$ and note by 
Corollary \ref{Gr"unschnabel} on page \pageref{Gr"unschnabel} that
$H^1(A^\bullet)$ is a pure Hodge structure of weight 1. Since 
$(\vec{J}\big/{\vec{J}\,}^3)^*$ has weights 1 and 2 and since $N$ and $M$
both lower the weights by 2, we see that $N$ and $M$ are 
both zero on $(\vec{J}\big/{\vec{J}\,}^3)^*$. This shows that the
nilpotent orbit of MHS 
\[
\left\{ \kaen^2(f_\lambda,\theta\ddt) 
\right\}_{\theta\in\cpxm^*}
\]
is constant. Our goal now is to show:
\begin{proposition}
The nilpotent orbit of MHS 
\[
\left\{ \kaen^3(f_\lambda,\theta\ddt) 
\right\}_{\theta\in\cpxm^*}
\]
is \underline{not} constant. 
\end{proposition}

\noindent
{\bf Proof:}
In order to prove this proposition we just have to show that $L\neq 0$.
Choose $\vec{v}\in{\mathcal S} (\ddt)$ and
$\vec{w}\in{\mathcal M} (\ddt)$. On the
elliptic curve $D_7$, choose a point $*\in D_7$ and a closed path
$\alpha$ based at $*$. Let $\omega^{(7)}$ be a holomorphic 1-form
on $D_7$ such that 
\[
\int_\alpha \omega^{(7)} =1.
\]
There exists a $\mu^{(7)}\in E^1(D_7\log p_{67})$ such that 
$\omega^{(7)}\wedge\bar{\omega}^{(7)} + d\mu^{(7)} = 0$
and $\rho:=\Res_{p_{67}} \mu^{(7)} \neq 0$. Now let $\eta_1^{(7)},\ldots,
\eta_6^{(7)}$ be meromorphic forms on $D_1,\ldots,D_6$ respectively
and $\eta_0^{(7)}=\rho_{01}\frac{dp}{p}\in{\bigwedge}^1 (\frac{dp}{p})$
such that:
\[
\rho_{i,i+1} = \Res_{p_{i,i+1}} \eta_i^{(7)}= -\Res_{p_{i-1,i}} \eta_i^{(7)}
= -\rho.
\]

Now consider the Chen-closed element in $(A^1\otimes A^1)\oplus A^1$
that is given by:
\[
\Omega=\omega^{(7)}\otimes\omega^{(7)}\otimes\bar{\omega}^{(7)}
+\omega^{(7)}\otimes\left(\mu^{(7)} +\eta_6^{(7)}+\cdots+\eta_0^{(7)}
+\rho_{67}\Theta+\cdots+\rho_{01}\Theta \right).
\]
And note that $N(\Omega)=\omega^{(7)}\otimes
\left(\rho_{67}d\xi+\cdots+\rho_{01}d\xi \right)$. There is an 
$h\in A^0$ such that $h(p_{01})=0$ and $dh=\rho_{67}d\xi+\cdots+\rho_{01}d\xi $.
It is easy to see that for this $h$ holds
\[
[N(\Omega)]= \left(\omega^{(7)}\big|
\left(\rho_{67}d\xi+\cdots+\rho_{01}d\xi \right)\right) = (-h\omega^{(7)})
=(7\rho \omega^{(7)})\in H^0\bar{B}_3(A^\bullet,a).
\]
Similarly we see 
$[M(\Omega)]=(\rho\;\omega^{(7)})\in H^0\bar{B}_3(A^\bullet,a)$.
From this we derive:
\[
[L(\Omega)] = \frac{1}{156}[N(\Omega)] - \frac{1}{24}[M(\Omega)] 
  = (-\frac{32}{156} \rho\;\omega^{(7)}) \neq 0
\]
\qed

\begin{remark} \rm
Since $\kaen^2(f_\lambda,\theta\ddt)
=\kaen^2(f_\lambda,\ddt)$ for all $\theta \in \cpxm^*$,
we may consider the map 
\begin{equation} \label{Theodorina}
\begin{matrix}
\cpxm^* & \rightarrow & \Extmhs\left( \bigoplus_{\vec{v},\vec{w}} 
H^1(A^\bullet);\; \kaen^2(f_\lambda,\ddt)\right) \\
\theta & \mapsto & \left\{ 0 \to 
\kaen^2(f_\lambda,\ddt)\to
\kaen^3(f_\lambda,\theta\ddt)\to
\bigoplus_{\vec{v},\vec{w}} H^1(A^\bullet)^{\otimes 3}\right\}.
\end{matrix}
\end{equation}
The fact that $[\omega^{(7)}]\otimes[\omega^{(7)}]\otimes[\bar{\omega}^{(7)}]
\in F^2(H^1(A^\bullet)^{\otimes 3})$ and $L(\Omega)\not\in F^2H^1(A^\bullet)=0$
can be used to prove that even the map \eqref{Theodorina} is
not constant. 
\end{remark}


\end{document}